\documentclass[11pt,reqno]{amsart}
\usepackage{indentfirst}
\usepackage{graphics,amsmath,amssymb,amsthm,mathrsfs,bm}
\usepackage{enumitem}
\usepackage{amsfonts}
\usepackage{xcolor}
\usepackage[colorlinks=true, linkcolor=red, urlcolor=blue, citecolor=blue]{hyperref}
\usepackage[leqno]{amsmath}
\usepackage{amsmath}
\usepackage{amssymb}
\usepackage{multicol}
\usepackage{stmaryrd}
\usepackage{setspace}
\usepackage{bbm}
\usepackage{cite}
\usepackage{epsfig}
\usepackage{color}
\usepackage{graphics}
\usepackage{graphicx}
\usepackage{epstopdf}
\usepackage{multicol,graphics}
\usepackage{mathtools}
\usepackage{indentfirst}
\usepackage{float}
\usepackage{hyperref}
\hypersetup{pdfstartview=Fit, pdfborder={0 0 0},breaklinks=true}

\newtheorem{theorem}{Theorem}[section]
\newtheorem{lemma}[theorem]{Lemma}

\newtheorem{definition}[theorem]{Definition}
\newtheorem{remark}[theorem]{Remark}
\newtheorem{proposition}[theorem]{Proposition}
\numberwithin{equation}{section}

\allowdisplaybreaks[4]
\makeatletter
\def\@textbottom{\vskip \z@ \@plus 30pt}
\let\@texttop\relax
\makeatother
\begin{document}
\title[Principal spectral theory of time-periodic advective system]{Principal spectral theory of time-periodic nonlocal advective cooperative systems and its applications}

\author[Wu, Fang, Li and Ni]{Hao Wu$^{1}$, Xiangdong Fang$^{2}$, Wan-Tong Li$^{1,*}$ and Wenjie Ni$^{3}$}
\thanks{\hspace{-.6cm}
$^1$School of Mathematics and Statistics, Lanzhou University, Lanzhou, Gansu 730000, China.\\
$^2$School of Mathematical Sciences, Dalian University of Technology, Dalian 116024, China.\\
$^3$School of Science and Technology, University of New England, Armidale, NSW 2351, Australia.\\
$^*${\sf Corresponding author} (wtli@lzu.edu.cn)}
\date{\today}

\begin{abstract}
We develop a new spectral theory for time-periodic nonlocal cooperative systems with advection on bounded intervals, allowing asymmetric dispersal kernels. By transforming the problem into a nonlocal moving-boundary problem, we show that the associated Poincar\'e map is strongly positive and power compact, which yields existence and uniqueness of the principal eigenvalue and several variational characterizations. For spatially homogeneous coefficients, we establish the limit of the principal eigenvalue as the interval expands to $\mathbb R$ and identify it as the minimum of the principal Floquet exponents of the associated exponentially weighted periodic systems. To the best of our knowledge, this spectral theory is new even in \textbf{the scalar case}. We then apply these results to nonlinear systems. On bounded intervals, the sign of the principal eigenvalue gives a sharp dichotomy: solutions either converge to a unique positive time-periodic solution or vanish as time tends to infinity. We also prove convergence of positive periodic solutions as the interval expands. For the Cauchy problem on the whole real line, we characterize the leftward and rightward spreading speeds, which may differ for asymmetric dispersal kernels, and establish necessary and sufficient conditions for accelerated spreading in either direction. In particular, under irreducible coupling, accelerated spreading in one component induces accelerated spreading in all components in the same direction.

		

\end{abstract}

\maketitle

\hypersetup{linkcolor=black}  
\tableofcontents
\hypersetup{linkcolor=red}

\section{Introduction}

\subsection{Principal eigenvalue problem}

The first objective of this paper is to establish a principal eigenvalue theory for the following time-periodic 	nonlocal cooperative system with advection:
\begin{equation}\label{P-draft}
\left\{
\begin{aligned}
& \partial_{t} \bm{\varphi} = c(t)\partial_{x}\bm\varphi + D(x,t)\mathcal{K}_{\Omega}[\bm\varphi](x,t) + A(x,t)\bm\varphi - \lambda\bm\varphi, &&(x,t)\in\Omega\times \mathbb{R},\\
&\bm\varphi(x,t+T)= \bm\varphi(x,t),&&(x,t)\in\Omega\times \mathbb{R},\\
&\bm\varphi(l_{1},t)=\mathbf0,  && \text{ if } \sup_{t} c(t)<0, \\
&\bm\varphi(l_{2},t)=\mathbf0,  &&\text{ if } \inf_{t} c(t)>0.
\end{aligned}
\right.
\end{equation}
Here, $\Omega=(l_1,l_2)\subset\mathbb R$ is a bounded interval, $\bm\varphi=(\varphi_1,\ldots,\varphi_m)^T$, and $c(t)$ is the advection coefficient. The matrices $D(x,t)=\operatorname{diag}(d_i(x,t))$ and $A(x,t)=(a_{ij}(x,t))_{m\times m}$ contain the dispersal rates and the coefficients of the linear reaction system, respectively. The nonlocal dispersal operator is given by $\mathcal K_\Omega=\operatorname{diag} (\mathcal K^1_\Omega,\ldots,\mathcal K^m_\Omega)$, where 
\begin{equation*}
\mathcal{K}^{i}_{\Omega}[\varphi_{i}](x,t) =\int_{\Omega} k_{i}(x-y,t)\varphi_{i}(y,t)dy- \varphi_{i}(x,t), \quad 1\le i \le m.
\end{equation*}
The precise assumptions on $c$, $D$, $A$, and $k_i$ are given below.

To put \eqref{P-draft} in context, first consider the case without advection, $c\equiv0$. Then \eqref{P-draft} reduces to
\begin{equation}\label{P-non-draft}
\left\{
\begin{aligned}
& \partial_{t} \bm{\varphi} = D(x,t)\mathcal{K}_{\Omega}[\bm\varphi](x,t) + A(x,t)\bm\varphi - \lambda\bm\varphi, &&(x,t)\in\Omega\times \mathbb{R},\\
&\bm\varphi(x,t+T)= \bm\varphi(x,t),&&(x,t)\in\Omega\times \mathbb{R}.
\end{aligned}
\right.
\end{equation}

Principal eigenvalues for nonlocal dispersal problems without advection have been studied extensively, including their existence, uniqueness, variational characterizations, and asymptotic behavior (see, e.g., \cite{BCV1,SLLW,SLLY,SWZ,SX,SX1,ShenVo,NV2,SZ,SLW,WZ,WFLS,WLSV,LZZ,LZZ1,LinWang,HMMV03,ZL1}). In particular, several approaches have been developed to establish the existence of principal eigenvalues. One uses a generalized Krein–Rutman theorem due to Nussbaum \cite{Nussbaum} and Edmunds, Potter, and Stuart \cite{Edmunds}; see Coville \cite{Co10}, Liang et al. \cite{LZZ,LZZ1}, and Lin and Wang \cite{LinWang}. Another approach is based on the ``touching lemma'' developed by Li et al. \cite{LCW}; see also Su et al. \cite{SLLW,SWZ}. An alternative way is to use an idea of perturbation theory for positive semigroups developed by B\"urger \cite{Burger}; see Bao and Shen \cite{BSS17} and Rawal and Shen \cite{RS}. The theory of resolvent-positive operators developed by Thieme \cite{Th1,Th2} has also been used to establish existence; see Kang and Ruan \cite{KR,KR4}, Feng et al. \cite{Feng2025, Feng}.

For nonlocal dispersal problems with advection, relatively few results are currently available on principal eigenvalues. In the autonomous scalar setting, Li et al. \cite{LCW} established the existence of principal eigenvalues on bounded intervals. Coville and Hamel \cite{CH} further studied generalized principal eigenvalues on both bounded and unbounded intervals. Du et al. \cite{DFN} derived an explicit formula for the principal eigenvalue limit as the interval expands to $\mathbb R$, allowing asymmetric dispersal kernels. Tong and Zhou \cite{TongZhou} obtained bounds for such limits in spatially heterogeneous settings and a variational formula under almost periodic assumptions. Vo \cite{Vo} studied the asymptotic behavior of generalized principal eigenvalues as the dispersal range tends to zero or infinity. Ducrot et al. \cite{DGL} showed that multiple principal eigenpairs may occur in spatially periodic scalar problems when the advection changes sign. Bentout and Vo \cite{BV1} studied the principal eigenvalue problem for an autonomous two-component cooperative system. In the non-autonomous setting, Cloez et al. \cite{CAG} recently established the existence of a principal Floquet eigenvalue for scalar time-periodic equations on $\mathbb R^d$. Recently, He \cite{He2026} considered cooperative systems on the circle.

However, the following two fundamental differences make it difficult to analyze the eigenvalue problem~\eqref{P-draft} using existing methods.

The first difficulty arises from the inflow boundary condition in \eqref{P-draft}, which is absent in nonlocal dispersal problems without advection and in advective problems on the whole space or in spatially periodic settings.  For \eqref{P-draft}, the eigenfunction must vanish at the inflow endpoint, while its nonlocal integral term is generally nonzero there. Consequently, the eigenvalue equation imposes a compatibility relation between the spatial derivative of the eigenfunction and the nonlocal integral. This makes the choice of function space delicate, as it must accommodate the compatibility condition while having a positive cone with nonempty interior. Moreover, precise boundary estimates are needed to compare positive functions near the inflow endpoint.

The second difficulty concerns the regularizing effect of advection in two different settings: autonomous problems and time-periodic problems. For an autonomous problem, the differential part of the eigenvalue equation is the first-order ordinary differential operator $-c\frac{d}{dx}$. Thus, the eigenvalue equation can be viewed as an ordinary differential equation (ODE) in $x$, with the nonlocal integral term treated as a perturbation. This ODE structure provides additional spatial regularity, which facilitates the analysis of the principal eigenvalue. In contrast, in the time-periodic setting, the differential part is the first-order hyperbolic operator $\partial_t-c(t)\partial_x$,
which provides regularity only along characteristic curves, rather than directly in the spatial variable $x$.
Consequently, the same ODE-based argument cannot be used to obtain the spatial regularity required for the spectral analysis.

To overcome these difficulties, the spatial variable is shifted by $\alpha(t)=\int_0^t c(s)\,ds$. This change of variables removes the advection term and transforms the problem on the fixed interval $(l_1,l_2)$ into a nonlocal problem on the moving interval $(l_1+\alpha(t),l_2+\alpha(t))$. More precisely, the transformation $\bm{\phi}(x,t)=\bm{\varphi}(x-\alpha(t),t)$ takes \eqref{P-draft} into
\begin{equation}\label{Free-EP}
\left\{
\begin{aligned}
&\partial_t\bm{\phi}
      =D^*(x,t)\mathcal J_{(g(t),h(t))}[\bm{\phi}](x,t)
        +A^*(x,t)\bm{\phi}-\lambda\bm{\phi},
&&g(t)<x<h(t),\, t\in\mathbb R,\\
&\bm{\phi}(x+\alpha(T),t+T)=\bm{\phi}(x,t),
&&g(t)<x<h(t),\, t\in\mathbb R,\\
&\bm{\phi}(g(t),t)=\mathbf0,
&& \text{if }\sup_t c(t)<0,\\
&\bm{\phi}(h(t),t)=\mathbf0,
&& \text{if }\inf_t c(t)>0.
\end{aligned}
\right.
\end{equation}
Here $g(t)=l_1+\alpha(t)$ and $h(t)=l_2+\alpha(t)$ for $t\in\mathbb R$. For $g(t)\leq x\leq h(t)$ and $t\in\mathbb R$,
$D^*(x,t)=D(x-\alpha(t),t)$, $A^*(x,t)=A(x-\alpha(t),t)$, and 
\begin{equation*}
    \left(\mathcal J_{(g(t),h(t))}[\bm{\phi}]\right)_i(x,t)
    =\int_{g(t)}^{h(t)}k_i(x-y,t)\bm{\phi}_i(y,t)\,dy-\bm{\phi}_i(x,t),
    \qquad 1 \le i\le m.
\end{equation*}
For the transformed problem, a Poincar\'e map is constructed by evolving solutions over one period and shifting them back to the original interval (Remark \ref{Well defined}). The fixed direction of advection allows a sufficiently large power of this map to be compact, as established in Theorem~\ref{thm:poincare-spectral-properties}. Together with strong positivity, this allows us to apply the spectral theory of power-compact positive operators to establish the existence and uniqueness of the principal eigenvalue of \eqref{P-draft}. The principal eigenvalue is then characterized variationally (Theorem~\ref{EX-PEV}). When $D$ and $A$ are independent of $x$, its large-domain limit is expressed using Floquet exponents of exponentially weighted periodic systems (Theorem~\ref{cor:general-kernel-large-domain-limit}).

\subsection{Persistence and spreading in nonlinear cooperative systems}
The eigenvalue problem above arises naturally in the study of spreading speeds, traveling waves, and dynamics in
shifting environments for the following nonlinear nonlocal cooperative system on $\mathbb R$.
\begin{equation}\label{1.4}
	\left\{
	\begin{aligned}
		&\partial_{t} U_{i} = d_{i}(x,t)\left( k_{i}*U_i-U_i\right)(x,t) +f_{i} (x,t,U_1,\cdots, U_m),&& x\in\mathbb{R},t>0,\\
		&U_{i}(x,0)=U_{0,i}(x), &&x\in \mathbb{R},\\
		& i=1,2,\cdots,m,
	\end{aligned}
	\right.
\end{equation}
where $(k_{i}*\varphi)(x,t)=\int_{\mathbb{R}} k_{i}(x-y,t)\varphi(y,t)dy$. An effective approach to these problems is to consider the system in moving coordinates. For a given function $c(t)$, let $\alpha(t)=\int_0^t c(s)\,ds$ and set
$V_i(x,t)=U_i(x+\alpha(t),t)$. Then $V_i$ satisfies
\begin{equation*}
	\left\{
	\begin{aligned}
		&\partial_{t} V_{i} = c(t) \partial_{x} V_{i}
		+  d_i^{*}(x,t)
		\left( k_{i}*V_i-V_i\right)(x,t)
		+ f_{i}^{*} (x,t,V_1,\cdots,V_m),
		&&t>0,\,x\in\mathbb{R},\\
		&V_{i}(x,0)= U_{0,i}(x),
		&&x\in\mathbb{R},\\
		&i=1,2,\cdots,m,
	\end{aligned}
	\right.
\end{equation*}
where $ d_i^{*}(x,t)=d_i(x+\alpha(t),t)$ and $ f_i^{*}(x,t,V)=f_i(x+\alpha(t),t,V)$. In the spatially homogeneous, time-periodic setting, the eigenvalue problem associated with the linearization of the above system at zero is the whole-line counterpart of \eqref{P-draft}. Thus, the large-domain limit of the principal eigenvalue provides a useful tool for studying the propagation dynamics of \eqref{1.4}. In this paper, we study the spreading speeds of \eqref{1.4}
in the spatially homogeneous, time-periodic setting.

For scalar equations in spatially and temporally homogeneous environments, spreading speeds have been well studied. For foundational results on nonlocal equations, see Weinberger \cite{Weninberger} and Yagisita \cite{Yagisita}. Lutscher et al. \cite{lpl} studied how asymmetric dispersal affects invasion speeds in stream populations. Xu et al. \cite{XLR} developed a forward--backward spreading method and further investigated the effects of kernel asymmetry on spreading speeds. Du et al. \cite{DFN} obtained another characterization through the large-domain limit of the principal eigenvalue in the spatially homogeneous scalar case. Spreading in heterogeneous environments and
cooperative systems has also received considerable attention; see \cite{SZ,Bao,LXZ,BSS,HKLL,LWL,RSZ}. These developments motivate a unified framework for propagation that incorporates temporal heterogeneity, component coupling, and asymmetric dispersal.

We also consider the corresponding nonlinear problem on a bounded interval $\Omega=(l_1,l_2)$, with an inflow boundary condition:
\begin{equation}\label{1.2}
	\left\{
	\begin{aligned}
		&\partial_{t} V_{i} = c(t) \partial_{x} V_{i}+ d_{i}(x,t)\mathcal{K}^{i}_{\Omega}[V_{i}](x,t) + f_{i} (x,t,V_1,\cdots,V_m), &&x\in\Omega,t>0,\\
		&V_{i}(l_1,t)=0,&& \text{if } \sup_{t} c(t)<0, \\
		&V_{i}(l_2,t)=0,&& \text{if } \inf_{t} c(t)>0,\\
		&V_{i}(x,0)= V_{0,i}(x), &&x\in\Omega,\\
		&i=1,2,\cdots,m.
	\end{aligned}
	\right.
\end{equation}
The advection term in \eqref{1.2} is motivated by the combined effects of water flow and active movement on aquatic populations; see \cite{MichalecEtAl2017,MasudEtAl2019}. In this formulation, advection describes continuous transport by water flow, while the kernels describe active dispersal relative to the water. Such dispersal may reinforce or oppose this transport. Retaining advection explicitly, rather than incorporating both mechanisms into a single effective kernel, allows us to study their separate and combined effects on population persistence.

Without advection, related nonlocal systems have been studied extensively; see  \cite{BCV1,BCV2,CDLL,DuNi,ShenVo,SLLY,WFLS,WZ} and the references therein. For systems with advection, fewer results are available. Tang and Dai \cite{TangDai1} established well-posedness for scalar nonlocal equations with time-dependent reaction terms and studied their long-time dynamics in the autonomous setting. Kong et al. \cite{KongEtAl2025} studied the spreading speeds for an autonomous nonlocal Fisher--KPP equation with advection on $\mathbb R$.

In this paper, we first establish the regularity of solutions to the initial value problem under different    smoothness and compatibility assumptions on the initial data (Lemmas~\ref{EX-free}-\ref{EX-free-linear}). We then obtain maximum principles and precise boundary estimates for the time-periodic problem and use them to prove the existence of periodic solutions. Combining these results, we establish a complete characterization of the long-time dynamics of \eqref{1.2} (Theorem \ref{Kaiguan} (i)-(iii)).

Our studies of \eqref{1.2} and \eqref{1.4} are closely connected. Indeed, in the spatially homogeneous setting, we show
that the positive time-periodic solution of \eqref{1.2} converges to that of \eqref{1.4} as the interval expands to $\mathbb R$ (Theorem \ref{Kaiguan} (iv)). Using this convergence and the principal eigenvalue theory developed above, we give a precise characterization of the spreading speeds of \eqref{1.4} and criteria for determining whether accelerated spreading occurs (Theorem \ref{thm:spread-periodic-cooperative} and Remark \ref{acce}).

\subsection{Basic notation and assumptions}

The following notation and assumptions will be used throughout the paper. The case $c(t) \le -c_0 < 0$ can be reduced to $c(t) \ge c_{0} > 0$ by reflecting the spatial variable, the spatially dependent coefficients, and the dispersal kernels. We therefore restrict attention to $c(t)\geq c_0>0$
and impose (H1) below. 

Let $m\geq1$ and set
\begin{equation*}
	\mathbb S=\{1,\ldots,m\},\qquad
	\mathbf0=(0,\ldots,0)^{T},\qquad
	\mathbf1=(1,\ldots,1)^{T}.
\end{equation*}
{Define}
\begin{equation*}
	\mathbb{R}^{m}=\{ u =(u_1,u_2,\cdots,u_m):\, u_{i} \in \mathbb R,\, \forall i\in \mathbb{S}  \}, \quad
	\mathbb{R}^{m}_{+}= \{ u \in \mathbb R^{m}: u_i \ge 0, \forall i \in \mathbb S\}.
\end{equation*}
For $u \in\mathbb R^m$, {define} $|u|=\sum_{i=1}^m|u_i|$. For $u,v\in\mathbb R^m$, {the notation} $u\geq v$ {means that} $u-v\in\mathbb R^m_+$, $u>v$ {means that} $u-v\in\mathbb R^m_+\backslash\{\mathbf0\}$, and $u\gg v$ {means that} $u_i>v_i$ for all $i\in\mathbb S$.

For the linear eigenvalue problem \eqref{P-draft}, the following assumptions are imposed on the coefficients and dispersal kernels.
\begin{itemize}[leftmargin=0.9cm]
\item[(H1)]
$c\in C(\mathbb R)$ is $T$-periodic, and there exists $c_0>0$ such that  
\begin{equation*}
c(t)\geq c_0 \enspace \text{for all} \enspace t\in\mathbb R.
\end{equation*}

\item[(H2)]
	For every $i\in\mathbb S$, $d_i\in C^{1,0}(\overline\Omega\times\mathbb R)$ is {strictly positive on $\overline\Omega\times\mathbb R$} and $T$-periodic in $t$.
	
\item[(H3)]
	For every $i\in\mathbb S$, $k_i\in C(\mathbb R\times\mathbb R)$ is nonnegative and $T$-periodic in $t$, and {satisfies} 
	\begin{equation*}
		k_i(0,t)>0,\qquad
		\int_{\mathbb R}k_i(z,t)\,dz=1,
		\qquad t\in\mathbb R.
	\end{equation*}
	
\item[(H4)]
	The entries of $A=(a_{ij})_{i,j\in\mathbb S}$ belong to $C^{1,0}(\overline\Omega\times\mathbb R)$ and are $T$-periodic in $t$. The matrix $A(x,t)$ is cooperative, that is, $a_{ij}(x,t)\geq0$ for $i\neq j$ 	and all $(x,t)\in\overline\Omega\times\mathbb R$. {Moreover,} there exists $x_0\in\Omega$ such that $A(x_0,t)$ is irreducible for every $t\in\mathbb R$. In view of the cooperativity of $A$, this irreducibility condition is equivalent to requiring that, for every nonempty proper subset $I\subset\mathbb S$ and every $t\in\mathbb R$,
	\begin{equation*}
		\sum_{i\in I}\sum_{j\in\mathbb S\setminus I}a_{ij}(x_0,t)>0.
	\end{equation*}
\end{itemize}
Here, $C^{1,0}(\overline\Omega\times\mathbb R)$ denotes the space of functions whose first-order derivatives with respect to $x$ exist and are jointly continuous in $(x,t)$, together with the functions themselves. No differentiability with respect to $t$ is required.

Set
\[
Q=\Omega\times\mathbb R,\qquad
Q_T=\Omega\times(0,T],\qquad
\overline Q_T=\overline\Omega\times[0,T].
\]
For the variational characterizations introduced below, let
\begin{equation*}
\mathcal D:=\left\{
u\in C(\overline Q;\mathbb R^m)\cap C^1(Q;\mathbb R^m):
u(x,t+T)=u(x,t),
u(x,t)\gg\mathbf0
\text{ for all }(x,t)\in Q
\right\}.
\end{equation*}
Set
\begin{equation}\label{1.7a}
	\Gamma_0=\{l_{\mathrm{2}}\}\times\mathbb R,\qquad
	\mathcal D_0=\{u\in\mathcal D:u=\mathbf0\text{ on }\Gamma_0\}.
\end{equation}
{By continuity, every $u\in\mathcal D$ is nonnegative on
	$\partial\Omega\times\mathbb R$. For $u\in\mathcal D$, define the operator}
\begin{equation}\label{1.8a}
	\mathcal Lu
	=-\partial_tu+c(t)\partial_xu
	+D(x,t)\mathcal K_\Omega[u]+A(x,t)u.
\end{equation}

\subsection{Main results}

\subsubsection{\textbf{{Principal eigenvalue and its variational characterizations}}}

A real eigenvalue of \eqref{P-draft} is called
principal if it admits a strongly positive $T$-periodic
eigenfunction $\varphi\in C^1(\overline Q;\mathbb R^m)$. Here, strong positivity
means $\varphi(x,t)\gg\mathbf0$ for all
$x\in\overline\Omega\setminus\{l_{2}\}$ and $t\in\mathbb R$.

\begin{theorem}[{Existence, uniqueness, and variational characterizations of the principal eigenvalue}]
	\label{EX-PEV}
	{
		Assume (H1)--(H4). Then the following statements hold.
		
		\begin{itemize}
			\item[{\rm (i)}]
			Problem \eqref{P-draft} admits a unique
			principal eigenvalue $\lambda_1$.
			The corresponding eigenfunction $\bm\varphi$
			is strongly positive and unique up to
			multiplication by a positive constant.
			
			\item[{\rm (ii)}]
			The principal eigenvalue satisfies
			\begin{equation*}
				\mu_p=\mu_p'=\lambda_p=\lambda_1,
			\end{equation*}
			where
			\begin{align*}
				&\mu_p=\inf_{u\in\mathcal D_0}
				\sup_{\substack{(x,t)\in Q_T\\i\in\mathbb S}}
				\frac{(\mathcal Lu)_i(x,t)}{u_i(x,t)},
				\qquad
				\mu_p'=\sup_{u\in\mathcal D_0}
				\inf_{\substack{(x,t)\in Q_T\\i\in\mathbb S}}
				\frac{(\mathcal Lu)_i(x,t)}{u_i(x,t)},\\
				&\lambda_p=\inf_{u\in\mathcal D}
				\sup_{\substack{(x,t)\in Q_T\\i\in\mathbb S}}
				\frac{(\mathcal Lu)_i(x,t)}{u_i(x,t)}.
			\end{align*}
		\end{itemize}
	}
\end{theorem}

The quantities $\mu_p$, $\mu_p'$ and $\lambda_p$
	also admit the following equivalent characterizations:
\begin{align*}
	&\mu_p =
	\inf\left\{
	\mu\in\mathbb R \,|\,
	\exists u\in\mathcal D_0 \text{ such that }
	\mathcal Lu\le \mu u \text{ in } Q
	\right\},\\
	&\mu_p^{\prime}=
	\sup\left\{
	\mu\in\mathbb R \,|\, \exists u\in\mathcal D_0
	\text{ such that } \mathcal Lu\ge \mu u
	\ \text{in }Q
	\right\},\\
	&\lambda_p=
	\inf\left\{
	\mu\in\mathbb R \,|\,
	\exists u\in\mathcal D \text{ such that }
	\mathcal Lu\le \mu u \text{ in }Q
	\right\}.
\end{align*}
These quantities are also referred to as
	generalized principal eigenvalues.

	The proof of Theorem~\ref{EX-PEV}(i) is based on
	transforming the problem with advection into
	a nonlocal evolution problem on a moving interval.
	More precisely, the change of variables
	$U(x,t)=V(x-\alpha(t),t)$, where
	$\alpha(t)=\int_0^t c(s)\,ds$,
	removes the advection term and transforms the
	original problem into a nonlocal problem on
	the moving interval
	$(l_1+\alpha(t),l_2+\alpha(t))$.
	This interval moves with velocity $c(t)$,
	while its length remains unchanged.
	
	For the transformed problem, a Poincar\'e map
	$\mathcal P$ is constructed by evolving the solution
	over one period and translating it back to the
	original interval.
	The key step is to establish the strong positivity
	and power compactness of $\mathcal P$.
	Strong positivity follows from the positivity
	properties of the dispersal kernels and the
	irreducibility of the cooperative system.
	To prove power compactness, a sufficiently high
	iterate of $\mathcal P$ is considered so that
	the moving interval has shifted completely
	beyond its initial position.
	Consequently, the direct contribution of the
	initial data through transport vanishes,
	and compactness can be established using
	the nonlocal dispersal term.
	
	The spectral properties of power-compact positive operators (see Lemma~\ref{power-positive}) imply that
	$r(\mathcal P)>0$ is an algebraically simple
	eigenvalue with a strongly positive eigenvector.
	Rescaling the corresponding solution yields
	a strongly positive $T$-periodic eigenfunction
	of \eqref{P-draft}, with the principal eigenvalue
	given by
	$\lambda_1=T^{-1}\log r(\mathcal P)$.

\subsubsection{ \textbf{Large-domain limit of the principal eigenvalue}}
Suppose now that $D(x,t)=D(t)$ and $A(x,t)=A(t)$, and let
$\lambda_l$ be the principal eigenvalue of \eqref{P-draft} on $(-l,l)$. 

For $\nu\in\mathbb R$ and $i\in\mathbb S$, define
\begin{equation*}
	m_i(\nu,t)=\int_{\mathbb R}k_i(z,t)e^{-\nu z}\,dz
	\in[0,+\infty],
\end{equation*}
and let
\begin{equation*}
	\mathcal E
	=\left\{\nu\in\mathbb R:
	\int_0^T m_i(\nu,t)\,dt<\infty
	\text{ for every }i\in\mathbb S\right\}.
\end{equation*}
By (H3), $m_i(0,t)=1$ and hence $0\in\mathcal E$.
For $\nu\in\mathcal E$, set
\begin{equation*}
	\mathsf M(\nu,t)
	=\operatorname{diag}\bigl(m_1(\nu,t)-1,\ldots, m_m(\nu,t)-1\bigr)
\end{equation*}
and
\begin{equation*}
	B_\nu(t)=A(t)+\nu c(t)I_m+D(t)\mathsf M(\nu,t),
\end{equation*}
where $I_m$ is the $m\times m$ identity matrix.
Let $\Phi_\nu(t,s)$ be the evolution matrix of $z'=B_\nu(t)z$,
whose coefficients are $T$-periodic and integrable over one period.
Define
\begin{equation}\label{FLOQUENT}
	\Lambda(\nu)=
	\begin{cases}
		\displaystyle\frac1T\log r\bigl(\Phi_\nu(T,0)\bigr),
		&\nu\in\mathcal E,\\[1mm]
		+\infty,&\nu\notin\mathcal E.
	\end{cases}
\end{equation}
Write $\Lambda_*:=\inf_{\nu\in\mathcal E}\Lambda(\nu)
=\inf_{\nu\in\mathbb R}\Lambda(\nu)$.

\begin{theorem}
	\label{cor:general-kernel-large-domain-limit}
	Assume (H1)--(H4), $D(x,t)=D(t)$, and $A(x,t)=A(t)$. Then
	$\Lambda$ has a unique minimizer $\nu_*\in\mathcal E$, and
	\begin{equation*}
		\lim_{l\to\infty}\lambda_l
		=\Lambda_*
		=\Lambda(\nu_*)\in\mathbb R.
	\end{equation*}
\end{theorem}

\subsubsection{\textbf{{Dynamics on bounded intervals}}}

{
	Consider the nonlinear problem \eqref{1.2}, where
	$F=(f_1,\ldots,f_m)^T$ satisfies the following assumptions.
}

\begin{itemize}[leftmargin=0.7cm]
	\item[(F)]
	\begin{itemize}
		\item[(i)]
		$F\in C^{1,0,2}(\mathbb R\times\mathbb R\times\mathbb R_+^m;
		\mathbb R^m)$ is $T$-periodic in $t$. Moreover, the mixed derivative $\partial_xD_uF$
        exists and is jointly continuous in $(x,t,u)$
        on $\mathbb R\times\mathbb R\times\mathbb R_+^m$.
		
		\item[(ii)]
		$F(x,t,\mathbf0)=\mathbf0$, and there exists
		$\mathbf M = (M_1,M_2,\cdots, M_m)\gg\mathbf0$
		such that $F(x,t,\mathbf M)\leq\mathbf0$
		for all $(x,t)\in\mathbb R^2$.
		
		\item[(iii)]
			$F$ is strictly subhomogeneous in the sense that
		for every $(x,t)\in\mathbb R^2$, $u>\mathbf0$,
		and $0<\theta<1$,
		\begin{equation*}
			F(x,t,\theta u)>\theta F(x,t,u).
		\end{equation*}
		
		\item[(iv)]
		$F$ is cooperative in the sense that
        $\partial_{u_j}f_i(x,t,u)\geq 0$ for all
        $(x,t)\in\mathbb R^2$, $u\in[\mathbf0,\mathbf M]$,
        and $i,j\in\mathbb S$ with $i\neq j$.
        Moreover, there exists $x_0\in\Omega$ such that
       the Jacobian $D_uF(x_0,t,u)$ is irreducible
       for all $t\in\mathbb R$ and
        $u\in[\mathbf0,\mathbf M]$.
	\end{itemize}
\end{itemize}

For the initial value problem associated with \eqref{1.2}, the following spaces are used. For $c(t)\geq c_0>0$, the inflow boundary is located at $x=l_2$. Define
\begin{equation*}
	\mathbb U_{(l_1,l_2)}
	=\{\phi\in C(\overline\Omega;\mathbb R^m):
	\phi(l_2)=\mathbf0\}.
\end{equation*}
For solutions that are continuously differentiable up to the initial boundary point $(l_2,0)$, the boundary condition must be consistent with the evolution equation at $t=0$. This requires the initial datum $\phi$ to satisfy the following compatibility condition:
\begin{equation}\label{eq:inflow-compatibility}
	\phi_i'(l_2)
	=-\frac{d_i(l_2,0)}{c(0)}
	\int_{l_1}^{l_2}k_i(l_2-y,0)\phi_i(y)\,dy,
	\qquad i\in\mathbb S.
\end{equation}
Accordingly, the space of compatible $C^1$ initial data is defined by
\begin{equation*}
	\mathbb V_{(l_1,l_2)}
	=\{\phi\in C^1(\overline\Omega;\mathbb R^m)
	\cap\mathbb U_{(l_1,l_2)}:
	\eqref{eq:inflow-compatibility}\text{ holds}\}.
\end{equation*}
The spaces $\mathbb U_{(l_1,l_2)}$ and $\mathbb V_{(l_1,l_2)}$ are closed linear subspaces of $C(\overline\Omega;\mathbb R^m)$ and $C^1(\overline\Omega;\mathbb R^m)$, respectively, and hence are Banach spaces equipped with the corresponding $C$-norm and $C^1$-norm.

Let $\mathcal D_c$ denote differentiation along the
characteristic curves with tangent vector $(-c(t),1)$
in the $(x,t)$-plane.
For a continuous function $V$, we define its
characteristic derivative by
\begin{equation*}
    \mathcal D_cV(x,t)
    :=
    \lim_{h\to0}
    \frac{
        V\left(x-\int_t^{t+h}c(s)\,ds,t+h\right)-V(x,t)
    }{h},
\end{equation*}
whenever this limit exists.
For $V\in C^1$, we have
$\mathcal D_cV=\partial_tV-c(t)\partial_xV$.
The definition of $\mathcal D_cV$ requires only
differentiability along characteristics;
$\partial_tV$ and $\partial_xV$ need not exist separately.

Throughout this paper, a solution in the characteristic
sense is understood to satisfy the equation with
$\partial_tV-c(t)\partial_xV$ interpreted as
$\mathcal D_cV$.
Under (H1)--(H3) and (F), for every
$V_0\in\mathbb U_{(l_1,l_2)}$ satisfying
$\mathbf0\leq V_0\leq\mathbf M$, problem \eqref{1.2}
admits a unique global solution $V(x,t;V_0)$ in the
characteristic sense.
If, in addition, $V_0\in\mathbb V_{(l_1,l_2)}$,
then the solution is classical. 
This existence result will be proved in Lemma \ref{EX-draft}.

\begin{theorem}\label{Kaiguan}
		Assume (H1)--(H3) and (F).
		Let $V(x,t;V_0)$ be the solution of \eqref{1.2}
with initial function $V_0\in\mathbb U_{(l_1,l_2)}$
satisfying $\mathbf0\leq V_0\leq\mathbf M$
and $V_0\not\equiv\mathbf0$.
		Let $\lambda_{(l_1,l_2)}$ denote the principal
		eigenvalue of the linearization of \eqref{1.2}
		at $\mathbf0$, that is, the principal eigenvalue
		of \eqref{P-draft} with
		$A(x,t)=\mathrm D_uF(x,t,\mathbf0)$.
		Then the following statements hold.
	\begin{itemize}
		\item[(i)]
			If $\lambda_{(l_1,l_2)}>0$, then
			the problem 
			\begin{equation}\label{draft-P}
				\left\{
				\begin{aligned}
					&\partial_tV
					=c(t)\partial_xV+D(x,t)\mathcal K_\Omega[V]+F(x,t,V),
					&&(x,t)\in Q_T,\\
					&V(l_2,t)=\mathbf0,
					&&t\in[0,T],\\
					&V(x,0)=V(x,T),
					&&x\in\overline\Omega
				\end{aligned}
				\right.
			\end{equation}
			admits a unique positive
			$T$-periodic solution $V^*$ in the characteristic sense, satisfying
        $\mathbf0\leq V^*\leq\mathbf M$ and
      $V^*,\mathcal D_cV^* \in C(\overline Q_T;\mathbb R^m)$.
         Moreover,
		\begin{equation*}
			\lim_{n\to\infty}
			\sup_{(x,t)\in\overline Q_T}
			|V(x,t+nT;V_0)-V^*(x,t)|=0.
		\end{equation*}
		
		\item[(ii)]
			If $\lambda_{(l_1,l_2)}<0$, then
			there exist constants $C,\sigma>0$ such that
		\begin{equation}\label{eq-2}
			\mathbf0\leq V(x,t;V_0)
			\leq Ce^{-\sigma t}\mathbf1,
			\qquad (x,t)\in \overline\Omega\times[0,\infty).
		\end{equation}
		
		\item[(iii)]
			If $\lambda_{(l_1,l_2)}=0$, then
		\begin{equation*}
			\lim_{n\to\infty}
			\sup_{(x,t)\in\overline Q_T}
			|V(x,t+nT;V_0)|=0.
		\end{equation*}
		
		\item[(iv)]
			Suppose further that $D(x,t)=D(t)$
			and $F(x,t,u)=F(t,u)$. If
		\begin{equation*}
			\lambda_\infty:=
			\lim_{l\to\infty}\lambda_{(-l,l)}>0,
		\end{equation*}
			then there exists $l_0>0$ such that,
			for every $l\geq l_0$, problem \eqref{draft-P}
			on $(-l,l)$ admits a unique positive
			$T$-periodic solution $V_{l,c}^*$ in the characteristic sense, satisfying
            $\mathbf0\leq V^*_{l,c}\leq\mathbf M$ and
            $V^*_{l,c},\mathcal D_cV^*_{l,c} \in C(\overline Q_T;\mathbb R^m)$.
			Moreover, for every $R>0$,
		\begin{equation*}
			\lim_{l\to\infty}
			\sup_{|x|\leq R,\,t\in[0,T]}
			|V_{l,c}^*(x,t)-U^*(t)|=0,
		\end{equation*}
			where $U^*$ is the unique positive
			$T$-periodic solution of $U'=F(t,U)$
			in $[\mathbf0,\mathbf M]$.
	\end{itemize}
\end{theorem}

\subsubsection{\textbf{{Spreading dynamics in the whole space}}}

	Finally, we investigate the spreading dynamics of
	\eqref{1.4} in a time-periodic environment.
	Assume that $D(x,t)=D(t)$ and $F(x,t,u)=F(t,u)$.
	
	Let
	$\Lambda_0(\nu)$ denote the Floquet exponent
	introduced above, with $A(t)=\mathrm D_uF(t,\mathbf0)$ and $c=0$.
	Assume further that $\Lambda_0(0)>0$,
	which ensures that the zero solution of
	$U'=F(t,U)$ is linearly unstable.
	
	To connect the spectral theory with spreading,
	consider the moving coordinate
	$V^c(x,t)=U(x+ct,t)$, where $c$ is a constant.
	This transformation introduces the advection term
	$c\partial_xV^c$.
	For $c\neq0$, let $\lambda_l(c)$ denote the
	principal eigenvalue of the corresponding
	linearized problem on $(-l,l)$.
	The large-domain limit established above gives
\begin{equation}\label{eq:large-domain-limit-speed}
	\lambda_\infty(c)
	:=\lim_{l\to\infty}\lambda_l(c)
	=\min_{\nu\in\mathcal E}
	\{\Lambda_0(\nu)+c\nu\}.
\end{equation}

	The leftward and rightward spreading speeds are given by
\begin{equation*}
	c_*^-:=\sup_{\mu<0}\frac{\Lambda_0(-\mu)}{\mu},
	\qquad
	c_*^+:=\inf_{\mu>0}\frac{\Lambda_0(-\mu)}{\mu}.
\end{equation*}
	The values $c_*^-=-\infty$ and $c_*^+=+\infty$
	are allowed. When finite, these speeds satisfy
	$\lambda_\infty(c_*^-)=0$ and
	$\lambda_\infty(c_*^+)=0$, respectively.
	The following theorem describes the spreading
	dynamics of \eqref{1.4}.
\begin{theorem}[{Spreading speeds}]
	\label{thm:spread-periodic-cooperative}
		Assume (H2), (H3), and (F), with
		$D(x,t)=D(t)$, $F(x,t,u)=F(t,u)$,
		and $\Lambda_0(0)>0$.
		Let $U(x,t;U_0)$ be the solution of
		\eqref{1.4} with
		$U_0\in C_c(\mathbb R;\mathbb R_+^m)$,
		$\mathbf0\leq U_0\leq\mathbf M$, and
		$U_0\not\equiv\mathbf0$.
		Let $U^{*}$ be the unique positive
		$T$-periodic solution of $U'=F(t,U)$
		in $[\mathbf0,\mathbf M]$.
		Then the following statements hold.
	\begin{itemize}
		\item[(i)]
			For any $c_*^-<a_1\leq b_1<c_*^+$,
		\begin{equation*}
			\lim_{t\to\infty}
			\sup_{x\in[a_1t,b_1t]}
			|U(x,t;U_0)-U^*(t)|=0.
		\end{equation*}
		
		\item[(ii)]
			If $c_*^->-\infty$, then for any
			$a_2<c_*^-$,
		\begin{equation*}
			\lim_{t\to\infty}
			\sup_{x\leq a_2t}
			|U(x,t;U_0)|=0.
		\end{equation*}
		
		\item[(iii)]
			If $c_*^+<+\infty$, then for any
			$b_2>c_*^+$,
		\begin{equation*}
			\lim_{t\to\infty}
			\sup_{x\geq b_2t}
			|U(x,t;U_0)|=0.
		\end{equation*}
	\end{itemize}
\end{theorem}

\begin{remark}\label{acce}{\rm 
(i) The spreading speeds are closely related to the large-domain principal eigenvalue. More precisely, 		$\lambda_\infty(c)$ satisfies
\begin{equation*}
		\lambda_\infty(c)
		\begin{cases}
			>0, & c\in(c_*^-,c_*^+),\\[1mm]
			=0, & c\in\{c_*^-,c_*^+\}\cap\mathbb R,\\[1mm]
			<0, & c<c_*^- \text{ or } c>c_*^+.
		\end{cases}
\end{equation*}
Here, $\lambda_\infty(0)$ is defined by the variational expression on the right-hand side of \eqref{eq:large-domain-limit-speed}. These sign relations will be proved in Proposition~\ref{prop-sign}. Thus, the finite spreading speeds are precisely the critical values at which $\lambda_\infty(c)=0$.
		
(ii) The finiteness of the spreading speeds depends on the tails of the dispersal kernels. Define
\begin{align*}
		\mathcal E^-&:=\left\{\mu>0:
		\int_0^T\int_{-\infty}^{0}
		k_i(z,t)e^{-\mu z}\,dz\,dt<\infty
		\text{ for every }1\le i\le m\right\},\\
		\mathcal E^+&:=\left\{\mu>0:
		\int_0^T\int_{0}^{\infty}
		k_i(z,t)e^{\mu z}\,dz\,dt<\infty
		\text{ for every }1\le i\le m\right\}.
\end{align*}
Then $\mathcal E=\mathcal E^-\cup\{0\}\cup(-\mathcal E^+)$. The corresponding one-sided thin-tail conditions are
\begin{equation*}
		(J_{\mathrm{thin}}^-):\quad
		\mathcal E^-\ne\varnothing,
		\qquad
		(J_{\mathrm{thin}}^+):\quad
		\mathcal E^+\ne\varnothing.
\end{equation*}
The leftward speed $c_*^-$ is finite if and only if $(J_{\mathrm{thin}}^-)$ holds, while the rightward speed $c_*^+$ is finite if and only if $(J_{\mathrm{thin}}^+)$ holds.
		
A kernel $k_i$ is called light-tailed if there exists $\mu>0$ such that
\begin{equation*}
		\int_0^T\int_{\mathbb R} k_i(z,t)e^{\mu|z|}\,dz\,dt<\infty;
\end{equation*}
otherwise, it is called heavy-tailed. If all the kernels are light-tailed, both spreading speeds are finite. If at least one kernel is heavy-tailed, accelerated spreading occurs in at least one direction.
		
In particular, due to the irreducible cooperative coupling, a heavy-tailed kernel in one component can induce 		accelerated spreading in the other components. This phenomenon was established in a homogeneous environment in \cite{XLL}.  In this paper, eigenvalue analysis provides a simple and direct way to extend this result to periodic environments.
}
\end{remark}

\smallskip
\noindent{\bf Organization of the paper.}
The rest of this paper is organized as follows.
Section~2 establishes the well-posedness of the
initial value problem and provides some preliminary results.
Section~3 develops the principal spectral theory,
including the existence, uniqueness, and variational
characterizations of the principal eigenvalue,
as well as its large-domain limit.
Section~4 investigates the long-time dynamics
of the nonlinear problem on bounded intervals
and the large-domain limit of positive periodic solutions.
Section~5 studies the spreading dynamics
of the corresponding Cauchy problem on the whole space.
Finally, the Appendix provides an auxiliary result
on the positive cone.

\section{Preliminaries}
\subsection{Existence, uniqueness, and regularity of solutions}
{This subsection establishes the existence, uniqueness, and regularity of solutions to the initial-boundary value problem \eqref{1.2}. Following the moving-coordinate approach in \cite{DFN}, set
\[
\alpha(t)=\int_0^t c(s)\,ds,\qquad
g(t)=l_1+\alpha(t),\qquad h(t)=l_2+\alpha(t).
\]
If $V$ is a solution of \eqref{1.2}, define $U(x,t)=V(x-\alpha(t),t)$ for $g(t)\leq x\leq h(t)$. By the chain rule, $U$ satisfies the following problem on the moving interval:}

\begin{equation}\label{IBP-free}
\left\{
\begin{aligned}
& \partial_{t} U_{i}(x,t) =  d_{i}^{*}(x,t)\mathcal{J}^{i}_{(g(t),h(t))}[U_{i}](x,t) + q_{i}(x,t,U),  &&x\in(g(t),h(t)),  t>0,\\
&U_{i}(h(t),t)=0, &&t > 0,\\
&U_{i}(x,0)= V_{0,i}(x), &&x\in[l_{1},l_{2}].\\
& i=1,2,\cdots,m.
\end{aligned}
\right.
\end{equation}
Here, for $i\in\mathbb S$,
$d_i^*(x,t)=d_i(x-\alpha(t),t)$ and
$q_i(x,t,u)=f_i(x-\alpha(t),t,u)$. The nonlocal operator is defined by
\begin{equation*}
\mathcal J^i_{(g(t),h(t))}[w](x,t)
=\int_{g(t)}^{h(t)}k_i(x-y,t)w(y,t)\,dy-w(x,t).
\end{equation*}
For $T_1>0$, write
$\Omega_{T_1}^{g,h}=\{(x,t):0<t\leq T_1,\,g(t)<x<h(t)\}$
and $\Omega_\infty^{g,h}=\{(x,t):t>0,\,g(t)<x<h(t)\}$. Their closures include the relevant initial and moving boundary points.

\begin{lemma}\label{EX-free}
Assume (H1)--(H3) and (F). For any
$V_0\in\mathbb U_{(l_1,l_2)}$ satisfying
$\mathbf0\leq V_0\leq\mathbf M$, problem \eqref{IBP-free}
admits a unique global solution
$U\in C(\overline\Omega_\infty^{g,h};\mathbb R^m)$,
with $\mathbf0\leq U\leq\mathbf M$.
If, in addition, $V_0\in\mathbb V_{(l_1,l_2)}$, then
$U\in C^1(\overline\Omega_\infty^{g,h};\mathbb R^m)$.
\end{lemma}
\begin{proof}
The proof follows the parameter-dependent ODE and fixed-point approach of \cite[Lemma~2.3]{CDLL} and \cite[Lemma~5.3]{DFN}, with the details needed for the present system.

\textbf{Step 1.} An ODE with the spatial variable as a parameter.

Fix any $T_{1}>0$. For $x\in[l_1,h(T_1)]$ and $i\in\mathbb S$, define
\begin{equation*}
\tilde{u}_{0,i}(x) =
\left\{
\begin{aligned}
& V_{0,i}(x),
&& x\in [l_{1},l_{2}],\\
& 0, && x>l_{2}.
\end{aligned}
\right.
\end{equation*}
Since $c(t)>0$ on $[0,T_1]$, both $g$ and $h$ are strictly
increasing $C^1$ functions. Hence
\begin{equation*}
h:[0,T_1]\to[l_2,h(T_1)]
\quad\text{and}\quad
g:[0,T_1]\to[l_1,g(T_1)]
\end{equation*}
are $C^1$ bijections. Denote their inverses by $h^{-1}$ and $g^{-1}$. Thus,
$h^{-1}(x)$ and $g^{-1}(x)$ are the unique times in $[0,T_1]$
such that
\begin{equation*}
    h\left(h^{-1}(x)\right)=x
    \quad\text{and}\quad
    g\left(g^{-1}(x)\right)=x,
\end{equation*}
respectively.
For $x\in[l_1,h(T_1)]$, define
\begin{equation*}
t_0(x)
=
\begin{cases}
0,
&x\in[l_1,l_2],\\
h^{-1}(x),
&x\in(l_2,h(T_1)],
\end{cases}
\quad \text{and} \quad
t_{\max}(x)
=
\begin{cases}
g^{-1}(x),
&x\in[l_1,g(T_1)],\\
T_1,
&x\in(g(T_1),h(T_1)].
\end{cases}
\end{equation*}
For a fixed $x$, the interval $[t_0(x),t_{\max}(x)]$ is exactly the time interval during which $x$ lies in $[g(t),h(t)]$; that is,
$g(t)\leq x\leq h(t)
\leftrightarrow
t_0(x)\leq t\leq t_{\max}(x)$.

Define
\begin{align*}
&\mathbb{X}_{T_1} = \!\bigg\{ \phi\in C(\overline{\Omega}_{T_1}^{g,h},\mathbb{R}^{m}) : 0\leq \phi \leq \mathbf{M} ,\, \phi(x,0) = V_{0}(x),\\
&\quad\quad\ \ \ \  \forall x\in [l_{1},l_{2}],\,\phi(h(t),t)=0,  \forall t\in [0,T_{1}]  \bigg\},\\
&Q(x,t,v)=(q_{1}(x,t,v),\cdots,q_{m}(x,t,v))^{T},\ \ \ \  \tilde{v}_{i} = v_{i}^{+},\\ 
&\tilde{v}= (\tilde{v}_{1},\tilde{v}_{2},\cdots,\tilde{v}_{m}),\,
\tilde{q}_{i}(x,t,v) = q_{i}(x,t,\tilde{v}).
\end{align*}

Fix $\phi\in\mathbb X_{T_1}$. For each spatial point $x$, consider the following ODE on $[t_0(x),t_{\max}(x)]$:
\begin{equation}\label{ODE}
\left\{
\begin{aligned}
& v_{i,t} = d_{i}^{*}(x,t) \left( \int_{ g(t) }^{ h(t) } k_{i}(x-y,t)\phi_{i}(y,t)dy - v_{i}(x,t)  \right)  + \tilde{q}_{i}(x,t,v),  && t_{0}(x )<t\le t_{max}(x),\\
&v_{i}(x,t_{0}(x))= \tilde{u}_{0,i}(x), &&x\in[l_{1},h(T_{1})),\\
& i=1,2,\cdots,m.
\end{aligned}
\right.
\end{equation}
Let $G_{i}(x,t,v) = d_{i}^{*}(x,t) \left( \int_{ g(t) }^{ h(t) } k_{i}(x-y,t)\phi_{i}(y,t)dy - v_{i} \right) + \tilde{q}_{i}(x,t,v)$. 
Since $q_i$ is continuously differentiable in $v$, $G_i(x,t,v)$ is locally Lipschitz continuous in $v$ on $\Pi_{i=1}^{m} (-\infty,M_{i}+1]$ with a Lipschitz constant uniform for $(x,t)\in \overline{\Omega}^{g,h}_{T_1}$. Moreover, $G_{i}(x,t,v)$ is continuous in all its variables in this range. 
The standard local existence and uniqueness theorem for ODEs therefore implies that for every given $x\in [l_{1},h(T_{1}))$, problem \eqref{ODE} admits a unique solution, denoted by $V^{\phi}(x,t)$, which is initially defined on a time interval starting at $t_0(x)$. 

The following upper and lower bounds show that $V^\phi(x,\cdot)$ extends uniquely to $[t_{0}(x), t_{max}(x)]$. 
Note that
\begin{equation*}
d_{i}^{*}(x,t) M_{i} \int_{ g(t) }^{ h(t) } k_{i}(x-y,t)dy - d_{i}^{*}(x,t)M_{i}  + \tilde{q}_{i}(x,t,\mathbf{M}) \leq 0  
\end{equation*}
for all $x\in [l_1, h(T_1))$ and $t\in (t_0(x),t_{max}(x)]$,
and $\widetilde{u}_{0,i}(x)\leq M_i$ in $[l_{1},h(T_{1}))$.
Thus $\mathbf M$ is an upper solution of \eqref{ODE}. The vector field is cooperative on $[\mathbf0,\mathbf M]$, and $\mathbf0$ is a lower solution because $\phi\geq\mathbf0$ and $Q(x,t,\mathbf0)=\mathbf0$. The comparison principle for cooperative ODE systems gives $\mathbf0\leq V^\phi\leq\mathbf M$. These bounds prevent finite-time blow-up and extend the solution to $[t_0(x),t_{\max}(x)]$.

Each class $\mathbb X_s$ is equipped with the supremum metric inherited from $C(\overline\Omega_s^{g,h};\mathbb R^m)$. In particular,
$\|\phi-\psi\|_{\mathbb X_s}=\sup_{(x,t)\in\overline\Omega_s^{g,h}}|\phi(x,t)-\psi(x,t)|$. 

\textbf{Step 2.} A fixed-point argument.

By Step 1, $\mathbf0\leq V^\phi\leq\mathbf M$, $V^\phi(x,0)=V_0(x)$ on $[l_1,l_2]$, and $V^\phi(h(t),t)=\mathbf0$ for $0\leq t\leq T_1$.

Moreover, by continuous dependence of ODE solutions on the initial data and the parameter $x$, we also see that $V^{\phi}$ is continuous in $(x,t)\in \overline{\Omega}_{T_1}^{g,h}$.
Hence, $V^{\phi}\in \mathbb{X}_{T_1}$. 

For $0<s\leq T_1$, let
$\overline{\Omega}_{s}^{g,h}$ and $\mathbb X_s$ denote the
corresponding domain and function class with $T_1$ replaced by $s$, respectively. 
Applying Step~1 with $T_1=s$, for every
$\phi\in\mathbb X_s$ we obtain a function
$V^\phi\in\mathbb X_s$. Define
\begin{equation*}
    \Gamma_s:\mathbb X_s\to \mathbb X_s,
    \quad
    \Gamma_s\phi:=V^\phi.
\end{equation*}
Let $\varphi,\psi\in \mathbb{X}_{s}$, $W = \Gamma_{s}\varphi - \Gamma_{s}\psi$ and
\begin{equation*}
    H_{i}(x,t) =  d_{i}^{*}(x,t)\int_{g(t)}^{h(t)}k_{i}(x-y,t) (\varphi_{i}(y,t) - \psi_{i}(y,t))dy + 
     \tilde{q}_{i}(x,t,\Gamma_{s}\varphi)-\tilde{q}_{i}(x,t,\Gamma_{s}\psi).
\end{equation*}
Then, a direct calculation gives that
\begin{equation*}
\begin{aligned}
\sum_{i=1}^{m} |W_{i}(x,t)| 
\le& \sum_{i=1}^{m} \int_{t_{0}(x)}^{t} e^{-\int_{\xi}^{t} d_{i}^{*}(x,\tau)d\tau} |H_{i}(x,\xi)| d\xi 
\\
\le& sm \max_{(x,t,u)\in \overline{\Omega}_{T_1}^{g,h}\times [0, \bf{M}] }  \|\mathbf{D}_u Q(x,t,u)\| \| W \|_{C(\overline{\Omega}_{s}^{g,h},\mathbb{R}^{m})}\\
 &+ sm \max_{1\leq i \leq m}\| d_{i}^{*}\|_{\infty} \|\varphi - \psi \|_{C(\overline{\Omega}_{s}^{g,h},\mathbb{R}^{m})}
\end{aligned}
\end{equation*}
for $(x,t)\in \overline{\Omega}_{s}^{g,h}$. Let 
\begin{equation*}
  \overline{d}=\max_{1\leq i \leq m}\| d_{i}^{*}\|_{\infty},\enspace  
  C_{0}=\max_{(x,t,u)\in \overline{\Omega}_{T_1}^{g,h}\times [\bf{0},\bf{M}] }  \|\mathbf{D}_uQ(x,t,u)\|  \enspace\text{and}\enspace
  C_{1}= \frac{sm \overline{d}}{1-sm C_{0}}.
\end{equation*}
Hence,
\begin{equation*}
    \| \Gamma_{s}\varphi - \Gamma_{s}\psi\|_{\mathbb{X}_{s}} \leq C_{1} \| \varphi-\psi \|_{\mathbb{X}_{s}} \leq \frac{1}{2}\| \varphi-\psi \|_{\mathbb{X}_{s}} \enspace \text{for}\enspace s\in(0,\delta],
\end{equation*}
provided that $\delta\in (0,\delta_{1}]$ with $\delta_{1} \leq \frac{1}{(2\overline{d}+C_{0})m}$.
For such $s$, the Banach fixed-point theorem gives a unique fixed point of $\Gamma_s$ in $\mathbb X_s$. This fixed point is the unique solution of \eqref{IBP-free} on $0\leq t\leq s$ satisfying $\mathbf0\leq U\leq\mathbf M$.

\textbf{Step 3.} Extension to an arbitrary finite time.

The contraction time can be chosen uniformly on each fixed interval $[0,T_1]$, since the relevant coefficients and their Lipschitz bounds are uniformly bounded there. Restarting the construction at the end of each short interval extends the solution step by step up to $T_1$. As $T_1$ is arbitrary, the solution is global.

\textbf{Step 4.} Regularity for compatible $C^1$ initial data.

Consider the subclass
$\mathbb X^1_{T_1}:=\mathbb X_{T_1}\cap
C^1(\overline{\Omega}_{T_1}^{g,h};\mathbb R^m)$
of $\mathbb X_{T_1}$.
For $\phi \in \mathbb{X}^{1}_{T_1}$, 
\begin{equation*}
\int_{g(t)}^{h(t)} k_{i}(x-y,t)\phi_{i}(y,t)dy = \int_{x-h(t)}^{x-g(t)}k_{i}(z,t)\phi_{i}(x-z,t)dz
\text{ is } C^{1} \text{ in }  x \text{ for }  i \in \mathbb{S}.
\end{equation*}
Since $V^\phi\ge0$, we may replace $\tilde q_i$ by $q_i$
in \eqref{ODE}. By the assumed regularity of the coefficients
and the $C^1$ regularity of $V_0$, the standard theorem on
differentiable dependence of ODE solutions on parameters
implies that the solution $V^\phi$ constructed in Step 1 is
continuously differentiable in $x$ for $x\in[l_1,l_2]$ and
$t\in[0,t_{\max}(x)]$.
Similarly, since $t_0(x)=h^{-1}(x)$ is continuously differentiable
on $[l_2,h(T_1)]$ and $V^\phi(x,t_0(x))=0$,
the same conclusion holds for $x\in[l_2,h(T_1)]$ and
$t\in[t_0(x),t_{\max}(x)]$.
Here the differentiation in $x$ at the endpoints of each of the two intervals is understood as one-sided differentiation.
It remains to verify that the spatial derivatives match at $x=l_2$. More precisely, we show that
\begin{equation*}
\lim_{x\to l_{2}^{-}} \partial_{x} V^{\phi}(x,t_0(x)) = \lim_{x \to l_{2}^{+}} \partial_{x} V^{\phi}(x,t_{0}(x)).
\end{equation*}
For $x
\neq  l_2$, set $W^\phi(x,t)=\partial_xV^\phi(x,t)$. Differentiating the parameter-dependent ODE yields
\begin{equation}\label{W-ODE}
\partial_{t}W^{\phi}(x,t)=A^{\phi}(x,t)W^{\phi}(x,t)+F^{\phi}(x,t),
\end{equation}
for $x \in [l_{1},h(T_{1})] \backslash \{ l_{2} \}$ and $t \in (t_0(x),t_{max}(x))$, where
\begin{equation*}
A^{\phi}(x,t)= \mathbf{D}_{u} Q(x,t,V^{\phi}(x,t))-\operatorname{diag}(d_{1}^{*}(x,t),\dots,d_{m}^{*}(x,t)),
\end{equation*}
and
\begin{equation*}
\begin{aligned}
F_{i}^{\phi}(x,t)
=\,&
\partial_{x}d_{i}^{*}(x,t)\int_{g(t)}^{h(t)}k_{i}(x-y,t)\phi_{i}(y,t)\,dy
+
d_{i}^{*}(x,t)k_{i}(x-g(t),t)\phi_{i}(g(t),t)
\\
&+
d_{i}^{*}(x,t)\int_{x-h(t)}^{x-g(t)}k_{i}(z,t)\,\partial_{x}\phi_{i}(x-z,t)\,dz
-
\partial_{x}d_{i}^{*}(x,t)V_{i}^{\phi}(x,t)
+
\partial_{x}q_{i}(x,t,V^{\phi}(x,t)).
\end{aligned}
\end{equation*}

At the moving boundary $x=h(t)=l_2+\alpha(t)$,
differentiating the identity $V_i^\phi(h(t),t)=0$ with respect to $t$ yields the following one-sided derivative relation:
\begin{align*}
&\partial_{x} V^{\phi}_{i}(l_{2}+\alpha(t),t) \\
=& -\frac{1}{c(t)} \partial_{t} V^{\phi}_{i}(l_{2}+\alpha(t),t) \\ 
=& -\frac{1}{c(t)} \lim_{x \to l_{2}+\alpha(t)} \left( d_{i}^{*}(x,t)\int_{ g(t) }^{ h(t) } k_{i}(x-y,t)\phi_{i}(y,t)dy - d_{i}^{*}(x,t)V_{i}^{\phi}(x,t)  + q_{i}(x,t,V^{\phi}) \right)\\
=& -\frac{d_{i}^{*}( h(t),t)}{c(t)} \int_{ g(t) }^{ h(t) } k_{i}(h(t)-y,t)\phi_{i}(y,t)dy.
\end{align*}
It follows that, as $t$ decreases to $0$,
\begin{equation*}
\partial_{x} V^{\phi}_{i}(l_{2}+\alpha(t),t) \to -\frac{d_{i}^{*}(h(0),0)}{c(0)}\int_{l_{1}}^{l_{2}}k_{i}(l_{2}-y,0)V_{0,i}(y)dy.
\end{equation*}
Since $h(0)=l_2$ and $d_i^*(h(0),0)=d_i(l_2,0)$, the last limit becomes
\begin{equation*}
\partial_{x} V^{\phi}_{i}(l_{2}+\alpha(t),t) \to -\frac{d_{i}(l_{2},0)}{c(0)}\int_{l_{1}}^{l_{2}}k_{i}(l_{2}-y,0)V_{0,i}(y)dy = V^{\prime}_{0,i}(l_{2}).
\end{equation*}
Since $\partial_{x}V^{\phi}_{i}(x,0) = V_{0,i}^{\prime}(x)$ in $ [l_{1},l_{2}]$, the above identity implies that $W_{0,i}^{\phi}(x) = \partial_{x} V^{\phi}_{i}(x, t_{0}(x))$ is continuous over $[l_{1},h(T_{1}))$ for $1\leq i \leq m$. Moreover,
\begin{equation}\label{2.4}
|W_{0}^{\phi}(x)| \leq R_{0}:= \max\left\{ \max_{x\in[l_1,l_2]}|V_0'(x)| , |\mathbf{M} |\frac{\max_{1\le i \le m} \| d_{i}^{*}\|_{\infty}}{\min_{t\in[0,T]} c(t)} \right\} \enspace \text{for} \enspace x\in[l_{1}, h(T_{1})].
\end{equation}
Regarding $W^{\phi}(x,t)$ as the ODE solution of \eqref{W-ODE} with parameter $x$ and initial condition $$W_{i}^{\phi}(x, t_{0}(x)) = W_{0,i}^{\phi}(x), \quad 1\leq i \leq m,$$ 
we can conclude from standard ODE theory that $W^{\phi} \in C(\overline{\Omega}_{T_1}^{g,h};\mathbb{R}^{m})$. Moreover, we can use \eqref{ODE} to see that $\partial_{t}V^{\phi} \in C(\overline{\Omega}_{T_1}^{g,h};\mathbb{R}^{m})$. Hence $V^{\phi}\in \mathbb{X}^{1}_{T_1}$.

\textbf{Step 5.} A fixed-point argument.

Let $\mathbb{X}_{s}^{1} = \mathbb{X}_{s} \cap C^{1}(\overline{\Omega}_{s}^{g,h};\mathbb{R}^{m})$, where $\mathbb{X}_{s}$ and $\overline{\Omega}_{s}^{g,h}$ are defined in Step 2. From the above analysis and Step 2, we can obtain that $\Gamma_{s}$ maps $\mathbb{X}_{s}^{1}$ into itself. 
The next aim is to obtain a fixed point of $\Gamma_s$ in $\mathbb X_s^1$,   denoted by $\psi$, by applying the Schauder fixed-point theorem on a bounded subset.
 Since Step 2 gives a unique fixed point $\phi$ in $\mathbb X_s$ and $\mathbb X_s^1\subset\mathbb X_s$, we must have $\psi = \phi$ and thus $V^{\phi}=\phi$ is in particular continuously differentiable in $x$ for $t\in[0,s]$, $x\in [g(t),h(t)]$.

Let 
$$C_{A} =  \sup_{(x,t,v)\in \overline{\Omega}_{T_1}^{g,h}\times [0,\mathbf M]}\| \mathbf{D}_{v}Q(x,t,v)\|
+
\max_{1\le i\le m}
\|d_i^*\|_{L^\infty(\overline\Omega_{T_1}^{g,h})}.$$ 
Then $\| A^{\phi}\| \leq C_{A}$. By the definition of $F_{i}^{\phi}$, there are constants $C_2,C_3>0$, depending only on the coefficients, $\mathbf M$, and the fixed time interval, such that
\begin{equation*}
    |F^{\phi}(x,t)| \leq C_{2} + C_{3} \|\phi_{x}\|_{C(\overline{\Omega}_{s}^{g,h},\mathbb{R}^{m})}.
\end{equation*}
Using \eqref{ODE} and $V^{\phi} \in \mathbb{X}_{s}^{1}$ we easily deduce 
\begin{equation}\label{Vt}
\| \partial_{t} V^{\phi}(x,t) \|_{\infty} \leq C_{4}:= 2\sum_{i =1}^{m}  M_{i}\| d_{i}^{*} \|_{\infty} + \max_{(x,t,v)\in \overline\Omega_{T_1}^{g,h}\times \Pi_{i=1}^{m} [0, M_{i}] }\|Q(x,t,v) \|.
\end{equation}

Define $R_1=2e^{C_A}(R_0+C_2)$ with $R_0$ given by \eqref{2.4}, and 
\begin{equation*}
\enspace \mathbb{X}^{1}_{s,R_{1}} = \{ \phi\in \mathbb{X}_{s}^{1} : \| \phi_{t}\|_{C( \overline{\Omega}_{s}^{g,h},\mathbb{R}^{m})} \le C_{4},\, \|\phi_{x}\|_{C(\overline{\Omega}_{s}^{g,h},\mathbb{R}^{m})}\le R_{1}  \}.
\end{equation*}
For sufficiently small $s>0$, chosen independently of $R_1$,
$\phi \in \mathbb{X}^{1}_{s,R_1} \Rightarrow V^{\phi} \in \mathbb{X}^{1}_{s,R_1}$.
It suffices to show 
\begin{equation*}
 \| \partial_{x} V^{\phi}(x,t) \|_{C(\overline{\Omega}_{s}^{g,h},\mathbb{R}^{m})} \leq R_{1}.
\end{equation*}
Let $\Phi^\phi(x,t,\tau)$ be the evolution matrix of
${dY}/{dt}=A^\phi(x,t)Y$.
 Since $\|A^{\phi}\|\le C_{A}$,
standard ODE theory yields
\begin{equation*}
\|\Phi^{\phi}(x,t,\tau)\|\le e^{C_{A}(t-\tau)},
\qquad
t_{0}(x)\le \tau\le t\le \min\{ s, t_{\max}(x) \}.
\end{equation*}
Using the variation-of-constants formula, we get 
\begin{equation*}
\begin{aligned}
|W^\phi(x,t)| &\le
\|\Phi^\phi(x,t,t_0(x))\| |W_0^\phi(x)|+
\int_{t_0(x)}^t \|\Phi^\phi(x,t,\tau)\|
|F^\phi(x,\tau)| d\tau\\
&\le e^{C_A(t-t_0(x))}R_0 +(C_2+C_3R_1)
\int_{t_0(x)}^t e^{C_A(t-\tau)}\,d\tau\\
&\le e^{C_As}\bigl(R_0+sC_2+sC_3R_1\bigr).
\end{aligned}
\end{equation*}
Consequently, for $0<s\le\min\{T_1,1\}$,
\begin{equation*}
\begin{aligned}
\|\partial_xV^\phi\|_{\mathbb X_s}
\le e^{C_A}(R_0+C_2)+s e^{C_A}C_3R_1 =\frac{R_1}{2}+s e^{C_A}C_3R_1.
\end{aligned}
\end{equation*}
Set
\begin{equation*}
\delta_2 =
\min\left\{1,\frac{1}{2C_3e^{C_A}}\right\},
\qquad
\delta_0=\min\{T_1,\delta_1,\delta_2\}.
\end{equation*}
Then, for $0<s\le\delta_0$,
\begin{equation*}
    \|\partial_xV^\phi\|_{\mathbb X_s}
    \le\frac{R_1}{2}+\frac{R_1}{2}=R_1.
\end{equation*}
Together with \eqref{Vt}, this proves that
$\Gamma_s(\mathbb X_{s,R_1}^1)\subseteq\mathbb X_{s,R_1}^1$
for $0<s\le\delta_0$.

To see that $X^1_{s,R_1}$ is nonempty, set
$\phi^0(x,t)=V_0(x-\alpha(t))$.
Then $\phi^0\in X_s^1$ and
$$\|\partial_x\phi^0\|_{C(\overline{\Omega}_s^{g,h};\mathbb R^m)}\le R_0\le R_1.$$
The preceding estimates imply that
$\Gamma_s\phi^0\in X^1_{s,R_1}$ for $0<s\le\delta_0$.
The set $\mathbb X_{s,R_1}^1$ is a closed convex subset of the Banach space $\mathbb{Y}= C^{1}(\overline{\Omega}_{s}^{g,h},\mathbb{R}^{m})$ with norm 
\begin{equation*}
\| \psi \|_{\mathbb{Y}} := \| \psi \|_{C(\overline{\Omega}_{s}^{g,h},\mathbb{R}^{m})} + \| \psi_{t} \|_{C(\overline{\Omega}_{s}^{g,h},\mathbb{R}^{m})} + \| \psi_{x}\|_{C(\overline{\Omega}_{s}^{g,h},\mathbb{R}^{m})}.
\end{equation*}
Moreover, we derive from
\eqref{ODE} and \eqref{W-ODE} that $\Gamma_s$ is continuous on $\mathbb X_{s,R_1}^1$
with respect to the $C^1$ norm. 
It remains to show that the image of this set under $\Gamma_s$ is relatively compact in $C^1(\overline\Omega_s^{g,h};\mathbb R^m)$. The Schauder fixed-point theorem will then provide a fixed point in $\mathbb X_{s,R_1}^1$.

We now examine $V^{\phi}$ for a given $\phi\in \mathbb{X}^{1}_{s,R_1}$. 
The argument in \cite[Lemma~5.3, Step~4]{DFN} shows that the families $\{V^\phi\}$ and $\{\partial_tV^\phi\}$ are equicontinuous on $\overline\Omega_s^{g,h}$ as $\phi$ varies in $\mathbb X_{s,R_1}^1$.
It remains to prove equicontinuity of the family $W^\phi=\partial_xV^\phi$. By \eqref{W-ODE}, we have 
\begin{equation}\label{W2}
    W^{\phi}(x,t) = W_{0}^{\phi}(x) + \int^{t}_{t_{0}(x)} \left( A^{\phi}(x,\tau)W^{\phi}(x,\tau)+F^{\phi}(x,\tau)  \right)d\tau.
\end{equation}
Since $A^{\phi}(x,t)W^{\phi}(x,t)+F^{\phi}(x,t)$ is uniformly bounded in $\overline{\Omega}_{s}^{g,h}$, we derive from \eqref{W2} that $W^\phi(x,t)$ is
equicontinuous in $t$ on
$[t_0(x), \min\{s,t_{\max}(x)\}]$,
uniformly in $x$ and $\phi$.
On the other hand, since
\begin{equation*}
W_{0,i}^{\phi}(x) = 
\left\{
\begin{aligned}
& V_{0,i}^{\prime}(x),
&& x\in [l_{1},l_{2}],\\
& -\frac{d_{i}^{*}( x,t_{0}(x))}{c(t_{0}(x))}  \int_{ g(t_{0}(x)) }^{ h(t_{0}(x)) } k_{i}(x-y,t_0(x))\phi_{i}(y,t_{0}(x))dy, && x\in (l_{2},h(s)],
\end{aligned}
\right.
\end{equation*}
we see that $W_{0,i}^{\phi}(x)$ is equicontinuous in $x$ on $[g(t),h(t)]$,
uniformly in $t\in[0,s]$ and
$\phi\in X^1_{s,R_1}$.
Since $V^\phi$ is equicontinuous on
$\overline\Omega_s^{g,h}$ and $D_vQ$ is continuous
on $\overline\Omega_s^{g,h}\times[\mathbf0,\mathbf M]$,
the family $A^\phi(\cdot,t)$ is equicontinuous
on $[g(t),h(t)]$, uniformly in $t\in[0,s]$
and $\phi\in X^1_{s,R_1}$.
Moreover, since $\partial_{x}\phi_{i}$ is bounded for $\phi\in \mathbb{X}^{1}_{s,R_1}$, 
\begin{equation*}
d_{i}^{*}(x,t)\int_{x-h(t)}^{x-g(t)}k_{i}(z,t)\partial_{x}\phi_{i}(x-z,t)dz = d_{i}^{*}(x,t)\int_{g(t)}^{h(t)}k_{i}(x-y,t)\partial_{x}\phi_{i}(y,t)dy
\end{equation*}
is equicontinuous in $x$ on $[g(t),h(t)]$,
uniformly in $t\in[0,s]$ and
$\phi\in X^1_{s,R_1}$. 
By the definition of $F_i^\phi$, the family
$F_i^\phi(\cdot,t)$ is also equicontinuous
on $[g(t),h(t)]$, uniformly in $t\in[0,s]$
and $\phi\in X^1_{s,R_1}$, for every $i\in\mathbb S$.
Let $\omega_{0}$, $\omega_{F}$, and $\omega_{A}$ be the moduli of continuity of the functions  
$W_{0}^{\phi}$, $F^{\phi}$ and $A^{\phi}$, respectively. Set $C_{5}=C_{A}R_{1}+C_{2}+C_{3}R_{1}$.
Consequently, for any $g(t) \le x_{1} \le x_{2} \le h(t)$,
\begin{equation*}
\begin{aligned}
& |W^\phi(x_1,t)-W^\phi(x_2,t)|
\\
\le\,&
|W_0^\phi(x_1)-W_0^\phi(x_2)|
+
\int_{t_0(x_1)}^{t_0(x_2)}
|
A^\phi(x_1,\tau)W^\phi(x_1,\tau)+F^\phi(x_1,\tau)
|d\tau
\\
&+
\int_{t_0(x_2)}^{t}
|A^\phi(x_1,\tau)W^\phi(x_1,\tau)-A^\phi(x_2,\tau)W^\phi(x_2,\tau)|d\tau
+
\int_{t_0(x_2)}^{t}|F^\phi(x_1,\tau)-F^\phi(x_2,\tau)|d\tau\\
\le&\, \omega_{0}(x_1,x_2)+ (C_{A}R_{1}+C_{2}+C_{3}R_{1})|t_{0}(x_{2}) -t_{0} (x_{1})| + \int_{t_{0}(x_{2})}^{t} \omega_{F}(x_1,x_2) d\tau \\
&+ \int_{t_0(x_2)}^{t} \left( C_A |W^\phi(x_1,\tau)-W^\phi(x_2,\tau)|
+
R_1\omega_A(x_1,x_2) \right)
\,d\tau
\\
\le & \omega_{0}(x_1,x_2) + \frac{C_{5}}{c_{0}}|x_{2} -x_{1}| + s\omega_{F}(x_1,x_2) + sR_{1}\omega_{A}(x_1,x_2)\\
&+ C_A\int_{t_0(x_2)}^{t}  |W^\phi(x_1,\tau)-W^\phi(x_2,\tau)| d\tau.
\end{aligned}
\end{equation*}
Applying Gr\"onwall's inequality yields
\begin{equation*}
\begin{aligned}
&|W^\phi(x_1,t)-W^\phi(x_2,t)| \\
\leq&
e^{C_{A}t}\left( \omega_{0}(x_1,x_2) + \frac{C_{5}}{c_{0}}|x_{2} -x_{1}| + s\omega_{F}(x_1,x_2) + sR_{1}\omega_{A}(x_1,x_2) \right).
\end{aligned}
\end{equation*}
This proves equicontinuity in the spatial variable on the moving domain. Thus $V^{\phi}(x,t)$, $\partial_{t} V^{\phi}(x,t)$ and  $\partial_{x} V^{\phi}(x,t)$ are all uniformly bounded and equicontinuous over $\overline\Omega_s^{g,h}$ for $\phi\in \mathbb{X}^{1}_{s,R_1}$. This implies that 
$
\Gamma_{s}( \mathbb{X}^{1}_{s,R_1}) = \{  V^{\phi} : \phi\in \mathbb{X}^{1}_{s,R_1}\}
$
is precompact in $\mathbb{Y}$.
Hence, the unique solution of \eqref{IBP-free} is $C^1$ on $\overline\Omega_{\delta_0}^{g,h}$.

To extend this regularity beyond $\delta_0$, it is necessary to verify the compatibility condition at each later initial time.
 From $U(l_{2}+\alpha(t),t)=0$ we deduce 
$\partial_{t} U_{i}(l_{2}+\alpha(t),t) + c(t) \partial_{x}U_{i}(l_{2}+\alpha(t),t)=0$
and so 
\begin{align*}
\partial_{x}U_i(l_{2}+\alpha(t),t) 
=& -\frac{1}{c(t)}\partial_{t} U_i(l_{2}+\alpha(t),t)\\
=& -\frac{d_{i}^{*}(l_{2}+\alpha(t),t)}{c(t)}\int_{g(t)}^{h(t)}k_{i}(l_{2}+\alpha(t)-y,t)U_i(y,t)dy\\
=& -\frac{d_{i}(l_{2},t)}{c(t)}\int_{l_{1}}^{l_{2}}k_{i}(l_{2}-y,t)U_i(y+\alpha(t),t)dy
\end{align*}
for $i \in \mathbb{S}$. Hence, 
$U(\cdot+\alpha(\delta_0),\delta_0)$ satisfies the
compatibility condition at time $\delta_0$.
Regarding $t=\delta_{0}$ as the initial time of \eqref{IBP-free} we can show in the same way that the solution is continuously
differentiable in $(x,t)$ for $t \in [\delta_{0}, 2\delta_{0}]$ and $x \in [g(t), h(t)]$.
Repeating this argument finitely many times yields
$U\in C^1(\overline\Omega_{T_1}^{g,h};\mathbb R^m)$.
Since $T_1>0$ was arbitrary, we conclude that
$U\in C^1(\overline\Omega_\infty^{g,h};\mathbb R^m)$.
This completes the proof of the lemma.
\end{proof}

Applying the inverse change of variables
$V(x,t)=U(x+\alpha(t),t)$, we obtain the following existence,
uniqueness, and regularity result for \eqref{1.2} from
Lemma~\ref{EX-free}.

\begin{lemma}\label{EX-draft}
Assume (H1)--(H3) and (F).
For any $V_0\in\mathbb U_{(l_1,l_2)}$ satisfying
$\mathbf0\leq V_0\leq\mathbf M$, problem \eqref{1.2}
admits a unique global solution $V(x,t;V_0)$
in the characteristic sense, with
$
V, \mathcal{D}_{c} V \in C\left( \overline\Omega\times[0,\infty);\mathbb R^m \right), \mathbf0\leq V\leq\mathbf M.
$
If, in addition, $V_0\in\mathbb V_{(l_1,l_2)}$, then
$V\in C^1\left([l_1,l_2]\times[0,\infty);\mathbb R^m\right)$.
\end{lemma}

$\mathcal D_cV$ is defined in $\Omega\times(0,\infty)$ and extends continuously to $\overline\Omega\times[0,\infty)$, with the extension given by the right-hand side of the equation.

For later use, the linear initial-boundary value problem associated with \eqref{P-draft} is written explicitly as
\begin{equation}\label{EX-linear}
\left\{
\begin{aligned}
&\partial_tV=c(t)\partial_xV+D(x,t)\mathcal K_\Omega[V]+A(x,t)V,
&&x\in\Omega,\ t>0,\\
&V(l_2,t)=\mathbf0,&&t>0,\\
&V(x,0)=V_0(x),&&x\in\overline\Omega.
\end{aligned}
\right.
\end{equation}
Under the same change of variables as above, its moving-interval counterpart is
\begin{equation}\label{IBP-free-linear}
\left\{
\begin{aligned}
&\partial_tU=D^*(x,t)\mathcal J_{(g(t),h(t))}[U]
+A^*(x,t)U,
&&g(t)<x<h(t),\ t>0,\\
&U(h(t),t)=\mathbf0,&&t>0,\\
&U(x,0)=V_0(x),&&l_1\leq x\leq l_2,
\end{aligned}
\right.
\end{equation}
where $D^*(x,t)=D(x-\alpha(t),t)$, $A^*(x,t)=A(x-\alpha(t),t)$, and $\mathcal J_{(g(t),h(t))}$ is the diagonal operator with entries $\mathcal J^i_{(g(t),h(t))}$.

The above arguments also apply, with only minor simplifications,
to the corresponding linear cooperative systems \eqref{EX-linear} and \eqref{IBP-free-linear}. Indeed, if the
nonlinear term $F(x,t,V)$ is replaced by $A(x,t)V$, where
$A(x,t)$ is cooperative, then the parametrized ODE and fixed-point
arguments remain valid, while the linearity directly rules out
finite-time blow-up. The proofs of Lemmas~\ref{linear-advection-existence} and~\ref{EX-free-linear} are therefore omitted.

\begin{lemma}\label{linear-advection-existence}
Assume (H1)--(H4). For every $V_0\in\mathbb U_{(l_1,l_2)}$ satisfying
$ V_0\ge 0$, problem \eqref{EX-linear} admits a unique nonnegative global solution $V(x,t;V_0)$
in the characteristic sense, with
$
V, \mathcal{D}_{c} V \in C\left( \overline\Omega\times[0,\infty);\mathbb R^m \right).
$
If, in addition, $V_0\in\mathbb V_{(l_1,l_2)}$, then
$V\in C^1\left(\overline{\Omega}\times[0,\infty);\mathbb R^m\right)$.
\end{lemma}

\begin{lemma}\label{EX-free-linear}
Assume (H1)--(H4). For every $V_0\in\mathbb U_{(l_1,l_2)}$ satisfying
$ V_0\ge 0$, problem \eqref{IBP-free-linear} has a unique  nonnegative global solution $U\in C(\overline\Omega_\infty^{g,h};\mathbb R^m)$. If, in addition, $V_0\in\mathbb V_{(l_1,l_2)}$, then $U\in C^1(\overline\Omega_\infty^{g,h};\mathbb R^m)$.
\end{lemma}

\begin{remark}{\rm 
The proof of Lemma~\ref{EX-free} gives, for every $T_1>0$ and $V_0\in\mathbb V_{(l_1,l_2)}$ with $\mathbf0\leq V_0\leq\mathbf M$,
\begin{equation*}
    \|U\|_{C^1(\overline{\Omega_{T_1}^{g,h}};\mathbb R^m)}
    \leq C_{T_1}
    \left(1+\|V_0\|_{C^1([l_1,l_2];\mathbb R^m)}\right),
\end{equation*}
where $U$ is the solution of \eqref{IBP-free}.
The same argument applied to \eqref{IBP-free-linear} yields
\begin{equation*}
    \|U\|_{C^1(\overline{\Omega_{T_1}^{g,h}};\mathbb R^m)}
    \leq C_{T_1}\|V_0\|_{C^1([l_1,l_2];\mathbb R^m)}
\end{equation*}
for every $V_0\in\mathbb V_{(l_1,l_2)}$. Here $C_{T_1}>0$ may differ in the two estimates and depends
only on $T_1$ and the coefficients of the corresponding problem,
but not on $V_0$.}
\end{remark}

\begin{remark}\label{Well defined}{\rm 
Let $K=\{\phi\in\mathbb V_{(l_1,l_2)}:\phi\geq\mathbf0\}$ be the positive cone of $\mathbb V_{(l_1,l_2)}$.
For $V_0\in K$, set
\begin{equation*}
    \mathcal PV_0:=V(\cdot,T;V_0)
    =U(\cdot+\alpha(T),T;V_0),
\end{equation*}
where $V$ and $U$ are the solutions of \eqref{EX-linear}
and \eqref{IBP-free-linear} with initial data $V_0$.
Nonnegativity of solutions and the boundary compatibility
argument in Step~4 of the proof of Lemma~\ref{EX-free},
together with the $T$-periodicity of the coefficients,
ensure that $\mathcal P(K)\subset K$.

Since the equation is linear and solutions are unique, $\mathcal P$ is additive and positively homogeneous on $K$.
Since $K$ has nonempty interior in $\mathbb V_{(l_1,l_2)}$,
as proved in the Appendix, $K-K=\mathbb V_{(l_1,l_2)}$.
Thus $\mathcal P$ extends uniquely to a linear operator
\begin{equation}\label{Poincare-map}
    \mathcal P:\mathbb V_{(l_1,l_2)}
    \longrightarrow\mathbb V_{(l_1,l_2)}.
\end{equation}
Moreover, the corresponding
$C^1$ estimate gives
\[
\|\mathcal PV_0\|_{C^1([l_1,l_2];\mathbb R^m)}
\le C_T\|V_0\|_{C^1([l_1,l_2];\mathbb R^m)},
\]
where $C_T>0$ depends only on $T$ and the given coefficients.
Hence $\mathcal P$ is a bounded linear operator, 
which we call the Poincar\'e map.}
\end{remark}

\subsection{Maximum principles and boundary estimates for periodic problems}

Let  $c_{\min}=\min_{t\in[0,T]}c(t)>0$ and $c_{\max}=\max_{t\in[0,T]}c(t)$.
\begin{lemma}[Scalar strong maximum principle]
\label{Max-sca}
Assume (H1). Let $d$ be a continuous strictly positive coefficient on $\overline Q_T$, and let $k$ be a nonnegative continuous $T$-periodic kernel with $k(0,t)>0$ and $\int_\mathbb R k(z,t)\,dz=1$ for all $t$.
Let $b\in L^\infty(Q_T)$, and suppose that $u\in C(\overline Q_T;[0,\infty))$ with $\mathcal D_cu\in C(\Omega\times(0,T);\mathbb R)$, where $\mathcal D_c=\partial_t-c(t)\partial_x$.    
Suppose that
\begin{equation*}
\mathcal D_cu(x,t)
\geq d(x,t)\int_{l_1}^{l_2}k(x-y,t)u(y,t)\,dy
+b(x,t)u(x,t)
\enspace \text{in}\enspace  \Omega\times(0,T),
\end{equation*}
and
$u(x,0)\geq u(x,T)$, $x\in\overline\Omega$.
Then either $u\equiv0$ in $\overline Q_T$, or $u(x,t)>0$ in $[l_1,l_2)\times[0,T]$.
\end{lemma}
\begin{proof}
Assume that there exists
$(x_0,t_0)\in[l_1,l_2)\times[0,T]$
such that $u(x_0,t_0)=0$.
If $t_0=0$, the endpoint inequality and nonnegativity give
\[
0=u(x_0,0)\geq u(x_0,T)\geq0.
\]
Thus, replacing $t_0$ by $T$ if necessary, we may assume
that $t_0\in(0,T]$.

Set
\begin{equation*}
\beta=\|b\|_{L^\infty(\Omega\times(0,T))},
\qquad
\mathcal N[u](x,t)
=
\int_{l_1}^{l_2}k(x-y,t)u(y,t)\,dy.
\end{equation*}
By the kernel assumptions, there exist $\rho_0,\kappa_0>0$
such that
$k(\xi,t)\geq\kappa_0$ for $|\xi|<\rho_0$, $t\in[0,T]$.
Since $b\geq-\beta$ almost everywhere, we have
\begin{equation*}
\mathcal D_cu+\beta u-d\mathcal N[u]\geq0
\enspace \text{almost everywhere in} \enspace \Omega\times(0,T).
\end{equation*}
The left-hand side is continuous, so this inequality holds
everywhere in $\Omega\times(0,T)$.

\textbf{Step 1.} A zero at one spatial point forces $u(\cdot,t_0)\equiv0$.
Define
\begin{equation*}
Z_{t_0}:=
\{x\in[l_1,l_2):u(x,t_0)=0\}.
\end{equation*}
This set is nonempty because $x_0\in Z_{t_0}$,
and it is relatively closed in $[l_1,l_2)$ by continuity.
We prove that it is also relatively open.

Fix $z\in Z_{t_0}$ and choose $a\in(0,t_0)$ sufficiently close
to $t_0$ that
\begin{equation*}
X(r):=z+\int_r^{t_0}c(s)\,ds\in\Omega,
\qquad a\leq r<t_0.
\end{equation*}
The function $v(r):=u(X(r),r)$ is continuous on $[a,t_0]$
and continuously differentiable on $(a,t_0)$, with
\begin{equation*}
v'(r)+\beta v(r)
\geq d(X(r),r) \mathcal{N}[u](X(r),r)\geq0,
\qquad a<r<t_0.
\end{equation*}
Hence, for $a<r<s<t_0$,
\begin{equation*}
0\leq e^{\beta r}v(r)\leq e^{\beta s}v(s).
\end{equation*}
Letting $s\to t_0$ and using $v(t_0)=0$, we obtain
$v\equiv0$ on $[a,t_0]$.
Since $d>0$ in $Q_T$, the preceding differential inequality gives
$\mathcal N[u](X(r),r)=0, a<r<t_0$.
Letting $r\to t_0$, continuity yields
\begin{equation*}
0=\mathcal N[u](z,t_0)
\geq\kappa_0
\int_{\Omega\cap(z-\rho_0,z+\rho_0)}u(y,t_0)\,dy
\geq0.
\end{equation*}
Thus $u(y,t_0)=0$ for every
$y\in[l_1,l_2)\cap(z-\rho_0,z+\rho_0)$.
The set $Z_{t_0}$ is relatively open, so connectedness gives
$Z_{t_0}=[l_1,l_2)$.
Continuity at $l_2$ then implies
$u(\cdot,t_0)\equiv0$ on $\overline\Omega$.

The same argument shows that, for any $s\in(0,T]$, a zero at some $z\in[l_1,l_2)$ implies $u(\cdot,s)\equiv0$ on $\overline\Omega$.

\textbf{Step 2.} Propagation of the zero state in time.
Define
\begin{equation*}
s_*:=\inf\left\{
s\in[0,t_0]:
u\equiv0\text{ on }\overline\Omega\times[s,t_0]
\right\}.
\end{equation*}
The set is nonempty by Step 1.
Continuity gives
$u\equiv0$ on $\overline\Omega\times[s_*,t_0]$.
Suppose that $s_*>0$.
Fix $x_*\in\Omega$ and choose $h\in(0,s_*)$ small enough that
\begin{equation*}
X(r;x_*,s_*):=x_*+\int_r^{s_*}c(\xi)\,d\xi\in\Omega,
\qquad s_*-h\leq r\leq s_*.
\end{equation*}
Since $u(x_*,s_*)=0$, the characteristic argument in Step 1 gives
\begin{equation*}
u(X(r;x_*,s_*),r)=0,
\qquad s_*-h\leq r\leq s_*.
\end{equation*}
Applying the result from Step 1 at each
such $r$ yields $u(\cdot,r)\equiv0$ on $\overline\Omega$.
Thus $u\equiv0$ on $\overline\Omega\times[s_*-h,t_0]$,
contrary to the definition of $s_*$.
Hence $s_*=0$, and
\begin{equation*}
u\equiv0
\enspace \text{in} \enspace \overline\Omega\times[0,t_0].
\end{equation*}
Since $u(\cdot,0)\equiv0$,
$0\leq u(x,T)\leq u(x,0)=0$, $x\in\overline\Omega$.
Repeating the preceding argument with
$T$ in place of $t_0$, we obtain
$u\equiv0$ on $\overline Q_T$.
This completes the proof.
\end{proof}

\begin{lemma}[Cooperative strong maximum principle]\label{max-coo}
Assume (H1)--(H3).
Let $A=(a_{ik})\in C(\overline Q_T;\mathbb R^{m\times m})$
be cooperative, and assume that $A(x_0,t_*)$ is irreducible
for some $(x_0,t_*)\in \Omega\times(0,T)$.
Let $W\in C(\overline Q_T;\mathbb R^m)$ be nonnegative, with
$\mathcal D_cW\in C(\Omega\times(0,T);\mathbb R^m)$, and suppose that
$W(x,0)\geq W(x,T)$,
$x\in\overline\Omega$.
Assume that, for every $i\in\mathbb S$,
\begin{equation*}
\begin{aligned}
\mathcal D_cW_i(x,t)
\geq &
d_i(x,t) \mathcal{K}_{(l_{1},l_{2})}^{i}[W_{i}](x,t) + \sum_{k=1}^m a_{ik}(x,t)W_k(x,t)
\enspace \text{in} \enspace \Omega\times(0,T).
\end{aligned}
\end{equation*}
Then either $W\equiv\mathbf0$ on $\overline Q_T$, or
$W(x,t)\gg\mathbf0$,
$(x,t)\in [l_1,l_2)\times[0,T]$.
\end{lemma}

\begin{proof}
Cooperativity and $W\geq\mathbf0$ give
\begin{equation*}
\mathcal D_cW_i
\geq d_i(x,t)\int_{l_1}^{l_2}k_i(x-y,t)W_i(y,t)\,dy
+(a_{ii}-d_i)W_i
\enspace \text{in} \enspace \Omega\times(0,T).
\end{equation*}
Since $a_{ii}-d_i$ is bounded and
$W_i(\cdot,0)\geq W_i(\cdot,T)$, by
Lemma~\ref{Max-sca}, for every $i\in\mathbb S$,
either $W_i\equiv0$ on $\overline Q_T$, or
$W_i(x,t)>0$ in $[l_1,l_2)\times[0,T]$.
Set
\begin{equation*}
I_0=\{i\in\mathbb S:W_i\equiv0\text{ on }\overline Q_T\}.
\end{equation*}
If $I_0=\mathbb S$ or $I_0=\varnothing$, the conclusion follows.
Suppose that $I_0$ is a nonempty proper subset of $\mathbb S$.
For $i\in I_0$, substitution into the differential inequality gives
\begin{equation*}
0\geq\sum_{k\notin I_0}a_{ik}(x,t)W_k(x,t)\geq0
\text{ in } \Omega\times(0,T).
\end{equation*}
The irreducibility of $A(x_0,t_*)$ provides
$i_*\in I_0$ and $k_*\notin I_0$ with
$a_{i_*k_*}(x_0,t_*)>0$.
Since $W_{k_*}(x_0,t_*)>0$, we obtain
\begin{equation*}
0\geq\sum_{k\notin I_0}a_{i_*k}(x_0,t_*)W_k(x_0,t_*)
\geq a_{i_*k_*}(x_0,t_*)W_{k_*}(x_0,t_*)>0,
\end{equation*}
a contradiction.
Hence $I_0$ is empty or equal to $\mathbb S$.
\end{proof}

For later use, when $c(t)>0$ define the backward characteristic through $(x,t)$ by
\[
X(r;x,t)=x+\int_r^t c(s)\,ds,
\]
and let $\tau(x,t)<t$ be the unique time satisfying
$X(\tau(x,t);x,t)=l_2$. By $T$-periodicity and the positive lower bound on $c$, this time is well defined for every $(x,t)\in Q$.

\begin{lemma}[Upper boundary estimate]
Assume (H1)--(H3).
\label{LEM-upper-boundary}
Let $H\in C\left(
\overline\Omega\times\mathbb R\times\mathbb R_+^m;
\mathbb R^m
\right)$ be $T$-periodic in $t$.
Assume that $V\in C(\overline\Omega\times\mathbb R;\mathbb R_+^m)$ is
 $T$-periodic in $t$, with $\mathcal D_cV\in C(Q;\mathbb R^m)$,
and satisfies
\begin{equation*}
\left\{
\begin{aligned}
&\mathcal D_cV(x,t)
\leq
D(x,t)\mathcal K_\Omega[V](x,t)
+
H(x,t,V(x,t)),
&& (x,t)\in Q_T,\\
&V(l_{2},t)
= \mathbf0,
&& t\in[0,T].
\end{aligned}
\right.
\end{equation*}
Then, for each $i\in\mathbb S$, there exists a constant $\overline C_i>0$ such that
$V_i(x,t)\leq\overline C_i|x-l_{2}|$
for all $(x,t)\in\overline\Omega \times \mathbb{R}$.

\end{lemma}

\begin{proof}
Fix $i\in\mathbb S$.
By continuity and periodicity, choose $B_i>0$ such that
\begin{equation*}
d_i(x,t)\mathcal K_\Omega^i[V_i](x,t)
+H_i(x,t,V(x,t))\le B_i \text{ in }\overline\Omega\times\mathbb R.
\end{equation*}
Fix $(x,t)\in Q$, and write $X(r)=X(r;x,t)$ and
$\tau=\tau(x,t)$.
For $\tau<\sigma<t$, integration along the interior
characteristic segment gives
\begin{equation*}
\begin{aligned}
V_i(x,t)-V_i(X(\sigma),\sigma)
=\int_\sigma^t\mathcal D_cV_i(X(r),r)\,dr \leq B_i(t-\sigma).
\end{aligned}
\end{equation*}
Letting $\sigma\to\tau$ and using
$V_i(l_{2},\tau)=0$, we obtain
\begin{equation*}
V_i(x,t)\leq B_i(t-\tau)
\leq\frac{B_i}{ c_{\min} }|x-l_{2}|
\enspace \text{for all} \enspace (x,t)\in Q.
\end{equation*}
The estimate extends to the spatial boundary by continuity.
\end{proof}

\begin{lemma}[Lower boundary estimate]
Assume (H1)--(H3).
\label{LEM-lower-boundary}
Let $H\in C\left(
\overline\Omega\times\mathbb R\times\mathbb R_+^m;
\mathbb R^m
\right)$ be $T$-periodic in $t$ and 
$H_i(l_2,t,V(l_2,t))\geq0$
whenever $V_i(l_2,t)=0$.
Assume that $V\in C(\overline\Omega\times\mathbb R;\mathbb R_+^m)$ is $T$-periodic in $t$, with $\mathcal D_cV\in C(Q;\mathbb R^m)$,
and satisfies
\begin{equation*}
\left\{
\begin{aligned}
&\mathcal D_cV(x,t)
\geq
D(x,t)\mathcal K_\Omega[V](x,t)
+
H(x,t,V(x,t)),
&& (x,t)\in Q_T,\\
&V(l_{2},t)
\ge \mathbf0,
&& t\in[0,T].
\end{aligned}
\right.
\end{equation*}
If $V_i>0$ on $(\overline\Omega\setminus\{l_2\})\times\mathbb R$ for every $i\in\mathbb S$, then, for each $i$, there exists
$\underline C_i>0$ such that
\[
V_i(x,t)\geq\underline C_i |x-l_2|
\quad\text{for all }(x,t)\in\overline Q_T.
\]
\end{lemma}

\begin{proof}
By periodicity, the differential inequality holds in $Q$.
Fix $i\in\mathbb S$ and set
\begin{equation*}
\mathcal R_i(x,t)
:=d_i(x,t)\mathcal K_\Omega^i[V_i](x,t)
+H_i(x,t,V(x,t)).
\end{equation*}
This function is continuous on $\overline\Omega\times\mathbb R$
and $T$-periodic in $t$.
If $V_i(l_2,t)=0$, we have
\begin{equation*}
\begin{aligned}
\mathcal R_i(l_2,t)
&=d_i(l_2,t)\int_{l_1}^{l_2}
k_i(l_2-y,t)V_i(y,t)\,dy
+H_i(l_2,t,V(l_2,t))>0.
\end{aligned}
\end{equation*}
We claim that there exists $\gamma_i>0$ such that
\begin{equation*}
2\beta_i:=\min_{t\in[0,T]}
\left[\mathcal R_i(l_2,t)+\gamma_iV_i(l_2,t)\right]>0.
\end{equation*}
Indeed, otherwise there would be $t_n\in[0,T]$ with
$\mathcal R_i(l_2,t_n)+nV_i(l_2,t_n)\leq0$.
Boundedness of $\mathcal R_i(l_2,\cdot)$ gives
$V_i(l_2,t_n)\to0$.
Passing to a subsequence with $t_n\to t_*$, we obtain
$V_i(l_2,t_*)=0$ and $\mathcal R_i(l_2,t_*)\leq0$,
a contradiction.
By continuity and periodicity, choose small
$\rho_i\in(0,(l_2-l_1)/2)$ such that
\begin{equation*}
\mathcal R_i(x,t)+\gamma_iV_i(x,t)\geq\beta_i
\text{ in }[l_2-\rho_i,l_2]\times\mathbb R.
\end{equation*}
Thus
\begin{equation*}
\mathcal D_cV_i+\gamma_iV_i\geq\beta_i
\text{ in }(l_2-\rho_i,l_2)\times\mathbb R.
\end{equation*}

Fix $x\in[l_2-\rho_i,l_2)$ and $t\in\mathbb R$.
Write $X(r)=x+\int_r^t c(s)\,ds$ and let
$\tau=\tau(x,t)$ satisfy $X(\tau)=l_2$.
For $\tau\leq r\leq t$, we have $x\leq X(r)\leq l_2$.
Moreover,
\begin{equation*}
l_2-x=\int_\tau^t c(r)\,dr,
\qquad
\frac{l_2-x}{c_{\max}}\leq t-\tau
\leq\frac{l_2-x}{c_{\min}}\leq\frac{\rho_i}{c_{\min}}.
\end{equation*}
For $\tau<\sigma<t$, integration along the interior
characteristic segment gives
\begin{equation*}
\begin{aligned}
V_i(x,t)
=e^{-\gamma_i(t-\sigma)}V_i(X(\sigma),\sigma)+\int_\sigma^t e^{-\gamma_i(t-r)}
(\mathcal D_cV_i+\gamma_iV_i)(X(r),r)\,dr.
\end{aligned}
\end{equation*}
Letting
$\sigma\to \tau$, we obtain
\begin{equation*}
\begin{aligned}
V_i(x,t)
&\geq e^{-\gamma_i(t-\tau)}V_i(l_2,\tau)
+\beta_i\int_\tau^t e^{-\gamma_i(t-r)}\,dr\\
&\geq\beta_i e^{-\gamma_i\rho_i/c_{\min}}(t-\tau)\\
&\geq\frac{\beta_i}{c_{\max}}e^{-\gamma_i\rho_i/c_{\min}}(l_2-x)
\end{aligned}
\end{equation*}
for all $(x,t) \in [l_2-\rho_i,l_2)\times \mathbb{R}$. Moreover, the estimate holds at $x=l_2$ as well by the continuity of $V_i$.

On the remaining compact set, strict positivity gives
$m_i:=\min_{[l_1,l_2-\rho_i]\times[0,T]}V_i>0$.
Since $l_2-x\leq l_2-l_1$, we have
\begin{equation*}
V_i(x,t)\geq m_i
\geq\frac{m_i}{l_2-l_1}(l_2-x)
\quad\text{on }[l_1,l_2-\rho_i]\times[0,T].
\end{equation*}
The conclusion follows with
$\underline C_i
:=\min\left\{
\frac{\beta_i}{c_{\max}}e^{-\gamma_i\rho_i/c_{\min}},
\frac{m_i}{l_2-l_1}
\right\}>0$.
\end{proof}

\begin{remark}\label{rem:boundary-estimates-applications}
The estimates in Lemmas~\ref{LEM-upper-boundary}
and~\ref{LEM-lower-boundary} apply both to the nonlinear
problem with
\[
H(x,t,u)=F(x,t,u),
\]
where $F$ satisfies \textup{(F)}, and to the linear eigenvalue
problem with
\[
H(x,t,u)=A(x,t)u-\lambda u,
\]
where the matrix $A$ satisfies \textup{(H4)}.
\end{remark}

\section{Principal spectral theory}
\subsection{Existence and uniqueness of the principal eigenvalue}

This subsection proves the existence and uniqueness of the
principal eigenvalue in Theorem~\ref{EX-PEV}(i). The main step is to
establish the strong positivity and power compactness of the
Poincar\'e map $\mathcal P$ defined in \eqref{Poincare-map}.

Recall that
\begin{equation*}
K=\{\phi\in\mathbb V_{(l_1,l_2)}:\phi\geq\mathbf0\}
\end{equation*}
is the positive cone of the Banach space $\mathbb V_{(l_1,l_2)}$, whose interior is characterized in Proposition~\ref{cone-interior}.
A bounded linear
operator $L$ is called power compact if $L^n$ is compact for
some positive integer $n$. The following result for power compact operators
can be found in \cite[Lemma~3.1]{LiangZhao}
and \cite[Section~7.1]{AA}.
\begin{lemma}\label{power-positive}
Let $X$ be a Banach space ordered by a closed positive cone
$K$ with nonempty interior. Suppose that $L:X\to X$ is a
bounded positive linear operator and that, for some $n\geq1$,
$L^n$ is compact and strongly positive, namely,
$L^n(K\setminus\{0\})\subset\operatorname{int}K$.
Then $r(L)>0$ is an algebraically simple eigenvalue
with an eigenvector in $\operatorname{int}K$,
and every other spectral value $\mu$ satisfies
$|\mu|<r(L)$.
\end{lemma}

The next maximum principle is used to establish the
strong positivity of $\mathcal P$. As in Section~2,
$g(t)=l_1+\alpha(t)$ and $h(t)=l_2+\alpha(t)$.
\begin{lemma}\label{max-ini-1}
 Assume (H1)--(H3) and let $T_1>0$.
Let $A^*=(a_{ij}^*)_{i,j\in\mathbb S}\in
C(\overline{\Omega}_{T_1}^{g,h};\mathbb R^{m\times m})$
be cooperative, and suppose that there exists $x_0\in\Omega$
such that $A^*(x_0+\alpha(t),t)$ is irreducible
for every $t\in[0,T_1]$.
Suppose that
$U,\partial_tU\in
C(\overline{\Omega_{T_1}^{g,h}};\mathbb R^m)$
and that the following inequalities hold for every $i\in\mathbb S$:
	\begin{equation}\label{IBP-U}
		\left\{
		\begin{aligned}
			& \partial_{t} U_{i}(x,t) \ge d_{i}^{*}(x,t)\mathcal{J}^{i}_{(g(t),h(t))}[U_{i}](x,t) + \sum_{j=1}^{m} a_{ij}^{*}(x,t)U_{j}(x,t),  &&(x,t)\in {\Omega}_{T_1}^{g,h},\\
			&U_{i}(h(t),t) \ge 0, && 0<t\le T_1,\\
			&U_{i}(x,0) \ge 0, &&x\in[l_{1},l_{2}].
		\end{aligned}
		\right.
	\end{equation}
	Then $U\ge\mathbf0$ on $\overline\Omega_{T_1}^{g,h}$.
		If $U(\cdot,0)\not\equiv\mathbf0$, then  $U\gg \bm 0$ in $\{ (x,t): 0<t\le T_{1},\,g(t) \le x < h(t)  \}$.
\end{lemma}
\begin{proof}
	
	By an argument similar to that in
	\cite[Lemma~3.1]{DuNi}, it follows that
	$U\ge\mathbf0$ on $\overline\Omega_{T_1}^{g,h}$.
	Moreover, if $U_{i}(\cdot,0)\not\equiv\mathbf0$,
	then $U_{i} > 0$ in $\Omega_{T_1}^{g,h}$.
	It remains to establish positivity at the outflow boundary and to show that irreducible coupling makes every component strictly positive, including those with identically zero initial data.

	Set
	\begin{equation*}
		\hat\Omega_{T_1}^{g,h}
		=\{(x,t):0<t\le T_1,\,
		g(t)\le x<h(t)\}.
	\end{equation*}
	
		First, suppose that
		$U_i(\cdot,0)\not\equiv0$ on $[l_1,l_2]$.
		It will be shown that
		$U_i(g(t),t)>0$ for all $0<t\le T_1$.
		Suppose, to the contrary, that there exists
		$t_1\in(0,T_1]$ such that
		$U_i(g(t_1),t_1)=0$.
	Since there exists $\delta>0$ such that $(g(t_{1}),t)\in \Omega_{T_1}^{g,h}$ for $t\in (t_{1}-\delta,t_{1})$ and $U_{i}(g(t_{1}),t)>0$ for $t\in (t_{1}-\delta,t_{1})$, it follows that $\partial_{t}U_{i}(g(t_{1}), t_{1}) \le 0$.
	It then follows from the first inequality in \eqref{IBP-U} that
	\begin{equation*}
		\begin{aligned}
			0 \ge \lim_{x\to g(t_{1})} \partial_{t}U_{i}(x,t_{1}) \ge& 
			\lim_{x\to g(t_{1})} \left( d_{i}^{*}(x,t_{1})\mathcal{J}^{i}_{(g(t_{1}),h(t_{1}))}[U_{i}](x,t_{1}) + \sum_{j=1}^{m} a_{ij}^{*}(x,t_{1})U_{j}(x,t_{1}) \right)
			\\
			\ge&\, d_{i}^{*}(g(t_{1}),t_{1})\int_{g(t_{1})}^{h(t_{1})}k_{i}(g(t_{1})-y,t_{1})U_{i}(y,t_{1})dy>0,
		\end{aligned}
	\end{equation*}
	which is a contradiction. 
	
	If $U(x,0)\not\equiv \bm 0$ in $[l_{1},l_{2}]$, then there exists $i_{0}\in\mathbb{S}$ such that $U_{i_{0}}(\cdot,0)\not\equiv 0$ on $[l_{1},l_{2}]$. Hence,
	\begin{equation*}
		\mathbb{I}:= \{ i\in \mathbb{S} : U_{i}>0 \text{ in } \hat\Omega_{T_1}^{g,h}\} \neq \emptyset.
	\end{equation*}
	If $\mathbb{I}=\mathbb{S}$, then $U_{i}>0 \text{ in } \hat\Omega_{T_1}^{g,h}$ for $i \in \mathbb{S}$. The proof is complete. Assume that $\mathbb{J}=\mathbb{S}\backslash \mathbb{I} \neq \emptyset$.
	Then for any $j\in \mathbb{J}$
	there exists $(x_{j},t_{j})\in \hat\Omega_{T_1}^{g,h}$ such that $U_{j}(x_{j},t_{j})=0$. Since $U\ge\mathbf0$ on
		$\overline\Omega_{T_1}^{g,h}$, it follows that
	\begin{equation*}
		\begin{aligned}
			0 \ge \lim_{x\to x_{j}} \partial_{t}U_{j}(x,t_{j}) \ge& 
			\lim_{x\to x_{j}} \left( d_{j}^{*}(x,t_{j})\mathcal{J}^{j}_{(g(t_{j}),h(t_{j}))}[U_{j}](x,t_{j}) + \sum_{k=1}^{m} a_{jk}^{*}(x,t_{j})U_{k}(x,t_{j}) \right)
			\\
			\ge&\, d_{j}^{*}(x_{j},t_{j})\int_{g(t_{j})}^{h(t_{j})}k_{j}(x_{j}-y,t_{j})U_{j}(y,t_{j})dy \ge 0.
		\end{aligned}
	\end{equation*}
	Since $d_j^*(x_j,t_j)>0$, it follows that
	\begin{equation*}
		\int_{g(t_j)}^{h(t_j)}
		k_j(x_j-y,t_j)U_j(y,t_j)\,dy=0.
	\end{equation*}
		To show that $U_j(\cdot,t_j)\equiv0$,
		define $Z_j=\{x\in[g(t_j),h(t_j)):
		U_j(x,t_j)=0\}.$
		Clearly, $Z_j$ is nonempty and relatively
		closed in $[g(t_j),h(t_j))$. For any $x\in Z_j$, the nonnegativity of
		$U_j$ implies that
		$\partial_tU_j(x,t_j)\le0$.
		At $x=g(t_j)$, this follows by taking the
		time derivative from the left, since
		$g(t)$ is strictly increasing.
		The supersolution inequality then yields
	\begin{equation*}
		0\ge d_j^*(x,t_j)
		\int_{g(t_j)}^{h(t_j)}
		k_j(x-y,t_j)U_j(y,t_j)\,dy
		\ge0.
	\end{equation*}
		By (H3), there exists $\delta>0$ such that
		$k_j(z,t_j)>0$ for $|z|<\delta$.
		Consequently,
		$U_j(y,t_j)=0$ whenever
		$y\in[g(t_j),h(t_j))$ and $|y-x|<\delta$.
		Thus, $Z_j$ is also relatively open in
		$[g(t_j),h(t_j))$. Since $[g(t_j),h(t_j))$ is connected,
		$Z_j=[g(t_j),h(t_j))$.
		By continuity,
	\begin{equation*}
		U_j(x,t_j)=0,
		\quad x\in[g(t_j),h(t_j)].
	\end{equation*}
	It then follows from \cite[Lemma 3.1]{DuNi} that 
	\begin{equation*}
		U_{j}(x,t)=0 \enspace\text{for}\enspace j\in\mathbb{J} \enspace\text{and}\enspace (x,t)\in \Omega^{g,h}_{t_{j}}.
	\end{equation*}
	Set $t_0=\min_{j\in\mathbb J}t_j$. It follows that
	\begin{equation*}
		U_{j}(x,t)=0 \enspace\text{for}\enspace j\in\mathbb{J} \enspace\text{and}\enspace (x,t)\in \Omega^{g,h}_{t_{0}}.
	\end{equation*}
	For each $j\in\mathbb J$ and $0<t<t_0$, the supersolution inequality gives
	\begin{equation}
		\begin{aligned}
			0 =  \partial_{t}U_{j}(x,t) \ge&\, d_{j}^{*}(x,t)\mathcal{J}^{j}_{(g(t),h(t))}[U_{j}](x,t) + \sum_{i=1}^{m} a_{ji}^{*}(x,t)U_{i}(x,t)\\ 
			\ge& \sum_{i \in \mathbb{I}} a_{ji}^{*}(x,t)U_{i} (x,t) \ge0 \enspace \text{for} \enspace (x,t)\in \Omega^{g,h}_{t_{0}}.
		\end{aligned} 
	\end{equation}
	It follows that
	\begin{equation*}
		a_{ji}^{*}(x,t)=0 \enspace\text{for}\enspace i\in \mathbb{I},j\in\mathbb{J} \text{ and }(x,t)\in \Omega_{t_0}^{g,h}.
	\end{equation*}
	For any $s\in(0,t_0)$, these identities imply that
	$A^*(x_0+\alpha(s),s)$ is reducible,
    contradicting the irreducibility of
    $A^*(x_0+\alpha(s),s)$.
	Hence, $\mathbb I=\mathbb S$.
\end{proof}

    Returning to the fixed interval through
	$U(x,t)=V(x-\alpha(t),t)$ gives the following
	maximum principle.

\begin{lemma}\label{max-ini-2}
Assume (H1)--(H3).
Let $A=(a_{ij})_{i,j\in\mathbb S}
\in C(\overline Q_T;\mathbb R^{m\times m})$
be cooperative, and suppose that there exists $x_0\in\Omega$
such that $A(x_0,t)$ is irreducible for every $t\in[0,T]$.
Suppose that $V,\mathcal D_cV\in
C(\overline Q_T;\mathbb R^m)$ satisfy
	\begin{equation}
		\left\{
		\begin{aligned}
			& \mathcal D_c V(x,t) \ge D(x,t)\mathcal{K}_{(l_{1},l_{2})}[V](x,t) + A(x,t)V(x,t),  &&(x,t)\in Q_{T},\\
			&V(l_{2},t) \ge 0, &&t \in (0,T],\\
			&V(x,0) \ge 0, &&x\in [l_{1},l_{2}].
		\end{aligned}
		\right.
	\end{equation}
	Then $V\ge \bm{0}$ on $\overline{Q}_{T}$. Moreover, if $V(\cdot,0)\not\equiv 0$ on $[l_{1},l_{2}]$, then $V(x,t)\gg\mathbf0$ for  $ (x,t)\in[l_1,l_2)\times(0,T]$. 
\end{lemma}

\begin{theorem}\label{thm:poincare-spectral-properties}
Assume (H1)--(H4). Then the Poincar\'e map
$\mathcal P:\mathbb V_{(l_1,l_2)}\to\mathbb V_{(l_1,l_2)}$
is strongly positive with respect to $K$ and is power compact.
\end{theorem}
\begin{proof}
By Lemma~\ref{max-ini-2}, every
$V_0\in K\setminus\{0\}$ satisfies
$V_i(x,T;V_0)>0$ for all $x\in[l_1,l_2)$ and 
$i\in\mathbb S$.
Proposition~\ref{cone-interior} therefore gives
$\mathcal PV_0\in\operatorname{int}K$.
This proves strong positivity. It remains to prove power compactness.

	\textbf{Step 1. Integral representation of $\mathcal P^N$.}

	For $U=(U_1,\ldots,U_m)$, define the vector function
	$\mathcal M_{(g(t),h(t))}[U](x,t)$
	componentwise by
\begin{equation*}
	\bigl(\mathcal M_{(g(t),h(t))}[U](x,t)\bigr)_i
	=d_i^*(x,t)\int_{g(t)}^{h(t)}
	k_i(x-y,t)U_i(y,t)\,dy,
	\qquad i\in\mathbb S.
\end{equation*}
	Set $\tilde B(x,t):=A^*(x,t)-D^*(x,t)$.
	Then $U(x,t;V_0)$ satisfies
	\begin{equation}\label{M-B}
		\left\{
		\begin{aligned}
			& \partial_{t} U(x,t) =  \mathcal{M}_{(g(t),h(t))}[U](x,t) + \tilde{B}(x,t)U(x,t), &&x\in(g(t),h(t)),t>0,\\
			&U(h(t),t)=0,  &&t>0, \\
			&U(x,0)=V_{0}(x),&&x\in[l_{1},l_{2}].
		\end{aligned}
		\right.
	\end{equation}
	
	As in the proof of Lemma~\ref{EX-free},
	\eqref{M-B} can be viewed as an ODE initial value
	problem with parameter $x$.
	Since $g'(t)=h'(t)=c(t)\ge c_0>0$ and
	$g(t),h(t)\to\infty$ as $t\to\infty$,
		the functions $g$ and $h$ have continuously
		differentiable inverses
		\[
		I_g:[l_1,\infty)\to[0,\infty),
		\qquad
		I_h:[l_2,\infty)\to[0,\infty),
		\]
		respectively. More precisely,
		$I_g(x)$ and $I_h(x)$ are uniquely determined by
		\[
		g(I_g(x))=x \quad (x\ge l_1),
		\qquad
		h(I_h(x))=x \quad (x\ge l_2).
		\]
	For $x\ge l_1$, define
	\begin{equation*}
		t_{0}(x) =
		\left\{
		\begin{aligned}
			& 0,
			&& x\in [l_{1},l_{2}],\\
			& I_{h}(x), && x>l_{2},
		\end{aligned}
		\right.
		\quad
		\tilde{V}_{0}(x) =
		\left\{
		\begin{aligned}
			& V_{0}(x),
			&& x\in [l_{1},l_{2}],\\
			& 0, && x>l_{2}.
		\end{aligned}
		\right.
	\end{equation*}
	For each $x>l_1$, the function $t\mapsto U(x,t)$ satisfies
	\begin{equation}\label{ODE-M-B}
		\left\{
		\begin{aligned}
			&\partial_tU(x,t)
			=\mathcal M_{(g(t),h(t))}[U](x,t)+\tilde B(x,t)U(x,t),
			&&t_0(x)<t<I_g(x),\\
			&U(x,t_0(x))=\tilde V_0(x).
		\end{aligned}
		\right.
	\end{equation}

	For $x>l_2$ and $t_0(x)\leq s\leq t\leq I_g(x)$, let
	$\Phi^x(t,s)=(\phi^x_{ij}(t,s))_{m\times m}$
	denote the evolution matrix determined by
	\begin{equation*}
		\partial_t\Phi^x(t,s)=\tilde B(x,t)\Phi^x(t,s),
		\qquad \Phi^x(s,s)=I.
	\end{equation*}
Since $\tilde B$ and $\partial_x\tilde B$
are jointly continuous, the parameter-dependence theorem
for linear ODEs implies that both $\Phi^x(t,s)$ and
$\partial_x\Phi^x(t,s)$ depend continuously on $(x,t,s)$.
For $x>l_2$, the starting value is
$U(x,t_0(x))=\mathbf0$, because the point
$(x,t_0(x))$ lies on the moving inflow boundary.
Applying the variation-of-constants formula to
\eqref{ODE-M-B} gives
	\begin{align}\label{3.6b}
		U(x,t) =& \Phi^{x}(t,t_{0}(x))U(x,t_{0}(x)) + \int_{t_{0}(x)}^{t}  \Phi^{x}(t,s)\mathcal{M}_{(g(s),h(s))}[U](x,s) ds \nonumber\\
		=&  \int_{t_{0}(x)}^{t}  \Phi^{x}(t,s)\mathcal{M}_{(g(s),h(s))}[U](x,s) ds,\qquad t_0(x)\leq t\leq I_g(x).
	\end{align}
	
	Let $\left\lfloor r  \right\rfloor$ denote the greatest integer not exceeding $r$.
	Set
	\begin{equation*}
    N^{*}=  \left\lfloor \frac{l_{2}-l_{1}}{\alpha(T)} \right\rfloor +1 
	\end{equation*}
	and fix an integer $N\geq N^*$.
	The choice $N\alpha(T)>l_2-l_1$ implies
	$x+N\alpha(T)>l_2$ for every $x\in[l_1,l_2]$.
    The identity \eqref{3.6b} and uniqueness of solutions give, for $x\in [l_1,l_2]$,
	\begin{align*}
		\mathcal{P}^{N}[V_{0}](x) =& U(x+N\alpha(T),NT;V_{0})\\
		=& \int_{t_{0}(x+N\alpha(T))}^{NT}  \Phi^{x+N\alpha(T)}(NT,s)\mathcal{M}_{(g(s),h(s))}[U](x+N\alpha(T),s) ds.
	\end{align*}
	
	\textbf{Step 2. Equicontinuity of the spatial derivatives of $\mathcal P^NV_0^n$.}

	Let $\{V_0^n\}_{n\geq1}$ be a bounded sequence in
	$\mathbb V_{(l_1,l_2)}$, and write
	$U^n(x,t):=U(x,t;V_0^n)$ and $W_n:=\mathcal P^NV_0^n$.
	By the uniform $C^1$ estimate for linear solutions
established in Section~2, $\{U^n\}_{n\geq1}$ is bounded
in $C^1(\overline{\Omega_{NT}^{g,h}};\mathbb R^m)$,
uniformly in $n$.
Since $W_n(x)=U^n(x+N\alpha(T),NT)$,
it also follows that
	$$\sup_{n\geq1} \|W_n\|_{C^1([l_1,l_2];\mathbb R^m)}<\infty.$$

	Set $\gamma(x):=t_0(x+N\alpha(T))$, which is increasing and belongs
	to $C^1([l_1,l_2])$. For
	\begin{equation*}
		(x,s)\in E_N
		:=\{(x,s):l_1\leq x\leq l_2,\ \gamma(x)\leq s\leq NT\},
	\end{equation*}
	define
	\begin{equation*}
		\begin{aligned}
			\beta_{ij}(x,s)
			&:=\phi^{x+N\alpha(T)}_{ij}(NT,s)d_j^*(x+N\alpha(T),s),\\
			h_{ij}^n(x,s)
			&:=\beta_{ij}(x,s)\int_{g(s)}^{h(s)}
			k_j(x+N\alpha(T)-y,s)U_j^n(y,s)\,dy.
		\end{aligned}
	\end{equation*}
	Then
	\begin{equation*}
		W_{n,i}(x)=\sum_{j=1}^m\int_{\gamma(x)}^{NT}h_{ij}^n(x,s)\,ds.
	\end{equation*}
	Differentiating after the substitution
	$z=x+N\alpha(T)-y$ and using $U_j^n(h(s),s)=0$
    gives
	\begin{equation}\label{h_{ij}1}
		\begin{aligned}
			\partial_xh_{ij}^n(x,s)
			={}&\partial_x\beta_{ij}(x,s)
			\int_{g(s)}^{h(s)}k_j(x+N\alpha(T)-y,s)U_j^n(y,s)\,dy\\
			&+\beta_{ij}(x,s)
			\int_{g(s)}^{h(s)}k_j(x+N\alpha(T)-y,s)
			\partial_yU_j^n(y,s)\,dy\\
			&+\beta_{ij}(x,s)k_j(x+N\alpha(T)-g(s),s)U_j^n(g(s),s).
		\end{aligned}
	\end{equation}
	Differentiating with respect to $x$ gives
	\begin{equation}\label{dPN}
		W_{n,i}'(x)
		=\sum_{j=1}^m\left[
		\int_{\gamma(x)}^{NT}\partial_xh_{ij}^n(x,s)\,ds
		-h_{ij}^n(x,\gamma(x))\gamma'(x)\right].
	\end{equation}
	The two terms on the right-hand side are estimated separately.

	Set $L=l_2-l_1$. In \eqref{h_{ij}1}, the factors involving
	$U^n$ and $\partial_yU^n$ are independent of $x$.
	Their uniform boundedness, together with the uniform continuity
	of $\beta_{ij}$ and $\partial_x\beta_{ij}$ on $E_N$ and of
	$k_j$ on $[-L,L]\times[0,NT]$, gives
	\begin{equation*}
		\begin{aligned}
			&\sup_n\sup_{(x,s)\in E_N}
			\bigl(|h_{ij}^n(x,s)|+|\partial_xh_{ij}^n(x,s)|\bigr)
			\leq C,\\
			&|\partial_xh_{ij}^n(x,s)-\partial_xh_{ij}^n(x',s)|
			\leq\omega(|x-x'|),
			\qquad
			\max\{\gamma(x),\gamma(x')\}\leq s\leq NT,
		\end{aligned}
	\end{equation*}
	where the constant  $C$ and the function $\omega$ are independent of $n,i,j$, and $\omega(\rho)\to0$ as $\rho\to0$.
	
	For $l_1\leq x\leq x'\leq l_2$, split the
difference at the larger lower limit $\gamma(x')$. This gives
	\begin{equation*}
		\begin{aligned}
			&\left|
			\int_{\gamma(x)}^{NT}\partial_xh_{ij}^n(x,s)\,ds
			-\int_{\gamma(x')}^{NT}\partial_xh_{ij}^n(x',s)\,ds
			\right|\\
			\leq&
			\int_{\gamma(x')}^{NT}
			|\partial_xh_{ij}^n(x,s)-\partial_xh_{ij}^n(x',s)|\,ds
			+\int_{\gamma(x)}^{\gamma(x')}
			|\partial_xh_{ij}^n(x,s)|\,ds\\
			\leq&
			NT\omega(|x-x'|)+C|\gamma(x)-\gamma(x')|.
		\end{aligned}
	\end{equation*}
	Since $\gamma$ is uniformly continuous, the integral terms
	in \eqref{dPN} are equicontinuous.

	Since $h(\gamma(x))=x+N\alpha(T)$, the substitution
	$z=h(\gamma(x))-y$ gives
	\begin{equation*}
		h_{ij}^n(x,\gamma(x))
		=\beta_{ij}(x,\gamma(x))\int_0^L
		k_j(z,\gamma(x))
		U_j^n(h(\gamma(x))-z,\gamma(x))\,dz.
	\end{equation*}
	For $s,s'\in[0,NT]$, the uniform $C^1$ bound for $U^n$ gives
	\begin{equation*}
		\begin{aligned}
			&|U_j^n(h(s)-z,s)-U_j^n(h(s')-z,s')|\\
			=&\left|\int_{s'}^s
			\bigl[\partial_tU_j^n+c(r)\partial_xU_j^n\bigr]
			(h(r)-z,r)\,dr\right|
			\leq C|s-s'|,
			\qquad 0\leq z\leq L.
		\end{aligned}
	\end{equation*}
	Together with the uniform continuity of $\beta_{ij}$, $k_j$,
	and $\gamma$, this proves that
	$\{h_{ij}^n(\cdot,\gamma(\cdot))\}_n$ is equicontinuous.
	Since this family is also uniformly bounded and $\gamma'$ is
	continuous, $\{h_{ij}^n(\cdot,\gamma(\cdot))\gamma'(\cdot)\}_n$
	is equicontinuous as well.
	Thus \eqref{dPN} and the preceding estimates imply that $\{W_n'\}_n$ is
	equicontinuous on $[l_1,l_2]$.
	
	\textbf{Step 3. Compactness of $\mathcal P^N$.}
	
	By Step 2, $\{W_n'\}_{n\geq1}$ is uniformly bounded and
	equicontinuous on $[l_1,l_2]$. The Arzel\`a--Ascoli theorem
	therefore gives, after passing to a subsequence, 
	\begin{equation*}
		W_n'\to Z
		\quad\text{uniformly on }[l_1,l_2],
	\end{equation*}
	where $Z\in C([l_1,l_2];\mathbb R^m)$.
	Since $W_n(l_2)=\mathbf0$, it follows that
	$W_n(x)=-\int_x^{l_2}W_n'(r)\,dr$.
	Set $W(x):=-\int_x^{l_2}Z(r)\,dr$. Then $W'=Z$ and
	\begin{equation*}
		\|W_n-W\|_{C^1([l_1,l_2];\mathbb R^m)}
		\leq (1+l_2-l_1)
		\|W_n'-Z\|_{C([l_1,l_2];\mathbb R^m)}
		\longrightarrow0.
	\end{equation*}
	The closedness of $\mathbb V_{(l_1,l_2)}$
    in $C^1([l_1,l_2];\mathbb R^m)$ gives
	$W\in\mathbb V_{(l_1,l_2)}$.
	As the bounded sequence $\{V_0^n\}_{n\geq1}$ was arbitrary,
	$\mathcal P^N$ is compact. Hence $\mathcal P$ is power-compact.
\end{proof}

\begin{proof}[Proof of Theorem~\ref{EX-PEV}(i)]
By Theorem~\ref{thm:poincare-spectral-properties},
$\mathcal P$ is strongly positive and power compact.
Applying Lemma~\ref{power-positive}, we obtain
$r(\mathcal P)>0$ and an eigenvector
$V_0\in\operatorname{int}K$ such that
$\mathcal PV_0=r(\mathcal P)V_0$.
Moreover, $r(\mathcal P)$ is algebraically simple,
and no other spectral value has the same modulus.

Define
\begin{equation*}
\lambda_1=\frac1T\log r(\mathcal P),\qquad
\varphi(x,t)=e^{-\lambda_1t}V(x,t;V_0).
\end{equation*}
Since $V_0\in\operatorname{int}K
\subset\mathbb V_{(l_1,l_2)}$, Lemma~2.3 gives
$\varphi\in
C^1(\overline\Omega\times[0,\infty);\mathbb R^m)$.
By the definition of $\varphi$ and \eqref{EX-linear},
we have
\begin{equation*}
\partial_t\varphi=c(t)\partial_x\varphi
+D(x,t)\mathcal K_\Omega[\varphi]
+A(x,t)\varphi-\lambda_1\varphi.
\end{equation*}
It also satisfies the inflow boundary condition.
The coefficients are $T$-periodic, so uniqueness of solutions
and the definition of $\mathcal P$ imply that
\begin{equation*}
V(x,t+T;V_0)
=V(x,t;\mathcal PV_0)
=r(\mathcal P)V(x,t;V_0).
\end{equation*}
Consequently,
\begin{equation*}
\varphi(x,t+T)
=e^{-\lambda_1(t+T)}r(\mathcal P)V(x,t;V_0)
=\varphi(x,t).
\end{equation*}
The periodicity identity and the $C^1$ regularity for $t\ge0$
ensure that the periodic extension of $\varphi$ belongs to
$C^1(\overline Q;\mathbb R^m)$.
Together with the strong positivity given by Lemma~3.3,
this shows that $(\lambda_1,\varphi)$ is a principal
eigenpair of (1.1).

To prove uniqueness, let
$(\widetilde\lambda,\widetilde\varphi)$ be another principal
eigenpair and put
$\widetilde V_0=\widetilde\varphi(\cdot,0)$.
The $C^1$ regularity and the zero inflow condition give
$\partial_t\widetilde\varphi(l_2,0)=\mathbf0$.
Taking the boundary limit in the eigenvalue equation at $t=0$
yields the compatibility condition (1.9), so
$\widetilde V_0\in\mathbb V_{(l_1,l_2)}$.
Uniqueness of the linear initial-value problem and
$T$-periodicity give
\[
\mathcal P\widetilde V_0
=e^{\widetilde\lambda T}\widetilde V_0.
\]
Since this eigenvalue admits a nonzero nonnegative eigenvector,
it must equal $r(\mathcal P)$, and hence
$\widetilde\lambda=\lambda_1$.
The simplicity of $r(\mathcal P)$ and uniqueness of the linear
initial-value problem imply that $\widetilde\varphi$ is a
positive scalar multiple of $\varphi$.
This proves Theorem~\ref{EX-PEV}(i).
\end{proof}

\subsection{Variational characterizations of the principal eigenvalue}

The principal eigenvalue established above can also
be described by the three variational quantities
$\mu_p$, $\mu_p'$ and $\lambda_p$ defined in
Theorem~\ref{EX-PEV}(ii). Two boundary estimates
will be used in the proof.

\begin{lemma}
	\label{upper-estimate}
	Let $U\in\mathcal D_0$ and suppose that
	$\mathcal LU\geq\eta U$ in $Q$ for some $\eta\in\mathbb R$, where $D_0$ and $\mathcal L$ are defined in \eqref{1.7a} and \eqref{1.8a}, respectively.
	Then, for each $i\in\mathbb S$, there exists $C_i>0$ such that
	\begin{equation*}
		U_i(x,t)\leq C_i(l_2-x),
		\quad (x,t)\in\overline\Omega\times\mathbb R.
	\end{equation*}
\end{lemma}
\begin{proof}
	The conclusion follows from
	Lemma~\ref{LEM-upper-boundary} and Remark~\ref{rem:boundary-estimates-applications}.
\end{proof}

\begin{lemma}
	\label{lower-estimate}
	Let $U\in\mathcal D$ and suppose that
	$\mathcal LU\leq\eta U$ in $Q$ for some $\eta\in\mathbb R$.
	Then, for each $i\in\mathbb S$, there exists $C_i>0$ such that
	\begin{equation*}
		U_i(x,t)\geq C_i(l_2-x),
		\quad (x,t)\in\overline\Omega\times\mathbb R.
	\end{equation*}
\end{lemma}

\begin{proof}
Applying Lemma 2.8 to $U$ with the coefficient matrix $A-\eta I$ gives $U_i(x,t)>0$ for all $(x,t)\in[l_1,l_2)\times\mathbb R$ and $i\in\mathbb S$. The claimed estimate then follows from Lemma 2.10 and Remark 2.11.
\end{proof}

\begin{proof}[Proof of Theorem~\ref{EX-PEV}(ii)]

Let $\varphi$ be the principal eigenfunction
from Theorem~\ref{EX-PEV}(i).
Since $\varphi\in\mathcal D_0$ and
$\mathcal L\varphi=\lambda_1\varphi$,
the definitions immediately give
	\begin{equation}\label{THM1.3-1}
		\lambda_p\le\mu_p\le\lambda_1\le\mu_p'.
	\end{equation}
	
	\textbf{Step 1. Proof that $\mu_p'\le\lambda_1$.}
	Suppose that $\mu_p'>\lambda_1$, and choose
	$\lambda_1<\mu<\mu_p'$.
	There exists $U\in\mathcal D_0$ such that
	\begin{equation*}
		\mathcal LU\geq\mu U\quad\text{in }Q.
	\end{equation*}
	By Lemmas~\ref{upper-estimate} and~\ref{lower-estimate},
	there exist $b_U,a_\varphi>0$ such that
	\begin{equation*}
		U_i(x,t)\leq b_U(l_2-x),
		\enspace
		\varphi_i(x,t)\geq a_\varphi(l_2-x)
		\enspace \text{in}\enspace \overline Q_T,\, i\in\mathbb S.
	\end{equation*}
	Hence $\varphi\geq(a_\varphi/b_U)U$ in $\overline{Q}_{T}$, and the set
	\begin{equation*}
		S_1:=\{\delta>0:\varphi\geq\delta U
		\text{ in }\overline Q_T\}
	\end{equation*}
	is nonempty.
	At any fixed $(x_0,t_0)\in Q_T$ and $i\in\mathbb S$, every
	$\delta\in S_1$ satisfies
	$\delta\leq\varphi_i(x_0,t_0)/U_i(x_0,t_0)$.
	Thus $0<\delta_*:=\sup S_1<\infty$.
	Passing to the limit along a sequence in $S_1$ gives
	$\varphi\geq\delta_*U$ on $\overline Q_T$.
	
	Set $W=\varphi-\delta_*U$.
	Then $W\geq\mathbf0$ on $\overline\Omega\times\mathbb R$,
	$W(l_2,t)=\mathbf0$, and
	\begin{equation*}
		\begin{aligned}
			(\mathcal L-\lambda_1)W
			&=(\mathcal L-\lambda_1)\varphi
			-\delta_*(\mathcal L-\lambda_1)U\\
			&\leq-\delta_*(\mu-\lambda_1)U
			\ll\mathbf0
			\quad\text{in }Q.
		\end{aligned}
	\end{equation*}
	If $W_i(x_0,t_0)=0$ at an interior point, then
	$\mathcal D_cW_i(x_0,t_0)=0$, since $W_i\geq0$.
	Consequently,
	\begin{equation*}
		\begin{aligned}
			\left ((\mathcal L-\lambda_1)W\right)_i(x_0,t_0)
			=d_i(x_0,t_0)\int_{l_1}^{l_2}
			k_i(x_0-y,t_0)W_i(y,t_0)\,dy+\sum_{j\neq i}a_{ij}(x_0,t_0)W_j(x_0,t_0)
			\geq0,
		\end{aligned}
	\end{equation*}
	a contradiction.
	Thus $W_i>0$ in $Q$ for every $i$, and $W\in\mathcal D$.
	Since $\mathcal LW\leq\lambda_1W$, Lemma~\ref{lower-estimate}
	gives $a_W>0$ such that
	\begin{equation*}
		W_i(x,t)\geq a_W(l_2-x)
		\geq\frac{a_W}{b_U}U_i(x,t)
		\enspace\text{in} \enspace \overline Q_T
		\enspace\text{for} \enspace
		i\in\mathbb S.
	\end{equation*}
	It follows that
	\begin{equation*}
		\varphi=\delta_*U+W
		\geq\left(\delta_*+\frac{a_W}{b_U}\right)U
		\enspace\text{in} \enspace \overline Q_T,
	\end{equation*}
	contrary to the definition of $\delta_*$.
	Hence $\mu_p'\leq\lambda_1$.

	\textbf{Step 2. Proof that $\lambda_p\geq\lambda_1$.}
	Suppose that $\lambda_p<\lambda_1$, and choose
	$\lambda_p<\mu<\lambda_1$.
	There exists $V\in\mathcal D$ such that
	\begin{equation*}
		\mathcal LV\leq\mu V  \enspace\text{in}\enspace
		Q.
	\end{equation*}
	Applying Lemma~\ref{lower-estimate} to $V$
and Lemma~\ref{upper-estimate} to $\varphi$ gives
	$a_V,b_\varphi>0$ such that
	\begin{equation*}
		V_i(x,t)\geq a_V(l_2-x),
		\enspace
		\varphi_i(x,t)\leq b_\varphi(l_2-x)
		\enspace\text{in}\enspace \overline Q_T,\,i\in\mathbb S.
	\end{equation*}
	Thus $V\geq(a_V/b_\varphi)\varphi$ and
	\begin{equation*}
		S_2:=\{\delta>0:V\geq\delta\varphi
		\text{ on }\overline Q_T\}\neq\varnothing.
	\end{equation*}
	As in Step 1,
	\begin{equation*}
		0<\delta^*:=\sup S_2<\infty,
		\quad V\geq\delta^*\varphi \enspace\text{in}\enspace \overline Q_T.
	\end{equation*}
	Set $W=V-\delta^*\varphi$.
	Then $W\geq\mathbf0$ on $\overline\Omega\times\mathbb R$ and $W(l_2,t)=V(l_2,t)\geq\mathbf0$.
	Moreover,
	\begin{equation*}
		\begin{aligned}
			(\mathcal L-\mu)W
			&=(\mathcal L-\mu)V
			-\delta^*(\mathcal L-\mu)\varphi\\
			&\leq-\delta^*(\lambda_1-\mu)\varphi
			\ll\mathbf0
			\quad\text{in }Q.
		\end{aligned}
	\end{equation*}
	The interior minimum argument in Step 1 gives $W_i>0$ in $Q$
	for every $i$, so $W\in\mathcal D$.
	Since $\mathcal LW\leq\mu W$, Lemma~\ref{lower-estimate}
	yields $a_W>0$ such that
	\begin{equation*}
		W_i(x,t)\geq a_W(l_2-x)
		\geq\frac{a_W}{b_\varphi}\varphi_i(x,t)
		\enspace\text{in}\enspace \overline Q_T,\,i\in\mathbb S.
	\end{equation*}
	Therefore,
	\begin{equation*}
		V=\delta^*\varphi+W
		\geq\left(\delta^*+\frac{a_W}{b_\varphi}\right)\varphi
		\enspace\text{in}\enspace \overline Q_T,
	\end{equation*}
	contradicting the definition of $\delta^*$.
	Hence $\lambda_p\geq\lambda_1$.
	
	Combining these inequalities with \eqref{THM1.3-1}
	yields $\mu_p=\mu_p'=\lambda_p=\lambda_1$.
\end{proof}

\subsection{Large-domain limit of the principal eigenvalue}

This subsection proves Theorem~\ref{cor:general-kernel-large-domain-limit}.
Throughout the proof, the spatial interval is $(-l,l)$, and
$D(x,t)=D(t)$ and $A(x,t)=A(t)$ are independent of $x$.
The principal eigenvalue on $(-l,l)$ is denoted by $\lambda_l$.
For clarity, $\mathcal L_l$, $\mathcal D_l$ and $\mathcal D_{0,l}$
denote the operator and the two classes of test functions
introduced in Section~1, with $\Omega=(-l,l)$.

The proof has two parts. First, the desired lower bound is
constructed for compactly supported dispersal kernels.
General kernels are then treated by truncation.

Recall the quantities introduced before
Theorem~\ref{cor:general-kernel-large-domain-limit}:
\begin{equation}\label{eq:LD-Floquet-matrix}
\begin{aligned}
&m_i(\nu,t)=\int_{\mathbb R}k_i(z,t)e^{-\nu z}\,dz,\\
&B_\nu(t)=A(t)+\nu c(t)I_m
+D(t)\operatorname{diag}\left(m_i(\nu,t)-1\right)_{i=1}^m.
\end{aligned}
\end{equation}
The admissible exponential weights form the set
$\mathcal E=\{\nu\in\mathbb R:
\int_0^T m_i(\nu,t)\,dt<\infty
\text{ for every }i\in\mathbb S\}$. For each \(\nu\in\mathcal E\), let \(\Phi_\nu(t,s)\) be the evolution matrix of \(z'=B_\nu(t)z\). On every bounded time interval, \(\Phi_\nu(\cdot,s)\) belongs to \(W^{1,1}\) and satisfies
\[
\partial_t\Phi_\nu(t,s)=B_\nu(t)\Phi_\nu(t,s)
\text{ for a.e. }t,\qquad
\Phi_\nu(s,s)=I_m,
\]
where $I_m$ is the $m\times m$ identity matrix.
The corresponding principal Floquet exponent is
\begin{align*}
\Lambda(\nu)=T^{-1}\log r\bigl(\Phi_\nu(T,0)\bigr),
\end{align*}
where $r$ denotes the spectral radius.
For $\nu\notin\mathcal E$, set $\Lambda(\nu)=+\infty$.

In the next two lemmas, a compactly supported kernel need not
have total mass one.

\begin{lemma}\label{Floquet-properties}
	Assume (H1), (H2) and (H4), with
	$D(x,t)=D(t)$ and $A(x,t)=A(t)$.
	Suppose that the nonnegative, $T$-periodic kernels
	$k_i\in C(\mathbb R^2)$ have uniformly compact support and are
	uniformly positive near the origin: for some
	$R,\delta,\kappa>0$,
	\begin{equation}\label{eq:LD-compact-kernel}
		k_i(z,t)=0 \text{ for } |z|>R,\quad
		k_i(z,t)\ge\kappa \text{ for } |z|\le\delta,\,
		i\in\mathbb S,\, t\in\mathbb R.
	\end{equation}
	Then $\mathcal E=\mathbb R$, and the following assertions hold.
	\begin{itemize}
		\item[(i)]
		The function $\Lambda:\mathbb R\to\mathbb R$ is
		continuous, strictly convex, and
		$\Lambda(\nu)\to+\infty$ as $|\nu|\to\infty$.
		Consequently, it has a unique minimizer $\nu_*$, and
		\begin{equation*}
			\Lambda_*:=\min_{\nu\in\mathbb R}\Lambda(\nu)
			=\Lambda(\nu_*)\in\mathbb R.
		\end{equation*}
		
		\item[(ii)]
		For each $\nu\in\mathbb R$, there is a strictly
		positive $T$-periodic vector $p^\nu\in C^1(\mathbb R;\mathbb R^m)$
		such that
		\begin{equation}\label{eq:LD-Floquet-vector}
			-(p^\nu)'(t)+B_\nu(t)p^\nu(t)
			=\Lambda(\nu)p^\nu(t).
		\end{equation}
		The vectors may be chosen so that
		$\sum_{i=1}^m p_i^\nu(0)=1$ and
		$(\nu,t)\mapsto p^\nu(t)$ is continuously
		differentiable in $\nu$ and $t$.
		For every $S>0$, there are constants $a_S,b_S>0$
		such that
\begin{equation}\label{eq:LD-periodic-vector-bounds}
a_S\le p_i^\nu(t)\le b_S,
\qquad |\partial_\nu p_i^\nu(t)|\le b_S,
\quad |\nu|\le S,\ t\in\mathbb R,\ i\in\mathbb S.
\end{equation}
\end{itemize}
\end{lemma}

\begin{proof}
The proof first constructs the positive periodic vector
in (ii). It then establishes the convexity and growth
of $\Lambda$ required in (i).
		
\medskip
\noindent\textbf{Step 1. Positive periodic vectors and their dependence on $\nu$.}

Since $B_\nu(t)$ is cooperative and irreducible
for every $t$, the standard positivity result
for cooperative linear systems
(see \cite[Section~3.1.2]{Zhao})
implies that $\Phi_\nu(t,s)$ is strictly positive
for all $t>s$.
In particular, $\Phi_\nu(T,0)$ is a strictly
positive matrix.
By the Perron--Frobenius theorem,
$r(\Phi_\nu(T,0))>0$ is an algebraically simple
eigenvalue with an eigenvector
$q^\nu\gg\mathbf0$.
We normalize $q^\nu$ by
$\sum_{i=1}^m q_i^\nu=1$.
Then
\begin{equation}\label{eq:LD-PF-eigenvector}
		\Phi_\nu(T,0)q^\nu
		=e^{T\Lambda(\nu)}q^\nu,
		\qquad q^\nu\gg\mathbf0,
		\qquad\sum_{i=1}^m q_i^\nu=1.
\end{equation}
The normalization makes $q^\nu$ unique.
		
Next, we justify the differentiability with respect
to the parameter $\nu$.
Compact support allows differentiation under the
integral sign:
\begin{equation}\label{eq:LD-moment-derivatives}
		\partial_\nu m_i(\nu,t)
		=-\int_{-R}^{R}z k_i(z,t)e^{-\nu z}\,dz,
		\qquad
		\partial_{\nu\nu}m_i(\nu,t)
		=\int_{-R}^{R}z^2k_i(z,t)e^{-\nu z}\,dz.
\end{equation}
		These derivatives are continuous in $(\nu,t)$,
		locally uniformly for $\nu$ in bounded intervals.
		Consequently, $B_\nu(t)$ and $\partial_\nu B_\nu(t)$
		are continuous in $(\nu,t)$.
		Differentiating the integral equation for
		$\Phi_\nu(t,0)$ with respect to $\nu$ gives
	\begin{equation}\label{eq:LD-evolution-parameter-derivative}
		\partial_\nu\Phi_\nu(t,0)
		=\int_0^t\Phi_\nu(t,s)
		\bigl(\partial_\nu B_\nu(s)\bigr)
		\Phi_\nu(s,0)\,ds,
		\qquad 0\le t\le T.
	\end{equation}
		In particular, the entries of $\Phi_\nu(T,0)$
		are continuously differentiable in $\nu$.
		Since its Perron eigenvalue is algebraically simple, the implicit function theorem applied to the eigenvalue equation and the normalization $\sum_{i=1}^m q_i^\nu=1$ shows that both $q^\nu$ and $r(\Phi_\nu(T,0))$ are continuously differentiable in $\nu$.
        Consequently, $\Lambda(\nu)=T^{-1}\log r(\Phi_\nu(T,0))$ is also continuously differentiable.
		
		Define
	\begin{equation*}
		p^\nu(t)=e^{-\Lambda(\nu)t}\Phi_\nu(t,0)q^\nu.
	\end{equation*}
  Since $q^\nu\gg\mathbf0$ and $\Phi_\nu(t,0)$
  preserves positivity, $p^\nu(t)\gg\mathbf0$
  for $t\geq0$.
  Moreover,
	\begin{equation*}
		(p^\nu)'(t)
		=(B_\nu(t)-\Lambda(\nu)I_m)p^\nu(t),
	\end{equation*}
		which is exactly \eqref{eq:LD-Floquet-vector}.
		Moreover, periodicity of $B_\nu$ implies
		$\Phi_\nu(t+T,0)=\Phi_\nu(t,0)\Phi_\nu(T,0)$.
		Using \eqref{eq:LD-PF-eigenvector}, we obtain
	\begin{equation*}
		\begin{aligned}
			p^\nu(t+T)
			&=e^{-\Lambda(\nu)(t+T)}
			\Phi_\nu(t,0)\Phi_\nu(T,0)q^\nu\\
			&=e^{-\Lambda(\nu)t}\Phi_\nu(t,0)q^\nu
			=p^\nu(t).
		\end{aligned}
	\end{equation*}
Thus $p^\nu$ extends to a positive $T$-periodic
$C^1$ function on $\mathbb R$.
The differentiability results established above
imply that $p^\nu(t)$ and
$\partial_\nu p^\nu(t)$ are jointly continuous
in $(\nu,t)$.
The same holds for $(p^\nu)'(t)$ by
\eqref{eq:LD-Floquet-vector}.
		
		Finally, fix $S>0$.
		Each component $p_i^\nu(t)$ is strictly positive
		on the compact set $[-S,S]\times[0,T]$.
		Thus their common minimum there is positive,
		while the components and their parameter derivatives
		have finite common upper bounds.
		By $T$-periodicity, the same bounds hold for every
		$t\in\mathbb R$, proving \eqref{eq:LD-periodic-vector-bounds}
		and completing the proof of (ii).
		
		\medskip
		\noindent\textbf{Step 2. Strict convexity of $\Lambda$.}
		The key point is that the exponential moments are
		uniformly strictly convex in $\nu$.
		Let
		$d_{\min}=\min_{i\in\mathbb S,\,t\in[0,T]}d_i(t)>0$.
		From \eqref{eq:LD-moment-derivatives} and the lower
		bound for $k_i$ in \eqref{eq:LD-compact-kernel},
	\begin{equation}\label{eq:LD-uniform-convexity}
		\begin{aligned}
			d_i(t)\partial_{\nu\nu}m_i(\nu,t)
			&\ge d_{\min}\kappa\int_{-\delta}^{\delta}
			z^2e^{-\nu z}\,dz\\
			&=d_{\min}\kappa\int_0^\delta
			z^2(e^{-\nu z}+e^{\nu z})\,dz\\
			&\ge 2d_{\min}\kappa\int_0^\delta z^2\,dz
			=\frac23d_{\min}\kappa\delta^3
			=:\gamma>0.
		\end{aligned}
	\end{equation}
		Here, the second inequality follows from
		$e^a+e^{-a}\ge2$.
		Since the second derivative of
		$\nu c(t)+d_i(t)(m_i(\nu,t)-1)$ is at least $\gamma$,
		subtracting $\gamma\nu^2/2$ from this function
		leaves a convex function.
		Therefore, for $\nu_1\ne\nu_2$,
		$0<\theta<1$, and
		$\nu_\theta=\theta\nu_1+(1-\theta)\nu_2$,
		we have
	\begin{equation}\label{eq:LD-weighted-moment-convexity}
		\begin{aligned}
			&\nu_\theta c(t)+d_i(t)\bigl(m_i(\nu_\theta,t)-1\bigr)\\
		\le&
			\theta\bigl[\nu_1c(t)+d_i(t)(m_i(\nu_1,t)-1)\bigr]
			+(1-\theta)\bigl[\nu_2c(t)+d_i(t)(m_i(\nu_2,t)-1)\bigr]\\
			&-\frac\gamma2\theta(1-\theta)(\nu_1-\nu_2)^2.
		\end{aligned}
	\end{equation}
		This is an inequality for the diagonal entries of
		$B_\nu(t)$. To obtain the corresponding inequality
		for the Floquet exponent, we use the positive
		periodic vectors constructed in Step 1.
		Set
	\begin{equation*}
		w_i(t)=\bigl(p_i^{\nu_1}(t)\bigr)^\theta
		\bigl(p_i^{\nu_2}(t)\bigr)^{1-\theta},\quad i\in\mathbb S.
	\end{equation*}
		Because both $p^{\nu_1}$ and $p^{\nu_2}$ are
		strictly positive and $T$-periodic, so is $w$.
		Differentiating the logarithm of $w_i$ gives
	\begin{equation*}
		\frac{w_i'(t)}{w_i(t)}
		=\theta\frac{(p_i^{\nu_1})'(t)}{p_i^{\nu_1}(t)}
		+(1-\theta)\frac{(p_i^{\nu_2})'(t)}{p_i^{\nu_2}(t)}.
	\end{equation*}
		We next compare the coupling terms.
		For $j\ne i$, the weighted arithmetic--geometric
		mean inequality gives
	\begin{equation*}
		\begin{aligned}
			\frac{w_j(t)}{w_i(t)}
			=\left(\frac{p_j^{\nu_1}(t)}{p_i^{\nu_1}(t)}\right)^\theta
			\left(\frac{p_j^{\nu_2}(t)}{p_i^{\nu_2}(t)}\right)^{1-\theta}\le\theta\frac{p_j^{\nu_1}(t)}{p_i^{\nu_1}(t)}
			+(1-\theta)\frac{p_j^{\nu_2}(t)}{p_i^{\nu_2}(t)}.
		\end{aligned}
	\end{equation*}
	 Since $A(t)$ is cooperative, summing over $j$ gives
	\begin{equation}\label{eq:LD-coupling-convexity}
		\frac{(A(t)w(t))_i}{w_i(t)}
		\le\theta\frac{(A(t)p^{\nu_1}(t))_i}{p_i^{\nu_1}(t)}
		+(1-\theta)\frac{(A(t)p^{\nu_2}(t))_i}{p_i^{\nu_2}(t)}.
	\end{equation}
		Now expand the $i$th component of
		$-w'(t)+B_{\nu_\theta}(t)w(t)$ and divide by $w_i(t)$:
	\begin{equation*}
		\begin{aligned}
			\frac{-w_i'(t)+(B_{\nu_\theta}(t)w(t))_i}{w_i(t)}=&-\frac{w_i'(t)}{w_i(t)}
			+\frac{(A(t)w(t))_i}{w_i(t)}
			+\nu_\theta c(t)+d_i(t)(m_i(\nu_\theta,t)-1)\\
			\le&
			\theta\frac{-(p_i^{\nu_1})'(t)
				+(B_{\nu_1}(t)p^{\nu_1}(t))_i}{p_i^{\nu_1}(t)}
			+(1-\theta)\frac{-(p_i^{\nu_2})'(t)
				+(B_{\nu_2}(t)p^{\nu_2}(t))_i}{p_i^{\nu_2}(t)}\\
			&-\frac\gamma2\theta(1-\theta)(\nu_1-\nu_2)^2\\
			=&\theta\Lambda(\nu_1)+(1-\theta)\Lambda(\nu_2)
			-\frac\gamma2\theta(1-\theta)(\nu_1-\nu_2)^2.
		\end{aligned}
	\end{equation*}
	By \eqref{eq:LD-Floquet-vector}, the preceding
inequality gives
$w'(t)+bw(t)\ge B_{\nu_\theta}(t)w(t)$,
where
\begin{equation*}
    b=\theta\Lambda(\nu_1)+(1-\theta)\Lambda(\nu_2)
    -\frac{\gamma}{2}\theta(1-\theta)(\nu_1-\nu_2)^2.
\end{equation*}
Since $w$ is $T$-periodic, the comparison principle
for cooperative systems yields
\begin{equation*}
    \Phi_{\nu_\theta}(T,0)w(0)
    \le e^{bT}w(T)=e^{bT}w(0).
\end{equation*}
Multiplying by a positive left Perron eigenvector
of $\Phi_{\nu_\theta}(T,0)$ and using $w(0)\gg\mathbf0$,
we obtain
$r(\Phi_{\nu_\theta}(T,0))\le e^{bT}$.
Consequently,
\begin{equation}\label{eq:LD-Lambda-strict-convexity}
    \Lambda(\nu_\theta)
    \le\theta\Lambda(\nu_1)+(1-\theta)\Lambda(\nu_2)
    -\frac{\gamma}{2}\theta(1-\theta)(\nu_1-\nu_2)^2.
\end{equation}
Since $\gamma>0$, this proves the strict
convexity of $\Lambda$.

\textbf{Step 3. Growth of $\Lambda(\nu)$ as $|\nu|\to\infty$.}
The lower bound for $k_i$ in
\eqref{eq:LD-compact-kernel} implies that
		$m_i(\nu,t)$ grows exponentially as $|\nu|\to\infty$.
		For $\nu\ge0$, integrate over
		$[-\delta,-\delta/2]$; on this interval,
		$e^{-\nu z}\ge e^{\nu\delta/2}$.
		For $\nu<0$, integrate over
		$[\delta/2,\delta]$; there,
		$e^{-\nu z}\ge e^{|\nu|\delta/2}$.
		Both cases give
	\begin{equation}\label{eq:LD-moment-growth}
		m_i(\nu,t)\ge
		\frac{\kappa\delta}{2}e^{|\nu|\delta/2},
		\qquad i\in\mathbb S,\ t\in\mathbb R.
	\end{equation}
		We now turn this scalar lower bound into a lower
		bound for $\Lambda(\nu)$.
		Let
	\begin{equation*}
		a_{\min}=\min_{i\in\mathbb S,\,t\in[0,T]}a_{ii}(t),
		\qquad
		d_{\max}=\max_{i\in\mathbb S,\,t\in[0,T]}d_i(t),
		\qquad
		\|c\|_\infty=\max_{t\in[0,T]}|c(t)|.
	\end{equation*}
		Since $A(t)$ is cooperative,
		$(A(t)\mathbf1)_i\ge a_{ii}(t)\ge a_{\min}$.
		Also, $\nu c(t)\ge-|\nu|\|c\|_\infty$.
		Using \eqref{eq:LD-moment-growth} and
		$d_{\min}\le d_i(t)\le d_{\max}$, we obtain
	\begin{equation}\label{eq:LD-B-growth}
		B_\nu(t)\mathbf1\ge b(\nu)\mathbf1,
		\qquad
		b(\nu):=a_{\min}-\|c\|_\infty|\nu|-d_{\max}
		+\frac{d_{\min}\kappa\delta}{2}e^{|\nu|\delta/2}.
	\end{equation}
	Thus $e^{b(\nu)t}\mathbf1$ is a subsolution
of $z'=B_\nu(t)z$.
By the comparison principle,
$\Phi_\nu(T,0)\mathbf1 \ge e^{b(\nu)T}\mathbf1$.
Multiplying by a positive left Perron eigenvector
of $\Phi_\nu(T,0)$ gives
\begin{equation*}
    \Lambda(\nu)\ge b(\nu)
    \longrightarrow+\infty
    \qquad\text{as }|\nu|\to\infty.
\end{equation*}
Since $\Lambda$ is continuous and strictly convex,
it attains its minimum at a unique point,
completing the proof of (i).
\end{proof}

\begin{lemma}\label{thm:large-domain-limit}
	Assume the hypotheses of Lemma~\ref{Floquet-properties}.
	For every $\varepsilon>0$, there exists $l_\varepsilon>0$
	such that for every $l\ge l_\varepsilon$ one can find
	$U_l\in\mathcal D_{0,l}$ satisfying
	\begin{equation}\label{eq:LD-lower-test-function}
		\mathcal L_lU_l\ge(\Lambda_*-\varepsilon)U_l
		\qquad\text{in }(-l,l)\times\mathbb R.
	\end{equation}
\end{lemma}

\begin{proof}
		The main idea is to construct a positive function whose
		spatial dependence is approximately exponential. Near
		$x=-l$ it will have the form $e^{Sx}p^S(t)$, and near
		$x=l$ it will have the form $e^{-Sx}p^{-S}(t)$, where
		$S>0$ is large. These two expressions are connected
		by allowing the exponent to vary slowly with $x$.
		A further factor will make the function vanish at $x=l$.
		The construction is arranged so that the errors arising
		from these two modifications, and from restricting the
		nonlocal integral to $(-l,l)$, can all be controlled.

	\textbf{Step 1. Choice of the exponent and the boundary cutoff.}
		Let $R,\delta,\kappa>0$ be as in
		\eqref{eq:LD-compact-kernel}, and fix $0<\rho<\delta$.
		Choose a nondecreasing function $b\in C^1([0,\infty))$
		such that
	\begin{equation*}
		\begin{gathered}
			b(0)=0,\qquad 0<b(r)\le1\quad(r>0),\\
			b(r)=r/\rho\quad(0\le r\le\rho/2),
			\qquad b(r)=1\quad(r\ge\rho).
		\end{gathered}
	\end{equation*}
		In particular, since $b$ is nondecreasing,
		$b(r)\ge1/4$ when $r\ge\rho/4$.
		The factor $b(l-x)$ will be equal to $1$ away from
		the right endpoint and will vanish at $x=l$.
		
		For later use, set
	\begin{equation*}
		K_+=\max_{i\in\mathbb S,\,t\in[0,T]}\int_0^R k_i(z,t)\,dz,
		\qquad
		K_-=\max_{i\in\mathbb S,\,t\in[0,T]}\int_{-R}^0 k_i(z,t)\,dz,
\qquad
        K_{0}=K_++K_-.
	\end{equation*}
		The constants $K_+$ and $K_-$ bound the parts of the
		nonlocal integral that may be lost at the left and
		right endpoints, respectively.
		Recall also that
		$d_{\min}=\min_{i\in\mathbb S,\,t\in[0,T]} d_i(t)>0$,
		$d_{\max}=\max_{i\in\mathbb S,\,t\in[0,T]}d_i(t)$, and
		$\|c\|_\infty=\max_{t\in[0,T]}|c(t)|$.
		By Lemma~\ref{Floquet-properties}(i),
		$\Lambda(S)$ and $\Lambda(-S)$ both tend to $+\infty$
		as $S\to\infty$. Moreover,
		$\int_0^{\rho/4}ze^{Sz}\,dz\to\infty$ as
		$S\to\infty$. Hence a single $S>0$ can be chosen
		large enough that
	\begin{equation}\label{eq:LD-endpoint-choice}
		\begin{aligned}
			&\Lambda(S)-d_{\max}K_+\ge\Lambda_*-\varepsilon,\\
			&\Lambda(-S)-d_{\max}K_-
			-4\|c\|_\infty\|b'\|_\infty\ge\Lambda_*-\varepsilon,\\
			&d_{\min}\kappa\int_0^{\rho/4}ze^{Sz}\,dz
			\ge\|c\|_\infty.
		\end{aligned}
	\end{equation}
		The first two conditions will control the endpoint
		errors when $b(l-x)$ stays bounded away from zero.
		The third will be used where $b(l-x)$ is small and
		its logarithmic derivative becomes large.

     \textbf{Step 2. Construction of the slowly varying exponential profile.}
		Let $L>0$ be a parameter to be chosen sufficiently
		large. Take $q\in C^1(\mathbb R)$ such that
	\begin{equation*}
		\begin{gathered}
			q(x)=S\quad(x\le-L),\qquad
			q(x)=-S\quad(x\ge L),\\
			|q(x)|\le S,\qquad
			\|q'\|_\infty\le C_S/L.
		\end{gathered}
	\end{equation*}
		For instance, let $\chi\in C^1(\mathbb R)$ satisfy
		$0\le\chi\le1$, $\chi(s)=0$ for $s\le-1$, and
		$\chi(s)=1$ for $s\ge1$, and set
		$q(x)=S[1-2\chi(x/L)]$.
		Then these bounds hold with
		$C_S=2S\|\chi'\|_\infty$.
		Notice that $S$ has already been fixed; therefore
		$C_S$ is independent of $L$ and $l$.
		
		Let $p^\nu(t)$ be the positive periodic Floquet vector
		from Lemma~\ref{Floquet-properties}(ii). Define
	\begin{equation}\label{eq:LD-slow-profile}
		W(x,t)=\exp\!\left(\int_0^xq(r)\,dr\right)p^{q(x)}(t).
	\end{equation}
		The function $W$ is positive, $C^1$ in $(x,t)$,
		and $T$-periodic in $t$. Since $q(x)=S$ for $x\le-L$,
		the integral in \eqref{eq:LD-slow-profile} equals
		$Sx+C_-$ in that region, for a constant $C_-$.
		Similarly, it equals $-Sx+C_+$ for $x\ge L$.
		Consequently,
	\begin{equation*}
		W(x,t)=
		\begin{cases}
			e^{C_-}e^{Sx}p^S(t),&x\le-L,\\
			e^{C_+}e^{-Sx}p^{-S}(t),&x\ge L.
		\end{cases}
	\end{equation*}
		We next estimate the error made by replacing the
		constant exponent $\nu$ with $q(x)$.
		Write $F(x)=\int_0^xq(r)\,dr$.
		For $|z|\le R$, the mean value theorem gives
	\begin{equation*}
		\begin{aligned}
			|F(x-z)-F(x)+zq(x)|
			&=\left|\int_{x-z}^x[q(x)-q(r)]\,dr\right|\\
			&\le\|q'\|_\infty\int_0^{|z|}s\,ds
			=\frac12\|q'\|_\infty z^2.
		\end{aligned}
	\end{equation*}
		Likewise, since $|q(x-z)-q(x)|\le\|q'\|_\infty|z|$,
		the $\nu$-derivative bound in
		\eqref{eq:LD-periodic-vector-bounds} gives
	\begin{equation*}
		|p_i^{q(x-z)}(t)-p_i^{q(x)}(t)|
		\le b_S\|q'\|_\infty|z|.
	\end{equation*}
		The ratio appearing in the nonlocal integral is
	\begin{equation*}
		\frac{W_i(x-z,t)}{W_i(x,t)}
		=e^{F(x-z)-F(x)}
		\frac{p_i^{q(x-z)}(t)}{p_i^{q(x)}(t)}.
	\end{equation*}
		The preceding two estimates, together with the bounds
		$a_S\le p_i^\nu(t)\le b_S$ for $|\nu|\le S$,
		give an estimate of this ratio.
		Indeed, write
		$F(x-z)-F(x)=-zq(x)+E$, where
		$|E|\le\|q'\|_\infty R^2/2$.
		If $L$ is large enough that $\|q'\|_\infty\le1$,
		then $|e^E-1|\le e^{R^2/2}|E|$ and
		$|e^{-zq(x)}|\le e^{SR}$.
		It follows that
	\begin{equation}\label{eq:LD-kernel-error}
		\begin{aligned}
		\left|\frac{W_i(x-z,t)}{W_i(x,t)}-e^{-q(x)z}\right|\le e^{SR}
			\left[e^{R^2/2}\frac{b_S}{a_S}\frac{R^2}{2}
			+\frac{b_S}{a_S}R\right]\|q'\|_\infty
			\le C\|q'\|_\infty,
			\quad |z|\le R.
		\end{aligned}
	\end{equation}
		Here and below, $C$ denotes a constant depending
		on $S,R,a_S,b_S$ and the bounds on the coefficients,
		but not on $L$, $l$, $x$ or $t$.
		
		The spatial derivative of $W$ is
	\begin{equation*}
		\partial_xW_i(x,t)
		=q(x)W_i(x,t)
		+e^{F(x)}q'(x)\partial_\nu p_i^{q(x)}(t),
	\end{equation*}
		whereas $\partial_tW_i(x,t)=e^{F(x)}(p_i^{q(x)})'(t)$.
		On the whole real line, the nonlocal term is
	\begin{equation*}
		\int_{\mathbb R}k_i(x-y,t)W_i(y,t)\,dy
		=W_i(x,t)\int_{-R}^R k_i(z,t)
		\frac{W_i(x-z,t)}{W_i(x,t)}\,dz.
	\end{equation*}
        Let \(\mathcal L_\infty\) denote the whole-line operator \(\mathcal L_\infty w=-\partial_t w+c(t)\partial_x w+D(t)\mathcal K_{\mathbb R}[w]+A(t)w\), where \((\mathcal K_{\mathbb R}[w])_i=\int_{\mathbb R}k_i(x-y,t)w_i(y,t)\,dy-w_i(x,t)\). 
		Substituting these expressions into
		$\mathcal L_\infty W$
		and using the Floquet equation
     	$-(p^{q(x)})'(t)+B_{q(x)}(t)p^{q(x)}(t)
		=\Lambda(q(x))p^{q(x)}(t)$,
	we obtain
	\begin{equation*}
		\begin{aligned}
			&\frac{(\mathcal L_\infty W)_i(x,t)}{W_i(x,t)}
			-\Lambda(q(x))\\
			=&c(t)q'(x)
			\frac{\partial_\nu p_i^{q(x)}(t)}{p_i^{q(x)}(t)}
			+d_i(t)\int_{-R}^R k_i(z,t)
			\left[\frac{W_i(x-z,t)}{W_i(x,t)}-e^{-q(x)z}\right]dz.
		\end{aligned}
	\end{equation*}
	By \eqref{eq:LD-periodic-vector-bounds},
	$|\partial_\nu p_i^{q(x)}|/p_i^{q(x)}\le b_S/a_S$.
	The estimate
	\eqref{eq:LD-kernel-error} therefore yields
	\begin{equation*}
		\begin{aligned}
		\left|(\mathcal L_\infty W)_i
			-\Lambda(q(x))W_i\right|\le\left(\|c\|_\infty\frac{b_S}{a_S}
			+d_{\max}K_0 C\right)\|q'\|_\infty W_i
			\le\frac{C C_S}{L}W_i.
		\end{aligned}
	\end{equation*}
		Enlarge $L$, if necessary, so that
		$\|q'\|_\infty\le1$ and $C C_S/L\le\varepsilon$.
		As $\Lambda(q(x))\ge\Lambda_*$, it follows that
	\begin{equation}\label{eq:LD-interior-subsolution}
		\mathcal L_\infty W\ge(\Lambda_*-\varepsilon)W
		\text{ in }\mathbb R\times\mathbb R.
	\end{equation}
		Finally, choose $l$ so large that
	\begin{equation*}
		-l+2R\le-L,\qquad l-(2R+\rho)\ge L.
	\end{equation*}
		Define
	\begin{equation}\label{eq:LD-test-profile}
		U_l(x,t)=b(l-x)W(x,t),\qquad -l\le x\le l.
	\end{equation}
		This function is $C^1$, $T$-periodic in $t$,
		positive for $-l\le x<l$, and zero at $x=l$.
		Thus $U_l\in\mathcal D_{0,l}$. It remains to verify
		\eqref{eq:LD-lower-test-function}. The interval is
		split into three regions so that the changes in the
		nonlocal integral can be calculated explicitly.
	
     \textbf{Step 3. Verification on the finite interval.}
	
	\emph{Left endpoint: $-l<x<-l+R$.}
		By the choice of $l$, $x+R\le-L$.
		All points $y\in[-l,l]$ with $|x-y|\le R$ therefore
		lie in the region where $q=S$. Moreover,
		$l-y\ge\rho$ there, so $b(l-y)=1$.
		Thus, by the endpoint expressions for $W$ given above,
		$U_l(y,t)=e^{C_-}e^{Sy}p^S(t)$ at every point
		contributing to the nonlocal integral.
		With the change of variable $z=x-y$,
	\begin{align*}
		\begin{aligned}
			\int_{-l}^l k_i(x-y,t)U_{l,i}(y,t)\,dy
			&=U_{l,i}(x,t)
			\int_{x-l}^{x+l}k_i(z,t)e^{-Sz}\,dz\\
		&=U_{l,i}(x,t)
			\int_{-R}^{x+l}k_i(z,t)e^{-Sz}\,dz.
		\end{aligned}
	\end{align*}
Since $0<x+l<R$ and $k_i$ is supported in $[-R,R]$,
the integral over the whole line differs only by
the contribution over $(x+l,R]$.
Using
$\mathcal L_\infty(e^{Sx}p^S)
=\Lambda(S)e^{Sx}p^S$, we obtain
	\begin{equation*}
		\begin{aligned}
			\frac{(\mathcal L_lU_l)_i}{U_{l,i}}
			&=\Lambda(S)
			-d_i(t)\int_{x+l}^R k_i(z,t)e^{-Sz}\,dz
			\ge\Lambda(S)-d_{\max}\int_0^R k_i(z,t)\,dz\\
			&\ge\Lambda(S)-d_{\max}K_+
			\ge\Lambda_*-\varepsilon,
		\end{aligned}
	\end{equation*}
	where $e^{-Sz}\le1$ was used for $z\ge0$.
		
		\emph{Interior: $-l+R\le x\le l-(R+\rho)$.}
		Every point $y$ with $|x-y|\le R$ belongs to
		$[-l,l-\rho]$. Hence $b(l-y)=1$, and also
		$b(l-x)=1$. Therefore $U_l=W$ at $x$ and at all
		points contributing to the nonlocal integral.
		The kernels vanish outside $[-R,R]$, so
	\begin{equation*}
		\begin{aligned}
			\int_{-l}^l k_i(x-y,t)U_{l,i}(y,t)\,dy=\int_{x-R}^{x+R}k_i(x-y,t)W_i(y,t)\,dy
			=\int_{\mathbb R}k_i(x-y,t)W_i(y,t)\,dy.
		\end{aligned}
	\end{equation*}
		The time and spatial derivatives agree as well,
		so \eqref{eq:LD-interior-subsolution} gives
	\begin{equation*}
		\mathcal L_lU_l=\mathcal L_\infty W
		\ge(\Lambda_*-\varepsilon)W
		=(\Lambda_*-\varepsilon)U_l.
	\end{equation*}

\emph{Right endpoint: $l-(R+\rho)<x<l$.}
Put $r=l-x$, so $0<r<R+\rho$.
By the choice of $l$, $x-R\ge L$.
Thus $q=-S$ at $x$ and at all the points contributing
to the nonlocal integral. In particular,
$W(y,t)=e^{C_+}e^{-Sy}p^{-S}(t)$ there.
Since $l-x=r$, we have
\begin{equation*}
		U_{l,i}(x,t)=e^{C_+}e^{-Sx}p_i^{-S}(t)b(r),
		\qquad
		\partial_xU_{l,i}
		=-S U_{l,i}-e^{C_+}e^{-Sx}p_i^{-S}(t)b'(r).
\end{equation*}
Under the change of variables $ z=x-y $, the condition $-l\le y\le l$  becomes $-r=x-l \le z \le x+l$. Since $x+l\ge R$
for the present choice of $l$, the kernel support
gives
	\begin{equation*}
		\frac{1}{U_{l,i}(x,t)}
		\int_{-l}^l k_i(x-y,t)U_{l,i}(y,t)\,dy
		=\int_{-\min\{R,r\}}^R
		k_i(z,t)e^{Sz}\frac{b(r+z)}{b(r)}\,dz.
	\end{equation*}
		For the whole-line mode $e^{-Sx}p^{-S}(t)$,
		the Floquet equation reads
	\begin{equation*}
		-(p^{-S})'(t)+B_{-S}(t)p^{-S}(t)
		=\Lambda(-S)p^{-S}(t).
	\end{equation*}
		Insert the formulas for $U_l$, its spatial derivative,
		and the preceding integral formula into $\mathcal L_lU_l$.
		Then subtract the Floquet equation, using
		$m_i(-S,t)=\int_{-R}^R k_i(z,t)e^{Sz}\,dz$.
		After separating the part of the integral omitted at
		$z<-\min\{R,r\}$, this gives the exact identity
	\begin{equation*}
		\begin{aligned}
			\frac{(\mathcal L_lU_l)_i}{U_{l,i}}
			=&\Lambda(-S)-c(t)\frac{b'(r)}{b(r)}+d_i(t)\int_{-\min\{R,r\}}^R
			k_i(z,t)e^{Sz}
			\left[\frac{b(r+z)}{b(r)}-1\right]dz\\
			&-d_i(t)\int_{-R}^{-\min\{R,r\}}
			k_i(z,t)e^{Sz}\,dz.
		\end{aligned}
	\end{equation*}
		Combining the negative-$z$ part of the first integral with the omitted-tail term, we obtain a lower bound of $-d_{\max} K_-$.
		Indeed, $b(r+z)\ge0$, and $e^{Sz}\le1$ for $z\le0$.
		The last term is nonpositive and,
		together with the negative-$z$ part just mentioned,
		is bounded below by $-d_{\max}K_-$.
		For $z\ge0$, monotonicity of $b$ gives
		$b(r+z)\ge b(r)$, so this part is nonnegative.
		Consequently,
	\begin{equation}\label{eq:LD-right-lower}
		\begin{aligned}
			\frac{(\mathcal L_lU_l)_i}{U_{l,i}}
			\ge\Lambda(-S)-d_{\max}K_-
			-c(t)\frac{b'(r)}{b(r)}+d_i(t)\int_0^R k_i(z,t)e^{Sz}
			\left[\frac{b(r+z)}{b(r)}-1\right]dz.
		\end{aligned}
	\end{equation}
		There are two cases, according to the distance $r$
		from the boundary.
		
		If $r\ge\rho/4$, then $b(r)\ge1/4$, and hence
		$0\le b'(r)/b(r)\le4\|b'\|_\infty$.
		Dropping the nonnegative integral in
		\eqref{eq:LD-right-lower} and applying the second
		inequality in \eqref{eq:LD-endpoint-choice}, we obtain
	\begin{equation*}
		\begin{aligned}
			\frac{(\mathcal L_lU_l)_i}{U_{l,i}}
			&\ge\Lambda(-S)-d_{\max}K_-
			-4\|c\|_\infty\|b'\|_\infty\\
			&\ge\Lambda_*-\varepsilon.
		\end{aligned}
	\end{equation*}
		If $0<r<\rho/4$, then $b(r)=r/\rho$ and
		$b'(r)/b(r)=1/r$. Also, for $0\le z\le\rho/4$,
		$r+z<\rho/2$, so
	\begin{equation*}
		\frac{b(r+z)}{b(r)}-1
		=\frac{(r+z)/\rho}{r/\rho}-1
		=\frac{z}{r}.
	\end{equation*}
		Since $k_i(z,t)\ge\kappa$ for
		$0\le z\le\rho/4<\delta$, the nonnegative integral
		in \eqref{eq:LD-right-lower} satisfies
	\begin{equation*}
		\begin{aligned}
			&d_i(t)\int_0^R k_i(z,t)e^{Sz}
			\left[\frac{b(r+z)}{b(r)}-1\right]dz\\
			\ge&\frac{d_{\min}\kappa}{r}
			\int_0^{\rho/4}ze^{Sz}\,dz
			\ge\frac{\|c\|_\infty}{r}
			\ge c(t)\frac{b'(r)}{b(r)},
		\end{aligned}
	\end{equation*}
		where the second inequality follows from the third
		condition in \eqref{eq:LD-endpoint-choice}.
		Thus \eqref{eq:LD-right-lower} yields
	\begin{equation*}
		\frac{(\mathcal L_lU_l)_i}{U_{l,i}}
		\ge\Lambda(-S)-d_{\max}K_-
		\ge\Lambda_*-\varepsilon.
	\end{equation*}
	Combining the estimates in the three regions, we obtain
$\mathcal L_lU_l\ge(\Lambda_*-\varepsilon)U_l$
in $(-l,l)\times\mathbb R$, as required.
\end{proof}

\begin{lemma}\label{lem:general-kernel-minimizer}
Assume (H1)--(H4), with $D(x,t)=D(t)$ and
$A(x,t)=A(t)$.
Then the following properties hold for the Floquet
exponent $\Lambda$ defined in \eqref{FLOQUENT}.
\begin{itemize}
\item[(i)]
There are compactly supported kernels $k_{i,n}$
satisfying the hypotheses of
Lemma~\ref{Floquet-properties} such that
\begin{equation*}
 0\le k_{i,n}\le k_{i,n+1}\le k_i,
 \qquad k_{i,n}\to k_i\quad\text{locally uniformly},
\end{equation*}
and $0<\int_{\mathbb R}k_{i,n}(z,t)\,dz\le1$.

\item[(ii)]
Let $\Lambda_n$ be the Floquet exponent obtained
by replacing $k_i$ with $k_{i,n}$ in
\eqref{eq:LD-Floquet-matrix}, and let $\nu_n$
be its unique minimizer. Then
$\Lambda_n(\nu)\uparrow\Lambda(\nu)$ for every
$\nu\in\mathbb R$ (allowing the value $+\infty$).
Moreover, $\Lambda$ is lower semicontinuous on
$\mathbb R$, strictly convex on $\mathcal E$, and
has a unique minimizer $\nu_*\in\mathcal E$.
In addition,
\begin{equation}\label{eq:LD-truncated-minimizers}
 \nu_n\to\nu_*,\qquad
 \min_{\nu\in\mathbb R}\Lambda_n(\nu)
 \to\Lambda(\nu_*)
 =\min_{\nu\in\mathcal E}\Lambda(\nu)\in\mathbb R.
\end{equation}
\end{itemize}
\end{lemma}

\begin{proof}
\textbf{Step 1. Truncate the kernels.}
By (H3), there exist constants $0<\delta\le1$ and $\kappa>0$
such that $k_i(z,t)\ge\kappa$ for
$|z|\le\delta$, all $t$ and $i$.
Choose an even function $\chi\in C_c^\infty(\mathbb R)$ that is
nonincreasing on $[0,\infty)$ and satisfies 
\begin{equation*}
    0\le\chi\le1, \enspace \chi=1 \enspace \text{for} \enspace x\in [-1,1], \enspace \chi=0 \enspace \text{for} \enspace |x|>2.
\end{equation*}
Set $k_{i,n}(z,t):=\chi(z/n)k_i(z,t)$ for $n\ge1$.
Then $k_{i,n}$ has the desired properties.

\textbf{Step 2. Pass to the limit in the Floquet exponents.}
For each $n$, let $B_{n,\nu}(t)$ and
$\Phi_{n,\nu}(t,s)$ denote the matrix in
\eqref{eq:LD-Floquet-matrix} and its evolution
matrix with $k_i$ replaced by $k_{i,n}$.
Since $k_{i,n}\le k_{i,n+1}$, comparison for
cooperative linear systems gives
$\Phi_{n,\nu}(T,0)\le\Phi_{n+1,\nu}(T,0)$,
and hence
$\Lambda_n(\nu)\le\Lambda_{n+1}(\nu)$.

First fix $\nu\in\mathcal E$.
The definition of $\mathcal E$ and the monotone
convergence theorem yield
\begin{equation*}
 \int_0^T|m_i(\nu,t)-m_{i,n}(\nu,t)|\,dt\to0,
 \qquad n\to\infty,
\end{equation*}
where $m_{i,n}$ is the exponential moment of
$k_{i,n}$. Hence $B_{n,\nu}\to B_\nu$ in
$L^1(0,T)$. The integral equations for the
evolution matrices and Gronwall's inequality imply
$\Phi_{n,\nu}(T,0)\to\Phi_\nu(T,0)$.
By continuity of the spectral radius,
$\Lambda_n(\nu)\to\Lambda(\nu)$.

Next, consider the case $\nu\notin\mathcal E$.
By the definition of $\mathcal E$, there exists
$i\in\mathbb S$ such that
$\int_0^T m_i(\nu,t)\,dt=+\infty$.
Since $k_{i,n}$ is obtained by truncating the
nonnegative kernel $k_i$, the monotone convergence
theorem gives
\begin{equation}\label{eq:LD-divergent-moment}
\int_0^T m_{i,n}(\nu,t)\,dt
\longrightarrow+\infty
\quad\text{as }n\to\infty.
	\end{equation}
	For each $n$, let $p^{n,\nu}$ be the strictly
		positive $T$-periodic Floquet vector associated
		with $B_{n,\nu}(t)$. It satisfies
		\[
		-(p^{n,\nu})'(t)+B_{n,\nu}(t)p^{n,\nu}(t)
		=\Lambda_n(\nu)p^{n,\nu}(t).
		\]

	Dividing the $i$th equation by $p_i^{n,\nu}(t)$
and integrating over $[0,T]$, the derivative term
vanishes by periodicity.
Since $A(t)$ is cooperative, the off-diagonal
terms are nonnegative. Hence
	\begin{equation}\label{eq:LD-divergent-Floquet-bound}
		\begin{aligned}
			\Lambda_n(\nu)
			&\ge\frac1T\int_0^T
			\bigl(a_{ii}(t)+\nu c(t)-d_i(t)\bigr)\,dt+\frac1T\int_0^T
			d_i(t)m_{i,n}(\nu,t)\,dt\\
			&\ge\frac1T\int_0^T
			\bigl(a_{ii}(t)+\nu c(t)-d_i(t)\bigr)\,dt
			+\frac{d_{\min}}T
			\int_0^T m_{i,n}(\nu,t)\,dt.
		\end{aligned}
	\end{equation}
	Since $d_{\min}>0$, the first term is independent
of $n$, while the second tends to $+\infty$
by \eqref{eq:LD-divergent-moment}.
Therefore,
$\Lambda_n(\nu)\to +\infty
=\Lambda(\nu)$.
Together with the result for $\nu\in\mathcal E$,
this proves pointwise convergence for all
$\nu\in\mathbb R$.
		
As already observed, $\Lambda_{n}$ is nondecreasing in $n$. Together with the pointwise convergence above, this gives
\begin{equation}\label{eq:LD-pointwise-Floquet-limit}
    \Lambda_n(\nu)\uparrow\Lambda(\nu),
    \qquad \nu\in\mathbb R.
\end{equation}
Since each $\Lambda_n$ is continuous by
Lemma~\ref{Floquet-properties} and
$\Lambda=\sup_{n\ge1}\Lambda_n$,
it follows that $\Lambda$ is lower semicontinuous
on $\mathbb R$.

\textbf{Step 3. Prove strict convexity and convergence
of minimizers.}
The truncated kernels share the positive lower
bound near zero from Step 1. Therefore the
strong-convexity estimate
\eqref{eq:LD-Lambda-strict-convexity} holds for
all $\Lambda_n$ with the same constant
$\gamma=2d_{\min}\kappa\delta^3/3>0$:
\begin{equation}\label{eq:Lambda-strong-convex}
 \Lambda_n\bigl(\theta\nu_1+(1-\theta)\nu_2\bigr)
 \le\theta\Lambda_n(\nu_1)+(1-\theta)\Lambda_n(\nu_2)
 -\frac\gamma2\theta(1-\theta)|\nu_1-\nu_2|^2.
\end{equation}
By H\"older's inequality,
$m_i(\theta\nu_1+(1-\theta)\nu_2,t)
\le m_i(\nu_1,t)^\theta
m_i(\nu_2,t)^{1-\theta}$.
Applying H\"older's inequality in $t$ shows
that $\mathcal E$ is convex. Passing to the limit in
\eqref{eq:Lambda-strong-convex} for
$\nu_1,\nu_2\in\mathcal E$ gives
\begin{equation*}
 \Lambda\left(\theta\nu_1+(1-\theta)\nu_2\right)
 \le\theta\Lambda(\nu_1)+(1-\theta)\Lambda(\nu_2)
 -\frac\gamma2\theta(1-\theta)|\nu_1-\nu_2|^2.
\end{equation*}
In particular, $\Lambda$ is strictly convex on
$\mathcal E$.

It remains to show that the minimizing parameters
$\nu_n$ converge and that their limit minimizes
$\Lambda$. Set $a_n:=\Lambda_n(\nu_n)$.
By the monotonicity in
\eqref{eq:LD-pointwise-Floquet-limit},
\begin{equation*}
 a_n\le a_{n+1}\le\Lambda_{n+1}(0)
 \le\Lambda(0)<\infty.
\end{equation*}
Thus $a_n$ increases to some finite number $a$.
For $n\ge r$, the minimality of $\nu_r$ and
\eqref{eq:Lambda-strong-convex} give
\begin{equation*}
\begin{aligned}
 a_r\le\Lambda_r((1-\theta)\nu_r+\theta\nu_n)\le(1-\theta)a_r+\theta\Lambda_r(\nu_n)
 -\frac\gamma2\theta(1-\theta)|\nu_n-\nu_r|^2.
\end{aligned}
\end{equation*}
Divide by $\theta$ and let $\theta\downarrow0$.
Since $\Lambda_r\le\Lambda_n$ and $\nu_n$ minimizes
$\Lambda_n$, this gives
\begin{equation*}
 \frac\gamma2|\nu_n-\nu_r|^2
 \le\Lambda_r(\nu_n)-a_r
 \le a_n-a_r.
\end{equation*}
Since $a_n$ converges,
the preceding inequality shows that
$\{\nu_n\}$ is Cauchy.

Let $\nu_n\to\nu_*$.
For each fixed $r$, the inequalities
$a_r\le\Lambda_r(\nu_n)\le a_n$ imply
$a_r\le\Lambda_r(\nu_*)\le a=\lim_{k\to\infty} a_k$.
Now let $r\to\infty$ and use
\eqref{eq:LD-pointwise-Floquet-limit} to obtain
$\Lambda(\nu_*)=a<\infty$.
In particular, $\nu_*\in\mathcal E$.
Finally, passing to the limit in
$a_n\le\Lambda_n(\nu)$, for any fixed $\nu$,
shows that $\Lambda(\nu_*)\le\Lambda(\nu)$.
Strict convexity makes the minimizer unique.
This proves (ii).
\end{proof}

\begin{proof}[Proof of Theorem~\ref{cor:general-kernel-large-domain-limit}]
By Lemma~\ref{lem:general-kernel-minimizer}, the minima
$\min_\nu\Lambda_n(\nu)$ converge to
$\min_\nu\Lambda(\nu)$. It remains to show that
$\lambda_l$ has this same limit.
The upper estimate is obtained from a positive test
function on $(-l,l)$; for the lower estimate we use
the test function of Lemma~\ref{thm:large-domain-limit}.

\textbf{Step 1. Upper bound.}
Fix $l>0$, and choose an integer $n>2l$.
Then
$k_{i,n}(x-y,t)=k_i(x-y,t)$ for all
$x,y\in[-l,l]$.
Let $\nu_n$ minimize $\Lambda_n$, and let
$p_n(t)$ be the positive periodic vector associated
with $B_{n,\nu_n}$. Define
$\Psi_n(x,t)=e^{\nu_nx}p_n(t)$.
For $x\in(-l,l)$,
\begin{equation*}
\begin{aligned}
 \int_{-l}^l k_i(x-y,t)\Psi_{n,i}(y,t)\,dy
 &=\Psi_{n,i}(x,t)
   \int_{x-l}^{x+l}k_{i,n}(z,t)e^{-\nu_nz}\,dz\\
 &\le m_{i,n}(\nu_n,t)\Psi_{n,i}(x,t).
\end{aligned}
\end{equation*}
Recall that $\Psi_n(x,t)=e^{\nu_nx}p_n(t)$ and
that $p_n$ satisfies the Floquet equation
\[
-p_n'(t)+B_{n,\nu_n}(t)p_n(t)
=\Lambda_n(\nu_n)p_n(t).
\]
By the definition of $\mathcal L_l$ and the
above integral estimate, for every $i\in\mathbb S$,
	\begin{align*}
		(\mathcal L_l\Psi_n)_i(x,t)
		={}&-\partial_t\Psi_{n,i}
		+c(t)\partial_x\Psi_{n,i}
		+d_i(t)\mathcal K^i_{(-l,l)}[\Psi_{n,i}]
		+(A(t)\Psi_n)_i\\
		\le{}&e^{\nu_nx}\big[
		-p_{n,i}'(t)+\nu_nc(t)p_{n,i}(t)
		+d_i(t)(m_{i,n}(\nu_n,t)-1)p_{n,i}(t)
		+(A(t)p_n(t))_i
		\big]\\
		={}&e^{\nu_nx}
		\big[-p_{n,i}'(t)
		+(B_{n,\nu_n}(t)p_n(t))_i\big]\\
		={}&\Lambda_n(\nu_n)\Psi_{n,i}(x,t).
	\end{align*}
	Consequently,
$\mathcal L_l\Psi_n\le\Lambda_n(\nu_n)\Psi_n$ in $(-l,l)\times\mathbb R$.
Because $\Psi_n\in\mathcal D_l$,  Theorem~\ref{EX-PEV}(ii) gives
$\lambda_l\le\Lambda_n(\nu_n)$.
Letting $n \to \infty$ and using \eqref{eq:LD-truncated-minimizers}, we obtain $\lambda_l\le\min_{\nu\in\mathcal E}\Lambda(\nu)$, and hence
\begin{equation}\label{eq:LD-upper-eigenvalue}
 \limsup_{l\to\infty}\lambda_l
 \le\min_{\nu\in\mathcal E}\Lambda(\nu).
\end{equation}

\textbf{Step 2. Lower bound.}
Fix $n\ge1$ and $\varepsilon>0$.
Let $\mathcal L_{l,n}$ be the 
operator obtained from $\mathcal L_l$ by replacing
$k_i$ with $k_{i,n}$. By Lemma~\ref{thm:large-domain-limit},
for all sufficiently large $l$, there exists
$U_{l,n}\in\mathcal D_{0,l}$ such that
\begin{equation*}
 \mathcal L_{l,n}U_{l,n}
 \ge\bigl(\min_{\nu\in\mathbb R}\Lambda_n(\nu)
       -\varepsilon\bigr)U_{l,n}.
\end{equation*}
Since $0\le k_{i,n}\le k_i$ and $U_{l,n}\ge0$,
we also have $\mathcal L_lU_{l,n}\ge\mathcal L_{l,n}U_{l,n}$.
The second variational characterization in
Theorem~\ref{EX-PEV}(ii) therefore implies
\begin{equation*}
\lambda_l\ge
\min_{\nu\in\mathbb R}\Lambda_n(\nu)-\varepsilon
\text{ for all sufficiently large }l.
\end{equation*}
First let $l\to\infty$, then $\varepsilon\to 0$,
and finally $n\to\infty$. By
\eqref{eq:LD-truncated-minimizers},
\begin{equation}\label{eq:LD-lower-eigenvalue}
 \liminf_{l\to\infty}\lambda_l
 \ge\min_{\nu\in\mathcal E}\Lambda(\nu).
\end{equation}
Combining \eqref{eq:LD-upper-eigenvalue} and
\eqref{eq:LD-lower-eigenvalue} proves the claimed
large-domain limit.
\end{proof}

\section{Persistence criterion for the nonlinear system on bounded domains}

This section relates the sign of the principal eigenvalue to the
long-time behavior of the nonlinear problem \eqref{1.2} on a bounded
interval. The argument has two parts. First, a periodic lower solution
is used to construct a positive periodic solution. The resulting
existence and uniqueness statement is then combined with the
principal eigenfunction to prove Theorem~\ref{Kaiguan}.

The comparison principle used below follows directly
from Lemma~\ref{max-ini-2}.
\begin{lemma}\label{comparsion principle}
Assume (H1)--(H3) and (F). Let $\overline V$, $\underline V$,
$\mathcal D_c\overline V$ and $\mathcal D_c\underline V$
$\in C(\overline Q_{T};\mathbb R^m)$ satisfy
$\bm 0\le \overline V,\underline V\le \mathbf{M}$ 
on $\overline Q_{T}$.
If $\overline V$ and $\underline V$ satisfy
\begin{equation*}
\left\{
\begin{aligned}
& \mathcal D_c \overline V(x,t) \ge   
D(x,t)\mathcal{K}_{(l_{1},l_{2})}[\overline{V}](x,t) + F(x,t,\overline{V}),  &&(x,t)\in Q_{T},\\
& 
\mathcal D_c \underline V(x,t) \le 
D(x,t)\mathcal{K}_{(l_{1},l_{2})}[\underline{V}](x,t) + F(x,t,\underline{V}),  &&(x,t)\in Q_{T},\\
&\overline{V}(l_{2},t) \ge \underline{V}(l_{2},t), &&t \in (0,T],\\
&\overline{V}(x,0) \ge \underline{V}(x,0), &&x\in \overline{\Omega}.
\end{aligned}
\right.
\end{equation*}
Then $\overline{V} \geq \underline{V}$ in $\overline{Q}_{T}$.
\end{lemma}

Indeed, the difference $W=\overline V-\underline V$ satisfies
a cooperative linear differential inequality, because
\[
F(x,t,\overline V)-F(x,t,\underline V)
=\left(\int_0^1\mathrm D_uF
\bigl(x,t,\underline V+\theta(\overline V-\underline V)\bigr)
\,d\theta\right)W.
\]
The matrix in parentheses is cooperative by (F), and the initial
and inflow inequalities for $W$ are nonnegative. Lemma~\ref{max-ini-2}
therefore gives $W\ge\mathbf0$.

\begin{definition}\label{def:periodic-characteristic-barriers}
A periodic lower solution of \eqref{draft-P}  in the characteristic sense
is a function $W$ satisfying
$W,\mathcal D_cW\in C(\overline Q_T;\mathbb R^m)$, $\mathbf0\leq W\leq\mathbf M$ in $\overline Q_T$, and
	\begin{equation*}
		\left\{
		\begin{aligned}
			&\mathcal D_cW(x,t)
			\leq
			D(x,t)\mathcal K_{(l_1,l_2)}[W](x,t)
			+F(x,t,W(x,t)),
			&&(x,t)\in Q_T,\\
			&W(l_{2},t)\leq\mathbf0,
			&&0<t\leq T,\\
			&W(x,0)\leq W(x,T),
			&&x\in[l_1,l_2].
		\end{aligned}
		\right.
	\end{equation*}
\end{definition}
Since $W\ge\mathbf0$ and $W(l_2,t)\le\mathbf0$,
necessarily $W(l_2,t)=\mathbf0$ for $0<t\le T$, and in
particular $W(l_2,T)=\mathbf0$.
Hence $W(\cdot,T)\in\mathbb U_{(l_1,l_2)}$ and
$\mathbf0\leq W(\cdot,T)\leq\mathbf M$ on $\overline\Omega$.
Therefore, $W(\cdot,T)$ can be used as initial data for the initial-value problem.

\begin{theorem}\label{l-U-m}
Assume (H1)--(H3) and (F). If $\underline V$ is a nonzero
continuous lower solution of \eqref{draft-P} satisfying
$\mathbf0\leq\underline V\leq\mathbf M$ on $\overline Q_T$,
then \eqref{draft-P} admits a unique positive $T$-periodic solution $V$
in the characteristic sense, satisfying
$V,\mathcal D_cV\in C(\overline Q_T;\mathbb R^m)$,
and
\begin{equation*}
\underline V\leq V\leq\mathbf M
\text{ in }\overline Q_T,
\quad
V\gg\mathbf0
\text{ in }\Omega\times\mathbb R.
\end{equation*}
\end{theorem}

\begin{proof}
The proof has three parts. The lower solution first generates
an increasing sequence of solutions over successive periods.
We then show that these solutions converge uniformly to a periodic
solution. Finally, strict subhomogeneity gives uniqueness.

For brevity, set
\begin{equation*}
\mathcal F_i[Z](x,t) =d_i(x,t)\mathcal K_\Omega^i[Z_i](x,t) +f_i(x,t,Z(x,t)),
\quad i \in \mathbb{S}.
\end{equation*}

\textbf{Step 1. Construction of an increasing sequence.}
Set $V^0=\underline V$. For $n\ge0$, 
let $V^{n+1}$ be the solution of the following initial problem in the characteristic sense
\begin{equation*}
\left\{
\begin{aligned}
&\mathcal D_cV_i^{n+1}(x,t) =
\mathcal F_i[V^{n+1}], && (x,t)\in Q_T,\\
&V^{n+1}(l_2,t)=\mathbf0, &&0\leq t\leq T,\\
&V^{n+1}(x,0)=V^n(x,T), &&x\in\overline\Omega.
\end{aligned}
\right.
\end{equation*}
Since the kernels have unit mass on $\mathbb R$,
$\mathcal K_\Omega[\mathbf M]\le\mathbf0$. Assumption (F)
also gives $F(x,t,\mathbf M)\le\mathbf0$; hence
$\mathcal F[\mathbf M]\le\mathbf0$. Thus $V^{n} \le \mathbf M$ in $\overline{Q}_{T}$ for $n \ge 0$.
 On the other hand, the lower-solution
inequality, together with
$\underline V(\cdot,0)\le\underline V(\cdot,T)$,
allows comparison first between $V^1$ and $V^0$, and then
between $V^{n+1}$ and $V^n$. Induction yields
\begin{equation*}
\underline V\leq V^n\leq V^{n+1}\leq \mathbf{M}
\enspace \text{in}\enspace \overline Q_T,
\quad n\geq0.
\end{equation*}
For each $(x,t)\in\overline Q_T$, the sequence
$V^n(x,t)$ is componentwise nondecreasing and bounded above by
$\mathbf M$. Therefore its pointwise limit
\begin{align*}
\widetilde V(x,t):=\lim_{n\to\infty}V^n(x,t)
\end{align*}
exists.

\textbf{Step 2. Uniform convergence and the periodic equation.}
We establish uniform convergence of the iterates
using the fact that every backward characteristic
reaches the inflow boundary within a uniformly
bounded time.

We first trace $V^n(x,t)$ back to the inflow boundary, keeping track
of the iterate on each time interval. Recall that $c$ is a positive $T$-periodic function
on $\mathbb R$ and
$\alpha(t)=\int_0^t c(r)\,dr$.
Choose an integer $N\ge1$ such that
$N\alpha(T)>l_2-l_1$.  
Fix $(x,t)\in\overline Q_T$, and follow its characteristic
backward in time:
\begin{align*}
X(r)=x+\alpha(t)-\alpha(r)=x+\int_r^t c(s)\,ds.
\end{align*}
In this parametrization,
$X(t)=x$ and $X'(r)=-c(r)<0$. Therefore, as $r$ decreases,
the characteristic moves to the right and eventually reaches
$x=l_2$. Denote this entry time by $\tau$, so that
$X(\tau)=l_2$. Since
$X(t-NT)=x+N\alpha(T)>l_2$, it follows that
\begin{equation}\label{eq:uniform-entry-time}
t-NT<\tau\leq t .
\end{equation}
	
For $n\geq1$, integration of the characteristic equation for $V^n$ gives
\begin{equation}\label{eq:one-period-characteristic}
V^n(x,t)=
\begin{cases}
\displaystyle\int_\tau^t\mathcal F[V^n](X(r),r)\,dr,
&\tau\geq0,\\
\displaystyle V^{n-1}(X(0),T)
+\int_0^t\mathcal F[V^n](X(r),r)\,dr,
&\tau<0 .
\end{cases}
\end{equation}
In the first case, the backward characteristic reaches the inflow boundary at a nonnegative time, so the representation involves only the current iterate.
In the second case, it reaches time~$0$ first. At this point the
initial condition for the $n$th iterate replaces $V^n(X(0),0)$
by $V^{n-1}(X(0),T)$. We can then repeat the same argument for
the preceding iterate.

Fix $n\geq N+1$ and suppose $\tau<0$. To evaluate
$V^{n-1}(X(0),T)$ in \eqref{eq:one-period-characteristic}, we apply
the characteristic formula to $V^{n-1}$ at $(X(0),T)$. At time $s\in[0,T]$, the position on this characteristic is
\begin{equation*}
X(0)+\int_s^T c(q)\,dq
=X(0)+\int_{s-T}^0 c(u)\,du
=x+\int_{s-T}^t c(u)\,du
=X(s-T),
\end{equation*}
where the first equality follows from $u=q-T$ and the periodicity of $c$, and the second uses $X(0)=x+\int_0^t c(u)\,du$.
If $-T\leq\tau<0$, the characteristic reaches $l_2$ at time
$s=\tau+T$, since $X(s-T)=X(\tau)=l_2$. Using
$V^{n-1}(l_2,\tau+T)=\mathbf0$, we integrate from $\tau+T$ to $T$ and then set $r=s-T$ to obtain
\begin{equation*}
\begin{aligned}
V^{n-1}(X(0),T) =\int_{\tau+T}^T
\mathcal F[V^{n-1}](X(s-T),s)\,ds=\int_\tau^0
\mathcal F[V^{n-1}](X(r),r+T)\,dr .
\end{aligned}
\end{equation*}
If $\tau<-T$, the characteristic is still inside  the interval $(l_1,l_2)$ at time $s=0$, since $X(-T)<X(\tau)=l_2$. We therefore integrate from $0$ to $T$. Using the initial condition
$V^{n-1}(X(-T),0)=V^{n-2}(X(-T),T)$ and the same change of
variables, we obtain
\begin{equation*}
\begin{aligned}
V^{n-1}(X(0),T) &=V^{n-1}(X(-T),0) +\int_0^T\mathcal F[V^{n-1}](X(s-T),s)\,ds\\
&=V^{n-2}(X(-T),T) +\int_{-T}^0\mathcal F[V^{n-1}](X(r),r+T)\,dr.
\end{aligned}
\end{equation*}

After repeating this procedure, the characteristic must reach
the inflow boundary in at most $N$ periods. More precisely,
choose the integer $q\in\{1,\ldots,N\}$ for which
$-qT\le\tau<-(q-1)T$.

For $1\leq j\leq q$ and $\max\{\tau,-jT\}\leq r\leq-(j-1)T$, the $T$-periodicity of $c$ gives
\begin{equation*}
X(-(j-1)T)+\int_{r+jT}^{T}c(\rho)\,d\rho
=X(-(j-1)T)+\int_r^{-(j-1)T}c(\rho)\,d\rho
=X(r).
\end{equation*}
Integrating the characteristic equation for $V^{n-j}$ from
$r+jT$ to $T$, we obtain
\begin{equation*}
V^{n-j}(X(-(j-1)T),T) =V^{n-j}(X(r),r+jT)
+\int_r^{-(j-1)T} \mathcal F[V^{n-j}](X(\rho),\rho+jT)\,d\rho .
\end{equation*}
For $1\leq j<q$, take $r=-jT$. The initial condition gives
\begin{equation*}
V^{n-j}(X(-jT),0)=V^{n-j-1}(X(-jT),T).
\end{equation*}
For $j=q$, take $r=\tau$. Since $X(\tau)=l_2$ and
$0\leq\tau+qT<T$, the inflow boundary condition gives
\begin{equation*}
V^{n-q}(X(\tau),\tau+qT)=V^{n-q}(l_2,\tau+qT)=\mathbf0.
\end{equation*}

At every intermediate time $-jT$, the terminal value of
one iterate is the initial value of the next. At the final time
$\tau$, the value is zero because the characteristic meets
the inflow boundary. Substituting these values successively into
\eqref{eq:one-period-characteristic} gives
\begin{equation*}
\begin{aligned}
V^n(x,t) =&\int_\tau^{-(q-1)T}
\mathcal F[V^{n-q}](X(r),r+qT)\,dr\\
&+\sum_{j=1}^{q-1}\int_{-jT}^{-(j-1)T}
\mathcal F[V^{n-j}](X(r),r+jT)\,dr
+\int_0^t\mathcal F[V^n](X(r),r)\,dr .
\end{aligned}
\end{equation*}

The preceding identity contains different iterates on
different time intervals. For the convergence argument, it is convenient to regard them as one function on the longer time interval $[-NT,T]$. Define $Z^n$ by
\begin{equation}\label{eq:finite-stitching}
Z^n(y,r)=V^{n-j}(y,r+jT),
\qquad -jT\leq r\leq(1-j)T,\quad j=0,\ldots,N .
\end{equation}
At the common endpoints of these time intervals, the
expressions for $Z^n$ agree: for example,
$V^n(y,0)=V^{n-1}(y,T)$, and the same identity holds
for every successive pair. Thus $Z^n$ is continuous.
The coefficients are $T$-periodic, so each piece solves the
same characteristic equation in its own time interval.
Combining the integral formulas and using the zero inflow value
gives, for $(x,t)\in\overline Q_T$,
\begin{equation}\label{eq:stitched-characteristic-at-t}
Z^n(X(s),s) =\int_\tau^s \mathcal F[Z^n](X(r),r)\,dr,
\qquad \tau\leq s\leq t.
\end{equation}
Note that at $s=t$, the left-hand side is $Z^n(x,t)=V^n(x,t)$.

We now estimate the difference between two such stitched
functions. Recall that $|z|=\sum_i|z_i|$ for $z\in\mathbb R^m$.
Since the nonlocal integral operator is bounded
and $F$ is uniformly Lipschitz continuous in $u$
on $[\mathbf0,\mathbf M]$, there exist constants
$B,L\ge0$, independent of $n$ and $p$, such that
\begin{equation}\label{eq:source-difference-bound}
|\mathcal F[Z](x,r)-\mathcal F[Y](x,r)|
\leq B\int_\Omega|Z(y,r)-Y(y,r)|\,dy
+L|Z(x,r)-Y(x,r)|
\end{equation}
whenever $\mathbf0\leq Z,Y\leq\mathbf{M}$. Fix $n\geq N+1$, $p\geq1$ and $(x,t)\in\overline Q_T$, and put
\begin{equation*}
\delta(s) =|Z^{n+p}(X(s),s)-Z^n(X(s),s)|,
\quad \tau\leq s\leq t.
\end{equation*}
Both stitched functions vanish at the same inflow
entry point $(l_2,\tau)$. Consequently, no initial-value term
appears when their characteristic formulas are subtracted.
Using \eqref{eq:source-difference-bound} then gives
\begin{equation*}
\begin{aligned}
\delta(s) &\leq \int_\tau^s
\left|\mathcal F[Z^{n+p}](X(r),r)
-\mathcal F[Z^n](X(r),r)\right|\,dr\\
&\leq B\int_\tau^s\int_\Omega
|Z^{n+p}(y,r)-Z^n(y,r)|\,dy\,dr
+L\int_\tau^s\delta(r)\,dr,
\qquad \tau\leq s\leq t.
\end{aligned}
\end{equation*}
Gronwall's inequality and \eqref{eq:uniform-entry-time} now give
\begin{equation*}
\begin{aligned}
|V^{n+p}(x,t)-V^n(x,t)|=\delta(t)&\leq
B\int_\tau^t e^{L(t-r)}
\int_\Omega
|Z^{n+p}(y,r)-Z^n(y,r)|\,dy\,dr\\
&\leq Be^{LNT}\int_{-NT}^T\int_\Omega
|Z^{n+p}(y,r)-Z^n(y,r)|\,dy\,dr.
\end{aligned}
\end{equation*}
Taking the supremum over
$(x,t)\in\overline Q_T$ and changing variables from each
interval $[-jT,(1-j)T]$ back to $[0,T]$ gives
\begin{equation}\label{eq:uniform-cauchy-estimate}
\begin{aligned}
\|V^{n+p}-V^n\|_{C(\overline Q_T)} &\leq Be^{LNT}\sum_{j=0}^N
\|V^{n+p-j}-V^{n-j}\|_{L^1(Q_T)}\\
&\leq (N+1)Be^{LNT} \|\widetilde V-V^{n-N}\|_{L^1(Q_T)}.
\end{aligned}
\end{equation}
To see the last inequality, note that for $0\le j\le N$,
monotonicity of the iterates gives
\[
\mathbf0\le V^{n+p-j}-V^{n-j}
\le\widetilde V-V^{n-j}
\le\widetilde V-V^{n-N}.
\]
Since $V^{n-N}\to\widetilde V$ pointwise and
$\mathbf0\le\widetilde V-V^{n-N}\le\mathbf M$,
the dominated convergence theorem gives
$\|\widetilde V-V^{n-N}\|_{L^1(Q_T)}\to0$.
As the estimate is independent of $p$,
$\{V^n\}$ is uniformly Cauchy on $\overline Q_T$.
Hence $V^n$ converges uniformly to the continuous
function $\widetilde V$.

Passing to the limit in the initial and boundary conditions gives
\begin{equation*}
\widetilde V(x,0)=\widetilde V(x,T)
\text{ for }x\in\overline\Omega,
\qquad
\widetilde V(l_2,t)=\mathbf0
\text{ for }t\in[0,T].
\end{equation*}
Extend $\widetilde V$ periodically in time. Since $N$ is fixed,
\eqref{eq:finite-stitching} and the uniform convergence of $V^n$ on $\overline Q_T$ give $Z^n\to\widetilde V$ uniformly on
$S=\overline\Omega\times[-NT,T]$. By
\eqref{eq:source-difference-bound},
$\mathcal F[Z^n]\to\mathcal F[\widetilde V]$ uniformly on $S$.
Uniform convergence also permits passage to the limit in
the characteristic integral, because $\mathcal F$ is uniformly
continuous on the bounded range under consideration.
Letting $n\to\infty$ in
\eqref{eq:stitched-characteristic-at-t} gives
\begin{equation}\label{IR-5}
\widetilde V(X(s;x,t),s) =\int_{\tau(x,t)}^s
\mathcal F[\widetilde V](X(r;x,t),r)\,dr,
\qquad \tau(x,t)\leq s\leq t .
\end{equation}
The characteristic formula extends to all time translates
by $T$-periodicity. Its integrand is continuous; hence the
fundamental theorem of calculus gives the derivative along the
characteristic. Recalling that $X'(s)=-c(s)$, this derivative
is precisely $\mathcal D_c=\partial_t-c(t)\partial_x$, so
\begin{equation*}
\mathcal D_c\widetilde V(x,t) =D(x,t)\mathcal K_\Omega[\widetilde V](x,t)
+F(x,t,\widetilde V(x,t)), \qquad (x,t)\in\Omega\times\mathbb R.
\end{equation*}
The right-hand side extends continuously to $\overline Q_T$. Thus $\widetilde V$ is a continuous $T$-periodic solution of \eqref{draft-P} in the characteristic sense.
Moreover, $\widetilde V$ satisfies
\begin{equation*}
\mathcal D_c\widetilde V
=
D(x,t)\mathcal K_\Omega[\widetilde V]
+B_{\widetilde V}(x,t)\widetilde V,
\end{equation*}
where $B_{\widetilde V}(x,t):=
\int_0^1D_uF(x,t,\theta\widetilde V(x,t))d\theta$. 
By (F)(iv),  $B_{\widetilde V}(x,t)$ is cooperative, and $B_{\widetilde V}(x_0,t)$ is irreducible for every $t$. Since
$\widetilde V\ge\underline V\not\equiv\mathbf0$,
Lemma~\ref{max-coo} and the $T$-periodicity of
$\widetilde V$ imply that
$\widetilde V\gg\mathbf0$ on
$[l_1,l_2)\times\mathbb R$.

\textbf{Step 3. Uniqueness of the positive periodic solution.}

Let $\widehat V$ be any nontrivial nonnegative continuous
$T$-periodic solution of \eqref{draft-P} in the characteristic sense satisfying $\mathbf0\leq\widehat V\leq\mathbf M$ on $\overline Q_T$.

Lemma~\ref{max-coo} applies to both solutions and gives
\begin{equation*}
\widetilde V,\widehat V\gg\mathbf0
\text{ in }[l_1,l_2)\times\mathbb R.
\end{equation*}
By Lemmas~\ref{LEM-lower-boundary} and~\ref{LEM-upper-boundary},
there exist constants $\underline C_i,\overline C_i>0$ such that
\begin{equation}\label{UL}
\widetilde V_i(x,t)\geq\underline C_i(l_2-x),
\,
\widehat V_i(x,t)\leq\overline C_i(l_2-x)
\text{ for } (x,t)\in \overline Q_T
\text{ and } i\in\mathbb S.
\end{equation}
	
Set
\begin{equation*}
\overline s :=
\sup\left\{ s>0:\widetilde V\geq s\widehat V
\text{ on }\overline Q_T \right\}.
\end{equation*}
The lower bound for $\widetilde V$ and the upper bound for
$\widehat V$ in \eqref{UL} show that $\overline s>0$.
Conversely, evaluation at any fixed interior point, where both
solutions are positive, shows that $\overline s<\infty$.
Moreover,
$\widetilde V\geq\overline s\widehat V$
in $\overline Q_T$.

Suppose that $\overline s<1$, and set
$W=\widetilde V-\overline s\widehat V\geq\mathbf0$.
Define
\begin{equation*}
B(x,t)=
\int_0^1
D_uF\left(x,t,\overline s\widehat V(x,t)+\theta W(x,t)
\right)\,d\theta .
\end{equation*}
By (F), $B$ is cooperative and $B(x_0,t)$ is irreducible.
Since $0<\overline s<1$ and $\widehat V\gg\mathbf0$
in the interior, strict subhomogeneity gives $F(x,t,\overline s\widehat V)-\overline sF(x,t,\widehat V)\geq\mathbf0$, with the difference nonzero at every interior point.
Subtracting the two equations and using the mean-value formula
for $F(x,t,\widetilde V)-F(x,t,\overline s\widehat V)$ gives
\begin{equation*}
\begin{aligned}
\mathcal D_cW-D(x,t)\mathcal K_\Omega[W]-B(x,t)W
=F(x,t,\overline s\widehat V)
-\overline sF(x,t,\widehat V)>\mathbf0
\enspace \text{in} \enspace Q_T.
\end{aligned}
\end{equation*}
The strict inequality shows that $W$ cannot vanish
identically. The cooperative strong maximum principle (Lemma \ref{max-coo}) therefore
gives
$W\gg\mathbf0$ in $[l_1,l_2)\times\mathbb R$.
By Lemma~\ref{LEM-lower-boundary}, there exist constants
$\omega_i>0$ such that
\begin{equation*}
W_i(x,t)\geq\omega_i(l_2-x)
\geq\frac{\omega_i}{\overline C_i}\widehat V_i(x,t)
\text{ in }\overline Q_T.
\end{equation*}
Therefore, for some $\delta>0$,
\begin{equation*}
\widetilde V
=W+\overline s\widehat V
\geq(\overline s+\delta)\widehat V
\text{ in }\overline Q_T.
\end{equation*}
This contradicts the definition of $\overline s$.
Thus $\overline s\ge1$, and hence
$\widetilde V\geq\widehat V$ on $\overline Q_T$.
Interchanging the two solutions gives the reverse inequality.
Therefore $\widetilde V=\widehat V$ on $\overline Q_T$.
This completes the proof.
\end{proof}

\begin{proof}[Proof of Theorem~\ref{Kaiguan}]
Throughout the proof, write
$\lambda=\lambda_{(l_1,l_2)}$.

Let $\varphi$ be the positive $T$-periodic principal eigenfunction associated with
$\lambda$. By Lemmas~\ref{LEM-upper-boundary} and~\ref{LEM-lower-boundary}, there exist $a_\varphi,b_\varphi>0$ such that
\begin{equation}\label{UD}
a_\varphi(l_2-x)\leq\varphi_i(x,t)\leq b_\varphi(l_2-x)
\text{ on }\overline Q_T,\quad i\in\mathbb S.
\end{equation}
By Lemma \ref{EX-draft},
for each $V_0\in\mathbb U_{(l_1,l_2)}$
with $\mathbf0\le V_0\le\mathbf M$,
the problem \eqref{1.2}  admits  a unique global solution, denoted by $V(x,t)=V(x,t;V_0)$,
satisfying $\mathbf0\leq V\leq\mathbf M$ in $\overline\Omega\times[0,\infty)$.

\textbf{(i) $\lambda>0$. }
In this case, the key is that a small multiple of
the principal eigenfunction is a periodic lower solution.
By \eqref{UD}, $|\sum_j\varphi_j|\leq C\varphi_i$ on $\overline Q_T$ for every $i\in\mathbb S$.
Since $F(x,t,\mathbf0)=\mathbf0$ and
$\mathrm D_uF(x,t,\mathbf0)=A(x,t)$,
Taylor's formula gives
$F(x,t,\varepsilon\varphi)
=\varepsilon A(x,t)\varphi+o(\varepsilon|\varphi|)$
uniformly in $(x,t)\in\overline Q_T$.
The ratio $|\varphi|/\varphi_i$ is uniformly bounded
by \eqref{UD}. Thus, for sufficiently small $\varepsilon>0$,
\begin{equation*}
\left|f_i(x,t,\varepsilon\varphi)
-\varepsilon(A(x,t)\varphi)_i\right|
\le \frac{\lambda}{2}\varepsilon\varphi_i(x,t),
\quad i\in\mathbb S.
\end{equation*} 
It follows that 
\begin{equation*}
\begin{aligned}
&\mathcal D_c(\varepsilon\varphi)
-D(x,t)\mathcal K_\Omega[\varepsilon\varphi]
-F(x,t,\varepsilon\varphi)\\
=&
\varepsilon \mathbf{D}_{u}F(x,t,\mathbf0)\varphi
-\lambda\varepsilon\varphi
-F(x,t,\varepsilon\varphi)\\
\leq&
-\tfrac{\lambda}{2}\varepsilon\varphi\leq\mathbf0
\quad\text{in }Q_T.
\end{aligned}
\end{equation*}
Thus, for small $\varepsilon>0$, $\varepsilon\varphi$ is a nonzero periodic lower solution of \eqref{draft-P} in $[\mathbf0,\mathbf M]$.
By Theorem~\ref{l-U-m}, the problem \eqref{draft-P} has a unique nontrivial nonnegative continuous $T$-periodic
solution $V^*$ in $[\mathbf0,\mathbf M]$.
Each component of $V^*$ is strictly positive on
$[l_1,l_2)\times\mathbb R$.

To prove attraction, fix nonzero initial data $V_0$.
We first show that the solution at time $T$ is comparable with
the periodic state all the way up to the inflow boundary.
Write
$F(x,t,V)=B(x,t)V$ with
$B(x,t)=\int_0^1\mathbf D_uF(x,t,rV(x,t))\,dr$ in $\overline Q_T$.
By cooperativity and continuity at $r=0$, $B(x_0,t)$ is irreducible for $t\in[0,T]$.
The strong maximum principle for the initial-value problem,
Lemma~\ref{max-ini-2}, then gives
$V_i(x,T)>0$ for $x\in[l_1,l_2)$ and $i\in\mathbb S$.

At the inflow boundary, $V(l_2,T)=\mathbf0$ and
$F(l_2,T,\mathbf0)=\mathbf0$. Under the regularity
established in Section~2, the equation can be evaluated by taking
the interior limit $x\to l_2$. It follows that
\begin{equation*}
\mathcal D_cV_i(l_2,T)
=
d_i(l_2,T)\int_{l_1}^{l_2} k_i(l_2-y,T)V_i(y,T)\,dy>0,\quad i\in\mathbb S.
\end{equation*}
For $x<l_2$ close to $l_2$, trace the backward characteristic
$X(r)=x+\int_r^T c(s)\,ds$ to its entry time $\tau$, defined
by $X(\tau)=l_2$. Since
$l_2-x=\int_\tau^T c(s)\,ds$, we have $\tau\to T$
as $x\to l_2$. Integrating the equation along this
characteristic, starting from $V_i(l_2,\tau)=0$, yields
\begin{equation*}
\frac{V_i(x,T)}{l_2-x}
=\frac{\int_\tau^T\mathcal D_cV_i(X(r),r)\,dr}{\int_\tau^T c(r)\,dr}
\longrightarrow\frac{\mathcal D_cV_i(l_2,T)}{c(T)}>0
\enspace\text{as} \enspace x\to l_2,\, i\in\mathbb S.
\end{equation*}
This positive limit controls $V_i(x,T)/(l_2-x)$ near the
inflow boundary. Away from the boundary, the ratio has a positive
minimum and finite maximum by continuity and strict positivity.
Thus there exist $a_0,b_0>0$ such that
\begin{equation*}
a_0(l_2-x)\leq V_i(x,T)\leq b_0(l_2-x),
\quad x\in\overline\Omega, i\in\mathbb S.
\end{equation*}
Together with Lemmas~\ref{LEM-upper-boundary}
and~\ref{LEM-lower-boundary} applied to $V^*$, this gives $0<\theta<1<R$ such that
\begin{equation*}
\theta V^*(x,0)\leq V(x,T)\leq RV^*(x,0)
\enspace \text{for all}\enspace x\in\overline\Omega.
\end{equation*}

The factors $\theta<1<R$ yield a lower and an upper
initial datum. Define
\begin{equation*}
\eta^-(x)=\theta V^*(x,0),
\qquad
\eta_i^+(x)=\min\{M_i,RV_i^*(x,0)\},
\quad i\in\mathbb S,
\end{equation*}
and let $Z^\pm(x,t)=V(x,t;\eta^\pm)$.
Then $\eta^\pm\in\mathbb U_{(l_1,l_2)}$,
$\mathbf0\le\eta^\pm\le\mathbf M$, and
$\eta^-\leq V(\cdot,T)\leq\eta^+$ in $\overline\Omega$.
By the comparison principle and periodicity,
\begin{equation}\label{SX}
Z^-(x,t)\leq V(x,T+t)\leq Z^+(x,t),
\quad x\in\overline\Omega,\, t\geq0.
\end{equation}
Note that
subhomogeneity gives $F(x,t,\theta V^*)\ge
\theta F(x,t,V^*)$ and
$F(x,t,R^{-1}Z^+)\ge R^{-1}F(x,t,Z^+)$.
Hence $\theta V^*$ and $R^{-1}Z^+$ are subsolutions in
$\Omega\times(0,\infty)$ with zero inflow condition. Since
$R^{-1}\eta^+\leq V^*(\cdot,0)$ on $\overline\Omega$, comparison gives
\begin{equation*}
\theta V^*\leq Z^-\leq V^*\leq Z^+
\leq\min\{\mathbf M,RV^*\}
\enspace\text{in}\enspace \overline\Omega\times[0,\infty).
\end{equation*}
In particular, $Z^-(\cdot,T)\geq\eta^-$ and $Z^+(\cdot,T)\leq\eta^+$ in
$\overline\Omega$. For $Z_n^\pm(x,t):=Z^\pm(x,nT+t)$ in $\overline Q_T$, comparison and periodicity give
\begin{equation*}
\theta V^* \leq Z_n^-
\leq Z_{n+1}^-
\leq V^*
\leq Z_{n+1}^+
\leq Z_n^+
\leq\mathbf M
\enspace \text{in} \enspace \overline Q_T,\, n\geq0.
\end{equation*}
The two sequences are monotone and bounded. The uniform
Cauchy estimate \eqref{eq:uniform-cauchy-estimate} in the
proof of Theorem~\ref{l-U-m} applies to both (for the decreasing
sequence, apply the estimate to the absolute differences).
Consequently,
$Z_n^\pm\to V^\pm$ uniformly in $\overline Q_T$, where $V^\pm$ are
continuous $T$-periodic solutions of \eqref{draft-P}.
Since $V^\pm\geq\theta V^*$ on $\overline Q_T$, both limits are nontrivial.
The uniqueness in Theorem~\ref{l-U-m} therefore gives $V^-=V^+=V^*$.
It then follows from \eqref{SX} that 
\begin{equation*}
\lim_{n\to\infty}
\sup_{(x,t)\in\overline Q_T}
\left|V(x,nT+t;V_0)-V^*(x,t)\right|=0.
\end{equation*}

\textbf{(ii) $\lambda<0$.}
For $\mathbf0\le u\le\mathbf M$,
subhomogeneity gives
$F(x,t,su)\ge sF(x,t,u)$ for $0<s<1$.
After dividing by $s$ and letting $s \to 0$,
differentiability at zero gives

$F(x,t,u)\leq A(x,t)u$, where $A(x,t)=\mathbf D_uF(x,t,\mathbf0)$. Hence
\begin{equation*}
\mathcal D_cV
\leq D(x,t)\mathcal K_\Omega[V]+A(x,t)V
\enspace \text{in} \enspace \Omega\times(0,\infty).
\end{equation*}
The upper boundary estimate for $V(\cdot,T)$ derived in (i)
does not require $\lambda>0$. Combining that estimate with
the lower bound for $\varphi$ in \eqref{UD}, we can choose
$K>0$ such that
\begin{equation}\label{UD1}
V(x,T)\leq K\varphi(x,0)
\enspace \text{in}\enspace \overline\Omega.
\end{equation}
Define $E(x,t)=Ke^{\lambda t}\varphi(x,t)$.
Then we have
\begin{equation*}
\mathcal D_cE
=
Ke^{\lambda t}(\mathcal D_c\varphi+\lambda\varphi)
=
D(x,t)\mathcal K_\Omega[E]+A(x,t)E
\enspace\text{in}\enspace \Omega\times(0,\infty).
\end{equation*}
By \eqref{UD1}, periodicity and the linear comparison principle,
\begin{equation*}
\mathbf0\leq V(x,T+t)\leq Ke^{\lambda t}\varphi(x,t),
\quad x\in\overline\Omega,\, t\geq0.
\end{equation*}
Since $\varphi$ is bounded and $V\leq\mathbf M$ on $\overline\Omega\times[0,T]$,
\begin{equation*}
\mathbf0\leq V(x,t)\leq Ce^{\lambda t}\mathbf1,
\quad x\in\overline\Omega,\,t\geq0,
\end{equation*}
for some $C>0$.

\textbf{(iii) $\lambda=0$.}
 First, we show that a nonzero
nonnegative periodic solution cannot exist. Suppose, to the contrary, that $\widehat V$ is a nontrivial nonnegative periodic solution in
$[\mathbf0,\mathbf M]$ in the characteristic sense.
Its upper boundary estimate and \eqref{UD} give
\begin{equation*}
k_*:=\inf\left\{
k>0:\widehat V\leq k\varphi\text{ on }\overline Q_T
\right\}\in(0,\infty).
\end{equation*}
The boundary estimates make $k_*$ finite and positive.
By continuity, the defining inequality remains valid at the infimum. Set $W=k_*\varphi-\widehat V\ge\mathbf0$.
Then $W(l_2,t)=\mathbf0$ for $t\in[0,T]$,
$W(x,0)=W(x,T)$ for $x\in\overline\Omega$, and $W$ satisfies
\begin{equation}\label{eq-W}
\mathcal D_cW-D(x,t)\mathcal K_\Omega[W]-A(x,t)W
=
A(x,t)\widehat V-F(x,t,\widehat V)
\geq\mathbf0\enspace \text{in}\enspace Q_T.
\end{equation}
If $W\equiv\mathbf0$, then $\widehat V=k_*\varphi\gg\mathbf0$ in $\Omega\times\mathbb R$. 
Applying $F(x,t,u)\leq A(x,t)u$ from (ii) with
$u=\widehat V/2$ and using strict subhomogeneity gives
\begin{equation*}
A(x,t)\widehat V \geq 2F(x,t,\widehat V/2)>F(x,t,\widehat V) \enspace \text{in} \enspace \Omega\times\mathbb R.
\end{equation*}
This contradicts the equality in \eqref{eq-W} when
$W\equiv\mathbf0$. Thus $W\not\equiv\mathbf0$.
The cooperative strong maximum principle (Lemma \ref{max-coo}) now yields
\begin{equation*}
W_i(x,t)>0
\enspace\text{in}\enspace[l_1,l_2)\times\mathbb R,\,
i\in\mathbb S.
\end{equation*}
Applying Lemma~\ref{LEM-lower-boundary} to $W$, we obtain $\omega>0$ such that
\begin{equation*}
W_i(x,t)\geq\omega(l_2-x)
\geq\frac{\omega}{b_\varphi}\varphi_i(x,t)
\enspace \text{in} \enspace \overline Q_T,\,i\in\mathbb S.
\end{equation*}
Choose $\varepsilon=\min\{k_*/2,\omega/b_\varphi\}>0$.
Then $\widehat V\leq(k_*-\varepsilon)\varphi$ in $\overline Q_T$, contrary to
the definition of $k_*$.
Thus $\mathbf0$ is the only nonnegative $T$-periodic solution of \eqref{draft-P} in $[\mathbf0,\mathbf M]$.

It remains to prove convergence of arbitrary solutions to
zero. Fix admissible initial data $V_0$.
The same boundary comparison used to obtain \eqref{UD1}
remains valid when $\lambda=0$.
Indeed, the upper boundary estimate for $V(\cdot,T)$
and the lower boundary estimate for $\varphi$
give a constant $K>0$ such that
$V(\cdot,T)\le K\varphi(\cdot,0)$.
For $x\in\overline\Omega$, set
\begin{equation*}
\eta_i(x)=\min\{M_i,K\varphi_i(x,0)\},
\quad i\in\mathbb S,
\end{equation*}
and let $Z(x,t)=V(x,t;\eta)$.
Since $V(\cdot,T)\leq\eta$ on $\overline\Omega$, the comparison principle gives
\begin{equation*}
\mathbf0\leq V(x,T+t)\leq Z(x,t)\leq\mathbf M,
\qquad x\in\overline\Omega,\quad t\geq0.
\end{equation*}
By $F(x,t,u)\leq A(x,t)u$ from (ii), $Z$ satisfies
\begin{equation*}
\mathcal D_cZ\leq D(x,t)\mathcal K_\Omega[Z]+A(x,t)Z
\text{ in }\Omega\times(0,\infty).
\end{equation*}
Since $\lambda=0$, $K\varphi$ satisfies the corresponding equality. Also, $Z(l_2,t)=K\varphi(l_2,t)=\mathbf0$ for $t\geq0$ and $Z(\cdot,0)=\eta\leq K\varphi(\cdot,0)$ on $\overline\Omega$.
The comparison principle therefore gives $Z\leq K\varphi$ on $\overline\Omega\times[0,\infty)$. Consequently,
$Z_i(x,T)\leq\min\{M_i,K\varphi_i(x,T)\}=\eta_i(x)$
for $x\in\overline\Omega$ and $i\in\mathbb S$.
	
For $Z_n(x,t):=Z(x,nT+t)$ on $\overline Q_T$, comparison and periodicity give $\mathbf0\leq Z_{n+1}\leq Z_n$ on
$\overline Q_T$ for $n\geq0$. 
The same uniform Cauchy argument used in Step~2 of the
proof of Theorem~\ref{l-U-m} applies to this decreasing
sequence. Since the sequence is bounded and converges
pointwise, the dominated convergence theorem shows that
the corresponding right-hand side of
\eqref{eq:uniform-cauchy-estimate} tends to zero.
Hence $\{Z_n\}$ converges uniformly
in $\overline Q_T$ to a nonnegative $T$-periodic solution
of \eqref{draft-P} in $[\mathbf0,\mathbf M]$.
Since $\mathbf0$ is the only such solution, the limit equals $\mathbf0$. Since
\begin{equation*}
\mathbf0\leq V(x,(n+1)T+t)\leq Z_n(x,t) \enspace \text{in} \enspace \overline Q_T 
\enspace \text{for} \enspace n\geq0,
\end{equation*}
we obtain $V(\cdot,t)\to\mathbf0$ uniformly in $\overline\Omega$ as $t\to\infty$.

\textbf{(iv) Limit on expanding intervals.}
Assume now that $D$ and $F$ are independent of $x$ and that
$\lambda_\infty>0$. The main point is to show that the
positive periodic state on $(-l,l)$ increases with $l$ and
approaches a function independent of $x$.
Write $I_l=(-l,l)$. Let $\Lambda_c$ denote the Floquet exponent defined in
\eqref{FLOQUENT} for the given advection coefficient $c(t)$. By the large-domain limit formula,
	\begin{equation*}
		\Lambda_0(0)=\Lambda_c(0)
		\ge \inf_{\nu\in\mathbb R}\Lambda_c(\nu)
		=\lambda_\infty>0.
	\end{equation*}
	Since $m_i(0,t)=1$, the linearization of $v'=F(t,v)$ at $\mathbf0$
	has principal Floquet multiplier $e^{T\Lambda_0(0)}>1$.
	Hence, by
	\cite[Theorem~3.1.2]{Zhao}, the ODE admits a unique positive
	$T$-periodic solution $U^*$ in $[\mathbf0,\mathbf M]$.
	
	Since $\lambda_{(-l,l)}\to\lambda_\infty>0$, there exists $l_0>0$
	such that $\lambda_{(-l,l)}>0$ for every $l\ge l_0$. By (i), \eqref{draft-P}
	on $I_l$ admits a unique positive $T$-periodic characteristic
	solution $V_l^*:=V_{(-l,l)}^*$ in $[\mathbf0,\mathbf M]$.

\textbf{Step 1. Monotonicity with respect to the interval.}
Fix $l_0\le l<L$. 
Since $V_L^*\gg\mathbf0$ on $\overline I_l\times[0,T]$, we have
\begin{equation*}
s_*:=\max\{s\in[0,1]:sV_l^*\le V_L^*
\text{ in }\overline I_l\times[0,T]\}>0.
\end{equation*}
Suppose that $s_*<1$ and set $W:=V_L^*-s_*V_l^*\ge\mathbf0$.
For $0<s_*<1$, subhomogeneity gives
$F(t,s_*V_l^*)\ge s_*F(t,V_l^*)$.
By the mean-value formula and cooperativity, the off-diagonal
terms are nonnegative. Since $\mathrm D_uF$ is bounded,
there exists a constant $C_F\ge0$ such that
$f_i(t,V_L^*)-f_i(t,s_*V_l^*)\ge -C_FW_i$.
Consequently,
in $I_l\times\mathbb R$,
\begin{equation*}
\begin{aligned}
\mathcal D_cW_i =&d_i(t)\mathcal K_{I_l}^i[W_i]
+f_i(t,V_L^*)-s_*f_i(t,V_l^*)+d_i(t)\int_{I_L\setminus I_l}
k_i(x-y,t)V_{L,i}^*(y,t)\,dy\\
\ge&d_i(t)\mathcal K_{I_l}^i[W_i] +f_i(t,V_L^*)-f_i(t,s_*V_l^*)
\ge -CW_i,
\end{aligned}
\end{equation*}
where $C=C_F+\max_{i\in\mathbb S,\;t\in[0,T]}d_i(t)$. For any $(x,t)\in\overline I_l\times\mathbb R$, the
backward characteristic reaches $x=l$ at a time $\tau \le t$
with
$0\le t-\tau\le2l/c_{\min}$, where
$c_{\min}:=\min_{r\in[0,T]}c(r)>0$.
Along this characteristic, $\mathcal D_cW_i\ge-CW_i$;
thus $r\mapsto e^{Cr}W_i(X(r),r)$ is nondecreasing.
At the entry point,
$W_i(l,\tau)=V_{L,i}^*(l,\tau)>0$.
Integration therefore gives
\begin{equation*}
W_i(x,t)\ge e^{-2Cl/c_{\min}}
\min_{r\in[0,T]}V_{L,i}^*(l,r)>0
\enspace \text{in}\enspace \overline I_l\times[0,T].
\end{equation*}
	As $V_l^*$ is bounded, $W\ge\varepsilon V_l^*$ for some
	$0<\varepsilon<1-s_*$, which contradicts the definition of $s_*$.
	Thus $V_l^*\le V_L^*$ on $\overline I_l\times\mathbb R$.
	
\textbf{Step 2. Identification of the spatially homogeneous limit.}
The monotonicity just proved and the uniform bound
$V_l^*\le\mathbf M$ show that
\begin{align*}
q(x,t):=\lim_{l\to\infty}V_l^*(x,t)
\end{align*}
exists for each
fixed $(x,t)$, once $l>|x|$. To show that the limit is spatially
constant, compare solutions on translated intervals.
Fix $h\in\mathbb R$.
By spatial homogeneity,
$V_l^*(x-h,t)$ solves the same problem on $h+I_l$. Whenever
	$L>l+|h|$, we have $\overline{h+I_l}\subset I_L$, so the preceding
	comparison applies on $h+I_l$, with inflow boundary $x=h+l$. Hence
	\begin{equation*}
		V_l^*(x-h,t)\le V_L^*(x,t)
		\enspace \text{in} \enspace \overline{h+I_l}\times\mathbb R.
	\end{equation*}
	Taking $L=l+|h|+1$ and letting $l\to\infty$ gives
	$q(x-h,t)\le q(x,t)$. Replacing $(x,h)$ by $(x-h,-h)$ gives the
	reverse inequality. Thus $q(x,t)=q(0,t)=:q(t)$ for all $x,t$.
	Moreover, $q$ is $T$-periodic and $q(t)\ge V_{l_0}^*(0,t)\gg\mathbf0$.

It remains to prove that the spatially constant limit $q(t)$
satisfies the periodic ODE $q'=F(t,q)$. Fix $s<t$ and
$x\in\mathbb R$. Along the characteristic through $(x,t)$,
write $X(r)=x+\int_r^t c(\xi)\,d\xi$, $s\le r\le t$.
For all sufficiently large $l$, the characteristic segment
	$\{X(r):s\le r\le t\}$ lies in $I_l$. Consequently,
	\begin{equation}\label{eq:characteristic-formula-Vl}
		\begin{aligned}
			V_{l,i}^*(x,t)-V_{l,i}^*(X(s),s)
			=&\int_s^t d_i(r)
			\left[\int_{I_l} k_i(X(r)-y,r)V_{l,i}^*(y,r)\,dy
			-V_{l,i}^*(X(r),r)\right]dr\\
			&+\int_s^t f_i(r,V_l^*(X(r),r))\,dr.
		\end{aligned}
	\end{equation}
	For each fixed $r$, the zero extensions of $V_l^*(\cdot,r)$ increase
	pointwise to $q(r)$. By monotone convergence and the normalization
	of $k_i$,
	\begin{equation*}
		\lim_{l\to\infty}
		\int_{I_l} k_i(X(r)-y,r)V_{l,i}^*(y,r)\,dy
		=\int_{\mathbb R}k_i(X(r)-y,r)q_i(r)\,dy=q_i(r).
	\end{equation*}
	Since $\mathbf0\le V_l^*\le\mathbf M$, the integrands in
	\eqref{eq:characteristic-formula-Vl} are bounded uniformly in $l$.
	Dominated convergence therefore yields
	\begin{equation}\label{eq:q-integral-ode}
		q(t)-q(s)=\int_s^t F(r,q(r))\,dr.
	\end{equation}
	Thus $q$ is locally absolutely continuous. Since $F$ is continuous,
	$q\in C^1(\mathbb R;\mathbb R^m)$ and $q'=F(t,q)$.
	By uniqueness of the positive $T$-periodic solution, $q=U^*$.

Fix $R>0$ and let $l_n>R$ increase to infinity.
Since $V_{l_n}^*\uparrow U^*$ pointwise and
all these functions are continuous on
$[-R,R]\times[0,T]$, Dini's theorem
yields uniform convergence along $l_n$,
componentwise.
Monotonicity in $l$ then extends this
convergence to all $l\to\infty$, proving (iv).

\end{proof}

\section{Spreading dynamics of the nonlinear system on the whole line}
This section proves the spreading-speed Theorem
\ref{thm:spread-periodic-cooperative}. We first determine the
sign of the large-domain spectral quantity
$\lambda_\infty(c)$ as the speed $c$ varies. Positive values
allow persistence in a moving frame, while negative values lead
to exponential supersolutions outside the spreading region.

Recall the Floquet exponent $\Lambda_0$ defined for the
linearization $A(t)=\mathrm D_uF(t,\mathbf0)$ with zero
advection. As in the Introduction, $\Lambda_0(\nu)=+\infty$
when the corresponding exponential moment is not admissible.
The candidate leftward and rightward speeds are

\begin{equation*}
	c_*^-
	=
	\sup_{\mu<0}\frac{\Lambda_0(-\mu)}{\mu},
	\qquad
	c_*^+
	=
	\inf_{\mu>0}\frac{\Lambda_0(-\mu)}{\mu}.
\end{equation*}

\begin{proposition}\label{prop-sign}
Assume that $\Lambda_0(0)>0$.
Then $-\infty\le c_*^-<c_*^+\le+\infty$.
Moreover, the mapping $c\mapsto\lambda_\infty(c)$ is finite, concave and continuous on $\mathbb R$, and, for every $c\in\mathbb R$,
\begin{equation*}
\lambda_\infty(c)
		\begin{cases}
			>0, & c_*^-<c<c_*^+,\\[1mm]
			=0, & c\in\{c_*^-,c_*^+\}\cap\mathbb R,\\[1mm]
			<0, & c<c_*^- \text{ or } c>c_*^+.
		\end{cases}
\end{equation*}
\end{proposition}
\begin{proof}
The starting point is the large-domain formula
\eqref{eq:large-domain-limit-speed}. 
Changing variables from $\nu$ to $\mu=-\nu$ gives
\begin{equation}\label{C5-1}
		\lambda_\infty(c)
		=
		\min_{\nu\in\mathcal E} \{\Lambda_0(\nu)+c\nu\}
		=
		\min_{\mu\in\mathbb R}
		\{\Lambda_0(-\mu)-c\mu\}
		\in\mathbb R \text{ for } c\in\mathbb R.
	\end{equation}
	Let $c_1,c_2\in\mathbb R$ and $0<\theta<1$. For every $\nu\in\mathcal E$,
	\eqref{C5-1} gives
\begin{equation*}
\begin{aligned}
\Lambda_0(\nu)+(\theta c_1+(1-\theta)c_2)\nu
&=\theta\bigl(\Lambda_0(\nu)+c_1\nu\bigr)
+(1-\theta)\bigl(\Lambda_0(\nu)+c_2\nu\bigr)\\
&\ge \theta\lambda_\infty(c_1)+(1-\theta)\lambda_\infty(c_2).
\end{aligned}
\end{equation*}
Because the inequality holds for every $\nu\in\mathcal E$,
we can take the minimum on the left. This gives
\begin{equation*}
\lambda_\infty(\theta c_1+(1-\theta)c_2)
\ge\theta\lambda_\infty(c_1)
+(1-\theta)\lambda_\infty(c_2).
\end{equation*}
Thus $\lambda_\infty$ is concave. A finite concave function on $\mathbb R$ is continuous.

We first show that $c_*^-<c_*^+$, including the possibility
of infinite speeds. The following estimate also explains why
a finite extremum in the speed formulas is attained.
Fix $M>0$. From \eqref{C5-1},
\begin{equation*}
\Lambda_0(-\mu)\ge
\lambda_\infty(M)+M\mu,
\qquad
\Lambda_0(-\mu)\ge
\lambda_\infty(-M)-M\mu.
\end{equation*}
Choose the first inequality if $\mu\ge0$ and the second
if $\mu<0$. Then
\begin{equation*}
	\Lambda_0(-\mu)
	\geq M|\mu|
	+\min\{\lambda_\infty(M),\lambda_\infty(-M)\},
	\qquad \mu\in\mathbb R.
\end{equation*}
As $M$ can be chosen arbitrarily large, this inequality
shows that $\Lambda_0(-\mu)$ grows faster than any prescribed
multiple of $|\mu|$ as $|\mu|\to\infty$.
Near $\mu=0$, lower semicontinuity and
$\Lambda_0(0)>0$ imply
$\Lambda_0(-\mu)\ge\tfrac12\Lambda_0(0)>0$
for all sufficiently small $|\mu|$. Together, these facts give
\begin{equation*}
\frac{\Lambda_0(-\mu)}{|\mu|}
	\longrightarrow+\infty
	\text{ as }\mu\to0
	\text{ or as }|\mu|\to\infty.
\end{equation*}
	The above estimates, together with the lower
semicontinuity of $\Lambda_0$, show that the extrema
defining $c_*^\pm$ are attained whenever finite.
Moreover, $c_*^-\in\mathbb R\cup\{-\infty\}$ and
$c_*^+\in\mathbb R\cup\{+\infty\}$.
	If either speed is infinite, then $c_*^-<c_*^+$ follows
	immediately. Otherwise, choose $\mu_-<0<\mu_+$ such that
	\[
	\Lambda_0(-\mu_-)=c_*^-\mu_-,
	\qquad
	\Lambda_0(-\mu_+)=c_*^+\mu_+.
	\]
Since $-\mu_+<0<-\mu_-$, the strict convexity
of $\Lambda_0$ gives
	\[
	\begin{aligned}
		0<\Lambda_0(0)
		<
		\frac{\mu_+\Lambda_0(-\mu_-)
			-\mu_-\Lambda_0(-\mu_+)}
		{\mu_+-\mu_-}=
		\frac{(-\mu_-)\mu_+}{\mu_+-\mu_-}
		(c_*^+-c_*^-).
	\end{aligned}
	\]
	Thus $c_*^-<c_*^+$ also holds in this case.

Suppose $c<c_*^-$. By the definition of $c_*^-$,
there exists $\mu<0$ such that
	$\frac{\Lambda_0(-\mu)}{\mu}>c$.
	It follows that
	\[
	\lambda_\infty(c)
	\leq \Lambda_0(-\mu)-c\mu
	=
	\mu\left(
	\frac{\Lambda_0(-\mu)}{\mu}-c
	\right)<0.
	\]
	Similarly, if $c>c_*^+$, we may choose $\mu>0$ with
	$\Lambda_0(-\mu)/\mu<c$ and obtain $\lambda_\infty(c)<0$.

Let $c\in(c_*^-,c_*^+)\cap\mathbb R$.
By \eqref{C5-1}, there exists $\mu_c\in\mathbb R$
	such that
	\[
	\lambda_\infty(c)
	=\Lambda_0(-\mu_c)-c\mu_c.
	\]
	Suppose, to the contrary, that $\lambda_\infty(c)\leq0$.
	Since $\Lambda_0(0)>0$,  $\mu_c\neq0$.
	If $\mu_c<0$, then
	$c \leq\frac{\Lambda_0(-\mu_c)}{\mu_c}
	\leq c_*^-$,
	contrary to $c>c_*^-$. If $\mu_c>0$, then
	$c \geq\frac{\Lambda_0(-\mu_c)}{\mu_c}
	\geq c_*^+$,
	contrary to $c<c_*^+$.
	Therefore $\lambda_\infty(c)>0$ throughout
	$(c_*^-,c_*^+)\cap\mathbb R$.

At any finite endpoint, the preceding sign assertions give
positive values of $\lambda_\infty$ on one side and
negative values on the other. Its continuity then forces
\begin{equation*}
\lambda_\infty(c)=0
\enspace \text{for}\enspace c\in\{c_*^-,c_*^+\}\cap\mathbb R.
\end{equation*}
This completes the proof.
\end{proof}

\begin{proof}[Proof of Theorem~\ref{thm:spread-periodic-cooperative}]
The proof of (i) has two stages. We first obtain convergence
on every bounded interval in a moving frame with a fixed speed
$c$ strictly between $c_*^-$ and $c_*^+$. We then use these
moving-frame estimates to obtain convergence uniformly on the
entire expanding interval $[at,bt]$. The two exterior-decay
assertions follow from exponentially decaying supersolutions.

	Since $F(t,\mathbf0)=\mathbf0$, $F(t,\mathbf M)\leq\mathbf0$, and
	$\int_{\mathbb R}k_i(z,t)\,dz=1$ for every $i\in\mathbb S$ and
	$t\in\mathbb R$, comparison with the constant vectors $\mathbf0$
	and $\mathbf M$ gives
	\begin{equation*}
		\mathbf0\leq U(x,t)\leq\mathbf M,
		\qquad x\in\mathbb R,\quad t\geq0.
	\end{equation*}

\textbf{Proof of (i).}
\textbf{Step 1. Convergence on bounded intervals in a moving frame.}
Fix $c\in(c_*^-,c_*^+)\setminus\{0\}$. By
Proposition~\ref{prop-sign}, $\lambda_\infty(c)>0$.
For $c<0$, the bounded-interval results below are applied
after the spatial reflection $x\mapsto-x$, which places the
inflow condition at the right endpoint. In particular,
$\Lambda_0(0)>0$. The periodic ODE theory
used in Theorem~\ref{Kaiguan}(iv) shows that
$v'=F(t,v)$ has a unique positive
$T$-periodic solution $U^*$ in $[\mathbf0,\mathbf M]$. Moreover,
	it attracts every solution with initial value in
	$[\mathbf0,\mathbf M]\setminus\{\mathbf0\}$.
	Let $\overline U$ solve
	$\overline U'=F(t,\overline U)$,
	$\overline U(0)=\mathbf M$.
	Since each kernel has unit mass, $(x,t)\mapsto\overline U(t)$
	also solves the Cauchy equation. 
The comparison principle gives
	\begin{equation*}
		\mathbf0\leq U(x,t)\leq\overline U(t)\leq\mathbf M,
		\quad x\in\mathbb R,\, t\geq0.
	\end{equation*}
	As $|\overline U(t)-U^*(t)|\to0$, it follows that, for every
	$i\in\mathbb S$,
	\begin{equation}\label{C5-3}
		\limsup_{t\to\infty}\sup_{x\in\mathbb R}
		[U_i(x,t)-U_i^*(t)]\leq0.
	\end{equation}
	
The upper bound has therefore been established uniformly
in $x$. To obtain the matching lower bound, it is enough to
work on a fixed bounded interval in the moving coordinate
$x+ct$. We claim that, for every $R>0$,
\begin{equation}\label{C5-4}
\lim_{t\to\infty}\sup_{|x|\leq R}
|U(x+ct,t)-U^*(t)|=0.
\end{equation}
Because $\lambda_\infty(c)>0$, the large-domain spectral
limit gives $l_0>0$ such that the principal eigenvalue
$\lambda_l(c)$ on $(-l,l)$ is positive for
$l>l_0$.
By Theorem~\ref{Kaiguan}, the time-periodic problem \eqref{draft-P} on $(-l,l)$ admits a unique positive $T$-periodic solution, denoted by $V_{l,c}^*$.
By Theorem~\ref{Kaiguan}(iv), for fixed
$R>0$ and $\varepsilon>0$, we can choose $l>\max \{R, l_0 \}$ such that
\begin{equation*}
\sup_{|x|\leq R,0\leq t\leq T}
|V_{l,c}^{*}(x,t)-U^*(t)|<\varepsilon.
\end{equation*}
By periodicity, this estimate holds for all $t\in\mathbb R$.

The solution from the nonzero compactly supported initial
datum is strictly positive at every spatial point after any
positive time, by the positivity of the kernels near zero
and the cooperativity and irreducibility of the system. In particular,
$$\min_{x\in[-l,l],\,i\in\mathbb S}U_i(x+cT,T)>0.$$
Let $\varphi_{l,c}$ be a positive $T$-periodic principal eigenfunction associated with $\lambda_l(c)$.
The eigenfunction $\varphi_{l,c}$ is continuous and
positive in the interior, and vanishes at the inflow boundary.
Since $U(\cdot+cT,T)$ has a strictly positive minimum on the
compact interval $[-l,l]$, a sufficiently small constant
$\delta_l>0$ gives
	\begin{equation}\label{C5-5}
		V_{l,0}^c:=\delta_l\varphi_{l,c}(\cdot,0),
		\qquad
		\mathbf0\not\equiv V_{l,0}^c(x)\leq U(x+cT,T) \le \mathbf{M}
		\enspace\text{for}\enspace x\in[-l,l].
	\end{equation}
	Let $W_l^c$ solve
	\begin{equation*}
		\left\{
		\begin{aligned}
			&\partial_tW_{l,i}^c(x,t)
			=
			d_i(t)\left[
			\int_{-l}^{l}k_i(x-y,t)W_{l,i}^c(y,t)\,dy-W_{l,i}^c(x,t) \right]\\
			&\hspace{2.5cm}+c\partial_xW_{l,i}^c(x,t)+f_i(t,W_l^c(x,t)),
			&&x\in(-l,l),\, t>0,\, i\in\mathbb{S},\\
			&W_l^c(l,t)=\mathbf0,
			&&t>0,\quad c>0, \\
			&W_l^c(-l,t)=\mathbf0,
			&&t>0,\quad c<0,\\
			&W_l^c(x,0)=V_{l,0}^c(x),
			&&x\in[-l,l].
		\end{aligned}
		\right.
	\end{equation*}
	By Theorem \ref{Kaiguan}, 
   $\lim_{t\to\infty}
		\|W_l^c(\cdot,t)-V_{l,c}^{*}(\cdot,t)\|_\infty=0.$

We now compare the bounded-interval solution with the
whole-line Cauchy solution. This requires translating the
bounded interval in physical space while retaining its
zero inflow boundary. Set
	\begin{equation*}
		I_l^c(t)=(-l+ct,l+ct),\quad
		\widehat U(y,t)=U(y+cT,t+T),\quad
		\widetilde W_l^c(y,t)=W_l^c(y-ct,t).
	\end{equation*}
For $y\in I_l^c(t)$, the translated solution
$\widehat U(y,t)=U(y+cT,t+T)$ solves the original
whole-line equation, because all coefficients are $T$-periodic.
Since $\widehat U$ satisfies the equation on
$\mathbb R$, restricting the nonlocal integral
to $I_l^c(t)$ leaves a nonnegative contribution
from the complement. Thus, for each component,
	\begin{equation*}
		\begin{aligned}
			&\partial_t\widehat U_i(y,t)
			-d_i(t)\left[
			\int_{I_l^c(t)}k_i(y-z,t)\widehat U_i(z,t)\,dz
			-\widehat U_i(y,t)
			\right]-f_i(t,\widehat U(y,t))\\
			=&
			d_i(t)\int_{\mathbb R\setminus I_l^c(t)}
			k_i(y-z,t)\widehat U_i(z,t)\,dz\geq0.
		\end{aligned}
	\end{equation*}
Together with
the initial comparison \eqref{C5-5} and the following inflow
boundary inequality,
	\begin{equation*}
		\begin{cases}
			\widetilde W_l^c(l+ct,t)=\mathbf0
			\leq\widehat U(l+ct,t),&c>0,\\
			\widetilde W_l^c(-l+ct,t)=\mathbf0
			\leq\widehat U(-l+ct,t),&c<0,
		\end{cases}
	\end{equation*}
the comparison principle on the moving interval yields
$\widehat U\ge\widetilde W_l^c$. Returning to the
moving coordinate gives
	\begin{equation*}
		U\left(x+c(t+T),t+T\right)\geq W_l^c(x,t),
		\quad x\in[-l,l],\, t\geq0.
	\end{equation*}
	Put $s=t+T$. By the periodicity of $V_{l,c}^*$ and the choice
	of $l$, for every $i\in\mathbb S$ and $s\geq T$,
	\begin{equation*}
		\begin{aligned}
			\inf_{|x|\leq R}[U_i(x+cs,s)-U_i^*(s)]
			&\geq \inf_{|x|\leq R}[W_{l,i}^c(x,s-T)-U_i^*(s)]\\
			&\geq
			-\|W_l^c(\cdot,s-T)-V_{l,c}^*(\cdot,s-T)\|_\infty
			-\varepsilon.
		\end{aligned}
	\end{equation*}
	Letting first $s\to\infty$ and then $\varepsilon \to 0$, we obtain
	\begin{equation*}
		\liminf_{s\to\infty}\inf_{|x|\leq R}
		[U_i(x+cs,s)-U_i^*(s)]\geq0,
		\quad i\in\mathbb S.
	\end{equation*}
	Together with \eqref{C5-3}, this gives \eqref{C5-4}.

\textbf{Step 2. Uniform convergence on an expanding interior interval.}
We first consider $[a,b]\subset(c_*^-,c_*^+)$ with $a<b$ and $a,b\neq0$.
For $0<\rho<\Lambda_0(0)$, consider the perturbed periodic ODE
	\begin{equation}\label{eq:spread-perturbed-ODE}
		z'=F(t,z)-\rho z.
	\end{equation}
Subtracting $\rho z$ only changes the diagonal of the
Jacobian, so cooperativity and irreducibility are preserved.
Strict subhomogeneity is also preserved, since
	\begin{equation*}
		[F(t,\theta z)-\rho\theta z]
		-\theta[F(t,z)-\rho z]
		=F(t,\theta z)-\theta F(t,z).
	\end{equation*}
	Moreover, $F(t,\mathbf M)-\rho\mathbf M\leq\mathbf0$.
	The linearization at zero is $A(t)-\rho I$, so its principal
	Floquet exponent is $\Lambda_0(0)-\rho>0$. Thus
	\eqref{eq:spread-perturbed-ODE} has a unique positive
	$T$-periodic solution $U_\rho^*$, which attracts every solution
	with initial value in $[\mathbf0,\mathbf M]\setminus\{\mathbf0\}$.
	
	Since $F(t,z)-\rho z$ is nonincreasing in $\rho$ for
	$z\geq\mathbf0$, the comparison principle gives
	\begin{equation*}
		U_{\rho_2}^*(t)\leq U_{\rho_1}^*(t)\leq U^*(t)
		\enspace\text{for}\enspace 0<\rho_1<\rho_2<\Lambda_0(0)
		\enspace\text{and}\enspace 
		t\in\mathbb R.
	\end{equation*}
In fact, $U^*_\rho(t)\ll U^*(t)$ for every $\rho\in(0,\Lambda_0(0))$ and $t\in\mathbb R$.  If
some component of $r:=U^*-U_\rho^*$ vanished at time $t_0$,
that component would attain a local minimum equal to zero,
so its time derivative would vanish. At this contact time,
cooperativity of $F$ implies that increasing the other
components cannot decrease $f_i$. Therefore
	\begin{equation*}
		0=r_i'(t_0)
		=f_i(t_0,U^*(t_0))-f_i(t_0,U_\rho^*(t_0))
		+\rho U_{\rho,i}^*(t_0)
		\geq\rho U_{\rho,i}^*(t_0)>0,
	\end{equation*}
	a contradiction. Consequently,
	\begin{equation}\label{eq:spread-perturbed-gap}
		\gamma_\rho:=
		\min_{\substack{0\leq t\leq T\\i\in\mathbb S}}
		[U_i^*(t)-U_{\rho,i}^*(t)]>0.
	\end{equation}

We next explain why the perturbed periodic state tends to
the unperturbed one as $\rho\to 0$. Fix
$\rho_0\in(0,\Lambda_0(0))$. For $0<\rho\le\rho_0$,
	$\mathbf0\ll U_{\rho_0}^*(t)\leq U_\rho^*(t)\leq U^*(t)$.
	Equation \eqref{eq:spread-perturbed-ODE} gives a uniform bound on
	$(U_\rho^*)'$ in $[0,T]$. Hence every sequence $\rho_n\to 0$
	has a subsequence along which $U_{\rho_n}^*$ converges uniformly
	in $[0,T]$. Passing to the limit in
	\begin{equation*}
		U_{\rho_n}^*(t)=U_{\rho_n}^*(0)
		+\int_0^t
		[F(s,U_{\rho_n}^*(s))-\rho_nU_{\rho_n}^*(s)]\,ds,
	\end{equation*}
	we see that the limit is a $T$-periodic solution of $z'=F(t,z)$
	bounded below by $U_{\rho_0}^*$. By uniqueness, it equals $U^*$.
	Therefore,
	\begin{equation*}
		\lim_{\rho\downarrow0}
		\sup_{0\leq t\leq T}|U_\rho^*(t)-U^*(t)|=0.
	\end{equation*}
	
	Fix $\varepsilon>0$ and choose $\rho>0$ sufficiently small that
	$\sup_{t\in\mathbb R}|U_\rho^*(t)-U^*(t)|<\varepsilon$.
To control the mass lost outside a finite spatial interval,
use the unit-mass normalization of the kernels. For each $i$,
the continuous nonnegative functions
$t\mapsto1-\int_{-n}^n k_i(z,t)\,dz$ decrease to zero
pointwise as $n\to\infty$. Since $[0,T]$ is compact and the
limit is continuous, Dini's theorem makes this convergence
uniform in $t$.
Hence we may choose $R>0$ such that
	\begin{equation}\label{eq:spread-tail-small}
		\max_{i\in\mathbb S}\sup_{0\leq t\leq T}
		d_i(t)\int_{|z|>R}k_i(z,t)\,dz<\rho.
	\end{equation}
	Set
	\begin{equation*}
		\Omega_t=[at,bt],\qquad
		\Omega_t^R=[at-R,bt+R].
	\end{equation*}
Step~1 applied with the two speeds $a$ and $b$ supplies
lower bounds in fixed-width neighborhoods of $at$ and $bt$.
Indeed, the sets
$[at-R,at]$ and $[bt,bt+R]$ correspond respectively to
$x-at\in[-R,0]$ and $x-bt\in[0,R]$.
Using the positive gap
	\eqref{eq:spread-perturbed-gap}, we can choose an integer $N$
	large enough that, with $t_1=NT$,
	\begin{equation}\label{eq:spread-collar-lower}
		U_i(x,t)\geq U_{\rho,i}^*(t)+\frac{\gamma_\rho}{2},
		\quad
		t\geq t_1,\,
		x\in[at-R,at]\cup[bt,bt+R],\, i\in\mathbb S.
	\end{equation}
At time $t_1$, each component of $U$ is strictly positive
on the compact interval $[at_1,bt_1]$. Because
$U_\rho^*(t_1)\gg\mathbf0$, a sufficiently small
$\sigma\in(0,1)$ satisfies
	\begin{equation*}
		\sigma U_\rho^*(t_1)\ll U(x,t_1),
		\quad x\in\Omega_{t_1}.
	\end{equation*}
	Let $v$ solve \eqref{eq:spread-perturbed-ODE} with
	$v(t_1)=\sigma U_\rho^*(t_1)$. By comparison and the attractivity
	of $U_\rho^*$,
	\begin{equation}\label{eq:spread-ODE-lower-properties}
		\mathbf0\ll v(t)\leq U_\rho^*(t),
		\, t\geq t_1,\qquad
		|v(t)-U_\rho^*(t)|\to0
		\text{ as }t\to\infty.
	\end{equation}

We now compare $U$ with the spatially constant function
$v(t)$ on the moving interval $\Omega_t=[at,bt]$.
The values of $U$ in the surrounding strips provide a
nonnegative contribution to the nonlocal integral. For
$x\in\Omega_t$ and $t\ge t_1$, denote that contribution by
	\begin{equation*}
		H_i(x,t)=d_i(t)\int_{\Omega_t^R\setminus\Omega_t}
		k_i(x-y,t)U_i(y,t)\,dy.
	\end{equation*}
	The Cauchy equation and the nonnegativity of $U$ give
	\begin{equation*}
		\partial_tU_i(x,t)\geq
		d_i(t)\left[
		\int_{\Omega_t}k_i(x-y,t)U_i(y,t)\,dy-U_i(x,t)
		\right]+f_i(t,U(x,t))+H_i(x,t).
	\end{equation*}
	By \eqref{eq:spread-collar-lower} and
	\eqref{eq:spread-ODE-lower-properties}, $U_i(y,t)\geq v_i(t)$
	in $\Omega_t^R\setminus\Omega_t$. 
For $x\in\Omega_t$, all points $y$ with
$|x-y|\le R$ belong to $\Omega_t^R$. Consequently,
$\int_{\Omega_t^R}k_i(x-y,t)\,dy
\ge1-\int_{|z|>R}k_i(z,t)\,dz$.
Using the tail estimate \eqref{eq:spread-tail-small}
gives for $x\in\Omega_t$,
	\begin{equation*}
		\begin{aligned}
			&d_i(t)\left[
			\int_{\Omega_t}k_i(x-y,t)v_i(t)\,dy-v_i(t)
			\right]+f_i(t,v(t))+H_i(x,t)\\
			\geq&
			d_i(t)v_i(t)\left[
			\int_{\Omega_t^R}k_i(x-y,t)\,dy-1
			\right]+f_i(t,v(t))\\
			> &
			-\rho v_i(t)+f_i(t,v(t))=v_i'(t).
		\end{aligned}
	\end{equation*}
The initial comparison and the boundary comparison in
(5.8)--(5.9) are strict. Any first contact would therefore
occur at an interior point $(x_0,t_0)$ in some component
$i\in\mathbb S$. At this point, cooperativity, nonnegativity
of the nonlocal integral difference, and the strict
inequality above give
$0\geq \partial_t(U_i-v_i)(x_0,t_0)>0$,
a contradiction. Thus $U(x,t)\geq v(t)$ for
$x\in\Omega_t$ and $t\geq t_1$.
Therefore, for every $i\in\mathbb S$,
	\begin{equation*}
		\inf_{x\in[at,bt]}[U_i(x,t)-U_i^*(t)]
		\geq v_i(t)-U_i^*(t)
		\geq-|v(t)-U_\rho^*(t)|-\varepsilon.
	\end{equation*}
	Letting $t\to\infty$, we obtain
	\begin{equation*}
		\liminf_{t\to\infty}\inf_{x\in[at,bt]}
		[U_i(x,t)-U_i^*(t)]\geq-\varepsilon.
	\end{equation*}
	Letting $\varepsilon\to 0$ and using \eqref{C5-3}, we conclude
	that
	\begin{equation*}
		\lim_{t\to\infty}\sup_{x\in[at,bt]}|U(x,t)-U^*(t)|=0.
	\end{equation*}
	
	Now let $[a,b]\subset(c_*^-,c_*^+)$ be arbitrary, with $a\leq b$.
	Choose $\widetilde a,\widetilde b\in(c_*^-,c_*^+)\setminus\{0\}$
	such that $\widetilde a<a\leq b<\widetilde b$. By the case just
	proved,
	\begin{equation*}
		\sup_{x\in[at,bt]}|U(x,t)-U^*(t)|
		\leq
		\sup_{x\in[\widetilde a t,\widetilde b t]}|U(x,t)-U^*(t)|
		\longrightarrow0
		\text{ as }t\to\infty.
	\end{equation*}
	This proves (i).

\textbf{Proof of (ii) and (iii).}
Suppose first that $c_*^+<+\infty$ and $b_2>c_*^+$.
	By the definition of $c_*^+$, we can choose $\mu>0$ with
	$-\mu\in\mathcal E$ such that
	\begin{equation*}
		s_\mu:=\frac{\Lambda_0(-\mu)}{\mu}<b_2.
	\end{equation*}
	Let $q^\mu$ be a positive $T$-periodic Floquet eigenfunction
	associated with $\Lambda_0(-\mu)$. 
	Since $m_i(-\mu,\cdot)\in L^1(0,T)$ for every $i\in\mathbb S$,
	we have $q^\mu\in W^{1,1}((0,T);\mathbb R^m)$ and
	\begin{equation*}
		(q^\mu)'(t)
		=\bigl[A(t)+D(t)\mathsf M(-\mu,t)\bigr]q^\mu(t)
		-\Lambda_0(-\mu)q^\mu(t)
		\text{ for a.e. }t\in(0,T).
	\end{equation*}
	For $C>0$, set
	\begin{equation*}
		\overline U^\mu(x,t)=Ce^{-\mu(x-s_\mu t)}q^\mu(t).
	\end{equation*}
To verify that $\overline U^\mu$ solves the linearized
equation, differentiate its time factor and change variables
$z=x-y$ in the nonlocal integral. This gives
\[
\partial_t\overline U_i^\mu
=Ce^{-\mu(x-s_\mu t)}
\bigl[\mu s_\mu q_i^\mu+(q_i^\mu)'\bigr],
\quad
\int_\mathbb R k_i(x-y,t)\overline U_i^\mu(y,t)\,dy
=m_i(-\mu,t)\overline U_i^\mu(x,t).
\]
Substituting these identities and the Floquet equation,
we obtain for every $x\in\mathbb R$ and almost every $t>0$,
	\begin{equation*}
		\begin{aligned}
			&\partial_t\overline U_i^\mu(x,t) -d_i(t)\left[
			\int_{\mathbb R}k_i(x-y,t)\overline U_i^\mu(y,t)\,dy
			-\overline U_i^\mu(x,t)\right]
			-\sum_{j=1}^m a_{ij}(t)\overline U_j^\mu(x,t)\\
			=&
			[\mu s_\mu-\Lambda_0(-\mu)]\overline U_i^\mu(x,t)=0.
		\end{aligned}
	\end{equation*}
	Since $U_0$ has compact support and $q^\mu(0)\gg\mathbf0$,
	we can choose $C>0$ such that
	$U_0(x)\leq\overline U^\mu(x,0)$ for all $x\in\mathbb R$.
	
	For $\mathbf0 < z\leq\mathbf M$, subhomogeneity gives
	$F(t,z)\leq\theta^{-1}F(t,\theta z)$ for $0<\theta<1$.
	Letting $\theta\to 0$ yields $F(t,z)\leq A(t)z$. Hence,
	\begin{equation*}
		\partial_tU_i(x,t)\leq
		d_i(t)\left[
		\int_{\mathbb R}k_i(x-y,t)U_i(y,t)\,dy-U_i(x,t)
		\right]+\sum_{j=1}^m a_{ij}(t)U_j(x,t).
	\end{equation*}
	Set $w=U-\overline U^\mu$. For any $\tau>0$, define
\[
h(t):=\sup_{x\in\mathbb R}\max_{i\in\mathbb S}
(w_i(x,t))_+,\qquad 0\leq t\leq\tau.
\]
Then $0\leq h(t)\leq\max_{i\in\mathbb S}M_i$.
For each fixed $x$, $w_i(x,\cdot)$ is absolutely continuous
on $[0,\tau]$, and its convolution term is integrable in time
because $m_i(-\mu,\cdot)\in L^1(0,T)$.
Writing $b_i(t)=a_{ii}(t)-d_i(t)$ and using
$w_i(x,0)\leq0$, we integrate the differential inequality
for $w_i$, which holds almost everywhere in time, to obtain
\begin{equation*}
\begin{aligned}
w_i(x,t)
&\leq
\int_0^t e^{\int_s^t b_i(r)\,dr}
\left[
d_i(s)\int_{\mathbb R}k_i(x-y,s)w_i(y,s)\,dy
+\sum_{j\ne i}a_{ij}(s)w_j(x,s)
\right]ds\\
&\leq
\int_0^t e^{\int_s^t b_i(r)\,dr}
\left[d_i(s)+\sum_{j\ne i}a_{ij}(s)\right]h(s)\,ds\\
&\leq C_\tau\int_0^t h(s)\,ds.
\end{aligned}
\end{equation*}
Here we used cooperativity and the unit mass of the kernels,
and $C_\tau>0$ is independent of $x$, $t$, and $i$.
Taking positive parts and then the supremum gives
$h(t)\leq C_\tau\int_0^t h(s)\,ds$. 
Gronwall's inequality yields $h\equiv0$ on $[0,\tau]$. As $\tau>0$ is
arbitrary, we obtain
\begin{equation*}
0\leq U(x,t)\leq Ce^{-\mu(x-s_\mu t)}q^\mu(t),
\quad x\in\mathbb R,\, t\geq0.
\end{equation*}
Writing $Q_\mu=\max_{0\leq t\leq T}|q^\mu(t)|$, we obtain
\begin{equation*}
    \sup_{x\geq b_2t}|U(x,t)|
    \leq C Q_\mu e^{-\mu(b_2-s_\mu)t}\to0
   \text{ as }t\to\infty.
\end{equation*}

If $c_*^->-\infty$ and $a_2<c_*^-$, the proof is similar, and we omit the details here.
\end{proof}

\appendix
\section{}

\begin{proposition}\label{cone-interior}
Assume (H1)–(H3), and let \(K=\{\phi\in\mathbb V_{(l_1,l_2)}:\phi\ge0\}\). Then K has nonempty interior relative to the C¹-topology of \(\mathbb V_{(l_1,l_2)}\), and 
\begin{equation*}
\operatorname{int}K=\left\{\phi\in\mathbb V_{(l_1,l_2)}:
		\phi_i(x)>0 \text{ for all } x\in[l_1,l_2),\
		i\in\mathbb S\right\}.
\end{equation*}
\end{proposition}

\begin{proof}
Write $\delta(x):=l_2-x$ and $\ell:=l_2-l_1$. Denote the set on the right-hand side by $K_+$.

\textbf{Step 1. In this step, we construct a function in $K_+$.}
For each $i\in\mathbb S$ and $\alpha\in\mathbb R$, set
\begin{equation*}
		p_i^\alpha(x)=\delta(x)e^{\alpha\delta(x)},\qquad
		G_i(\alpha)=\frac{d_i(l_2,0)}{c(0)}
		\int_{l_1}^{l_2}k_i(l_2-y,0)\delta(y)
		e^{\alpha\delta(y)}\,dy.
\end{equation*}
Clearly, $p_i^\alpha(l_2)=0$ and $(p_i^\alpha)'(l_2)=-1$. Consequently,
\begin{equation*}
		c(0)(p_i^\alpha)'(l_2)+d_i(l_2,0)\int_{l_1}^{l_2}k_i(l_2-y,0)p_i^\alpha(y)\,dy
		=c(0)\bigl(G_i(\alpha)-1\bigr).
\end{equation*}
Thus $p_i^\alpha$ satisfies the boundary compatibility condition if and only if $G_i(\alpha)=1$.
	
We now show that such an $\alpha$ exists. Since $0\leq\delta(y)\leq\ell$ and $k_i(l_2-\cdot,0)\in L^1(l_1,l_2)$, $G_i$ is continuous on $\mathbb R$. Moreover, by dominated convergence, $\lim_{\alpha\to-\infty}G_i(\alpha)=0$. Since $\delta(y)>0$ for $y\in(l_1,l_2)$ and $\int_{l_1}^{l_2}k_i(l_2-y,0)\,dy>0$, the monotone convergence theorem gives 	$\lim_{\alpha\to+\infty}G_i(\alpha)=+\infty$. Hence, by the intermediate value theorem, there exists 	$\alpha_i\in\mathbb R$ such that $G_i(\alpha_i)=1$.
	
Taking $p_i=p_i^{\alpha_i}$ and $p:=(p_1,\ldots,p_m)$, we obtain $p\in\mathbb V_{(l_1,l_2)}$ and $p_i(x)>0$ for 	every $x\in[l_1,l_2)$ and $i\in\mathbb S$. Thus $p\in K_+$, proving that $K_+$ is nonempty.
	
\textbf{Step 2. In this step, we prove $K_+\subset\operatorname{int}K$.}
Fix $\phi\in K_+$. Since $\phi_i(l_2)=0$, the boundary compatibility condition gives 
\begin{equation*}
		\lim_{x\to l_2}\frac{\phi_i(x)}{\delta(x)}
		=-\phi_i'(l_2)
		=\frac{d_i(l_2,0)}{c(0)}
		\int_{l_1}^{l_2}k_i(l_2-y,0)\phi_i(y)\,dy>0,
		\quad i\in\mathbb S.
\end{equation*}
Thus each $\phi_i/\delta$ extends to a strictly positive continuous function on $[l_1,l_2]$. Since $\mathbb S$ is finite, there exists $\eta_\phi>0$ such that
\begin{equation*}
		\phi_i(x)\geq\eta_\phi\delta(x),
		\quad x\in[l_1,l_2],\, i\in\mathbb S.
\end{equation*}
	
Let $h\in\mathbb V_{(l_1,l_2)}$ satisfy $\|h\|_{C^1}<\eta_\phi/2$. Since $h_i(l_2)=0$,
\begin{equation*}
		|h_i(x)|
		=\left|\int_x^{l_2}h_i'(y)\,dy\right|
		\leq\delta(x)\|h\|_{C^1}.
\end{equation*}
It follows that
\begin{equation*}
(\phi_i+h_i)(x) \geq\bigl(\eta_\phi-\|h\|_{C^1}\bigr)\delta(x)\geq\frac{\eta_\phi}{2}\delta(x)\geq0.
\end{equation*}
Hence $\phi+h\in K$, proving that $\phi\in\operatorname{int}K$. Since $K_+$ is nonempty by Step~1, $\operatorname{int}K$ is also nonempty.
	
\textbf{Step 3. In this step, we prove $\operatorname{int}K\subset K_+$.}
Let $\phi\in\operatorname{int}K$. By definition, there exists $r_\phi>0$ such that
\begin{equation*}
		\phi+h\in K
		\quad\text{whenever }h\in\mathbb V_{(l_1,l_2)}
		\text{ and }\|h\|_{C^1}<r_\phi.
\end{equation*}
Choose $\varepsilon>0$ with $\varepsilon\|p\|_{C^1}<r_\phi$. Since $-\varepsilon p\in\mathbb V_{(l_1,l_2)}$, we obtain
\begin{equation*}
		\phi-\varepsilon p\in K,\qquad
		\phi_i(x)\geq\varepsilon p_i(x)>0
		\quad\text{for }x\in[l_1,l_2),\ i\in\mathbb S.
\end{equation*}
Thus $\phi\in K_+$, completing the proof.
\end{proof}

\section*{Acknowledgments}
\noindent

Research of Li was partially supported by the National Natural Science Foundation of China (12531008; 12271226).

%
%
%
%


\end{document}